\documentclass[a4paper,12pt,leqno]{article}
\usepackage{amsmath,amsfonts,amsthm,amsrefs}
\usepackage[OT4]{fontenc}

\long\def\symbolfootnote[#1]#2{\begingroup%
\def\thefootnote{\fnsymbol{footnote}}\footnote[#1]{#2}\endgroup}

\newtheorem{theorem}{Theorem}[section]
\newtheorem{lemma}[theorem]{Lemma}
\newtheorem{proposition}[theorem]{Proposition}
\newtheorem{claim}[theorem]{Claim}
\newtheorem{cor}[theorem]{Corollary}

\newtheorem{lem}[theorem]{Lemma}

\newtheorem{tw}[theorem]{Theorem}
\newtheorem{wtw}{Theorem}

\newtheorem{prop}[theorem]{Proposition}
\newtheorem*{cor*}{Corollary}

\theoremstyle{definition}
\newtheorem{defin}[theorem]{Definition}
\newtheorem{rem}[theorem]{Remark}
\renewcommand{\proof}{\medskip\par\noindent\textbf{Proof.} \ignorespaces}
\newcommand{\dow}{\medskip\par\noindent\textbf{Proof.} \ignorespaces}

\newcommand{\Span}{\mathrm{span}}
\newcommand{\conv}{\mathrm{conv}}
\newcommand{\diam}{\mathrm{diam}}
\renewcommand{\qed}{\quad\hskip0pt\null\hfill$\square$\par}
\newcommand{\kon}{\quad\hskip0pt\null\hfill$\square$\par}

\newcommand {\mr}{\mathrm}
\newcommand {\lk}{\left\{}
\newcommand {\rk}{\right\}}

\usepackage[dvips]{graphicx}

\newcommand {\bdo}{\partial_{O}}
\newcommand {\bdoo}{\partial_{O'}}

\newcommand {\cx}{\overline{X}}

\begin{document}

\begin{center}
\large\bfseries Boundaries of systolic groups
\end{center}

\begin{center}\bf
Damian Osajda \symbolfootnote[1]{Partially supported by MNiSW
grant N201 012 32/0718. This research was supported by a Marie Curie European Reintegration Grant within the 6th European Community Framework Programme.}

\it Instytut Matematyczny, Uniwersytet Wroc\l awski,
pl. Grunwaldzki 2/4,
50--384 Wroc{\l}aw, Poland

\rm dosaj@math.uni.wroc.pl
\end{center}

\begin{center}\bf
Piotr Przytycki \symbolfootnote[2]{Partially supported by MNiSW
grant N201 003 32/0070, MNiSW
grant N201 012 32/0718, and the Foundation for Polish Science.}

\it Institute of Mathematics, Polish Academy of Sciences,
 \'Sniadeckich 8, 00-956 Warsaw, Poland

\rm pprzytyc@mimuw.edu.pl
\end{center}

\begin{abstract}
\noindent For all systolic groups we construct boundaries which are $EZ$--structures. This implies the Novikov conjecture for torsion--free systolic groups. The boundary is constructed via a system of distinguished geodesics in a systolic complex, which we prove to have coarsely similar properties to geodesics in $CAT(0)$ spaces.
\end{abstract}



\tableofcontents


\section{Introduction}
\label{introduction}

There are many notions of boundaries of groups
used for various purposes. In this paper we focus on
the notions of $Z$--structure and $EZ$--structure
introduced by Bestvina \cite{Be} and studied
e.g. in \cite{D}, \cite{FL}.
Our main result is the following.
\begin{wtw}[Theorem \ref{main}]
\label{1}
Let a group $G$ act geometrically by simplicial automorphisms
on a systolic complex $X$.
Then there exists a compactification $\cx=X\cup \partial X$ of $X$
satisfying the
following:
\begin{enumerate}
 \item $\cx$ is a Euclidean retract (ER),
 \item $\partial X$ is a $Z$--set in $\cx$,
 \item for every compact set $K\subset X$,
$(gK)_{g\in G}$ is a null sequence,
 \item the action of $G$ on $X$ extends to an action,
by homeomorphisms, of $G$ on $\cx$.
\end{enumerate}
\end{wtw}

A group $G$ as in Theorem \ref{1} is called a \emph{systolic group}.
It is a group acting \emph{geometrically} (i.e. cocompactly and properly discontinuously) by simplicial automorphisms
on a \emph{systolic complex}---contractible simplicial complex
satisfying some local combinatorial
conditions. Systolic complexes were introduced by Januszkiewicz--\'Swi\k{a}tkowski
\cite{JS} and, independently, by Haglund \cite{H} and by Chepoi \cite{C} (in Section \ref{Systolic complexes} we give some background on them).
Systolic complexes (groups) have many properties
of non--positively curved spaces (groups). There are systolic complexes
that are not $CAT(0)$ when equipped with the path metric in which
every simplex is isometric to the standard Euclidean simplex.
On the other hand, there are systolic groups that are not hyperbolic, e.g. $\mathbb{Z}^2$. Summarizing, systolic setting does not reduce to the $CAT(0)$ or to the hyperbolic one --- it turns out that systolic groups form a large family: allow various combinatorial constructions \cite{H},\cite{JS},\cite{JS+} and provide the discipline with new range of examples frequently with unexpected properties \cite{JS2},\cite{JS+}. We also believe that eventually both systolic complexes and $CAT(0)$ cubical ones will be placed among a wider family of combinatorially non--positively curved contractible cell complexes.
\medskip

Here we give the other definitions that appear in the statement of Theorem \ref{1}.
A compact
space is a {\it {Euclidean retract}} (or ER) if it
can be embedded in some Euclidean space as its retract.
A closed subset $Z$ of a Euclidean retract $Y$ is called a $Z$--{\it {set}}
if for every open set $U\subset Y$, the inclusion
$U\setminus Z\hookrightarrow U$ is a homotopy equivalence.
A sequence $(K_i)_{i=1}^{\infty}$ of subsets of a topological
space $Y$ is called a {\it {null sequence}} if for every open cover
${\cal U}=\lk U_i \rk _{i\in I}$ of $Y$ all but finitely many $K_i$
 are
\emph{$\cal U$--small}, i.e. for all but finitely many $j$ there
exist $i(j)$ such that $K_j\subset U_{i(j)}$.
\medskip

Conditions 1, 2 and 3 of Theorem \ref{1}
mean (following \cite{Be}, where only free actions are considered, and \cite{D})
that any systolic group $G$ admits a
\emph{$Z$--structure} $(\cx,\partial X)$.
The notion
of an \emph{$EZ$--structure}, i.e. a $Z$--structure with additional property
4 was explored by Farrell--Lafont \cite{FL} (in the case
of a free action).

Bestvina \cite{Be} showed that some local homological invariants
of the boundary $\partial X$ are related to cohomological
invariants of the group. In particular, the dimension of the boundary
is an invariant of the group i.e. it does not depend
on the $Z$--structure we choose. This was generalized by Dranishnikov
\cite{D} to the case of geometric actions. It should be emphasized
that the homeomorphism type of the boundary is not a group
invariant (but the shape is an invariant \cite{Be}),
as the Croke--Kleiner examples of
visual boundaries of some $CAT(0)$ spaces exhibit \cite{CrK}.
Carlsson--Pedersen \cite{CP} and Farrell--Lafont \cite{FL}
proved that existence of an $EZ$--structure on a torsion--free
group $G$ implies that the Novikov conjecture
is true for $G$. Thus, by Theorem \ref{1}, we get the following.

\begin{cor*}
Torsion--free systolic groups satisfy the Novikov conjecture.
\end{cor*}

There are only few classes of groups for which a $Z$--structure
$(\cx,\partial X)$ has been found (and even fewer for which an
$EZ$--structure is known).
The most important examples are: hyperbolic groups \cite{BeMe}
---with $X$ being the Rips complex and $\partial X$
being the Gromov boundary of $G$; $CAT(0)$ groups
---with $X$ a $CAT(0)$ space and $\partial X$ the visual boundary
of $X$; relatively hyperbolic groups whose parabolic subgroups
admit a $Z$--structure \cite{Da}.
Bestvina \cite{Be} asked whether every group $G$ with finite
$K(G,1)$ has a $Z$--structure.
\medskip

The question whether for
every systolic group there exists an $EZ$--structure was posed by Januszkiewicz and \'Swi\k{a}tkowski in 2004. Theorem \ref{1} answers
affirmatively this question.

We hope that, similarly to the hyperbolic and $CAT(0)$ cases,
our boundaries will be also useful for purposes other
than the ones mentioned above. In particular we think that
splittings of systolic groups can be recognized through the topology
of the boundary, as in e.g. \cite{B}, \cite{PS}.
Studying more refined structures on the boundary could help
in obtaining rigidity results for some systolic groups.
\medskip

The essential point of our construction is the choice of the system of \emph{good geodesics} (derived from the system of \emph{Euclidean geodesics}, the distinction being not important at this moment), which is coarsely closed under taking subsegments (Theorem B below), and which satisfies coarsely a weak form of $CAT(0)$ condition (Theorem C below).

Recall that Januszkiewicz and \'Swi\k{a}tkowski \cite{JS} considered a system of directed geodesics in a systolic complex (c.f. Definition \ref{directed geodesic}). One may try to define the boundary of a systolic complex by taking the inverse limit of the following system. Consider the sequence of combinatorial spheres around a fixed vertex $O$ and projections from larger to smaller spheres along the directed geodesics terminating at $O$. Unfortunately, the inverse limit of this system does not have, in general, property 3 from Theorem \ref{1}. Property 3 fails, for example, already for the flat systolic plane (c.f. Definition \ref{def_of_disc}).

Hence, instead of using directed geodesics, we introduce \emph{Euclidean geodesics}, which
behave like $CAT(0)$ geodesics with respect to the flat subcomplexes a systolic complex.
To define the Euclidean geodesic between two vertices, say $s,t$, in a systolic complex, we
consider the loop obtained by concatenating the two directed geodesics joining $s$ to $t$ and $t$ to $s$. Then we span a minimal surface $S$ on this loop. (Here Elsner's minimal surfaces theory \cite{E} comes in handy. To obtain some uniqueness properties on $S$ we complement Elsner's theory with our results on \emph{layers}, which span the union of all 1--skeleton geodesics between $t$ and $s$.)
The surface $S$ is isometric to a contractible subcomplex of the flat systolic plane and hence has a natural structure of a $CAT(0)$ space. The Euclidean geodesic is defined as a sequence of simplices in $S$, which runs near the $CAT(0)$ geodesic between $s$ and $t$.
\medskip

Now we pass to the more technical part of the exposition. Formally, the Euclidean geodesic is defined for a pair of simplices $\sigma,\tau$ in a systolic complex satisfying $\sigma\subset S_n(\tau), \tau\subset S_n(\sigma)$ for some $n\geq 0$ (where $S_n(\sigma)$ denotes the combinatorial sphere of radius $n$ around $\sigma$, c.f. Definition \ref{convex}).
The Euclidean geodesic is a certain sequence of simplices $\delta_k$, where $0\leq k\leq n$, such that $\delta_0=\sigma, \delta_n=\tau$,
which satisfies  $\delta_k\subset S_1(\delta_{k+1}),\delta_{k+1}\subset S_1(\delta_k)$ for
$0\leq k<n$ (c.f. Lemma \ref{properties of euclidean geodesics}(i)). The two most significant features of Euclidean geodesics are given by the following.

\begin{wtw}[Theorem \ref{corollary_to_main_theorem_on_euclidean_geodesics}]
\label{2}
Let $\sigma,\tau$ be simplices of a systolic complex $X$, such that for some natural $n$ we have
$\sigma \subset S_n(\tau), \tau \subset S_n(\sigma)$.
Let $(\delta_k)_{k=0}^n$ be the Euclidean geodesic between $\sigma$ and $\tau$. Take some $0\leq l<m\leq n$ and let $(r_k)_{k=l}^m$ be a 1--skeleton geodesic such that $r_k\in \delta_k$
for $l\leq k\leq m$.
Consider the simplices $\tilde{\delta}_l=r_l,\tilde{\delta}_{l+1},\ldots, \tilde{\delta}_m=r_m$
of the Euclidean geodesic between
vertices $r_l$ and $r_m$. Then for each $l\leq k\leq m$ we have $|\delta_k,\tilde{\delta}_k|\leq C$, where $C$ is a universal constant.
\end{wtw}

\begin{wtw}[Theorem \ref{contracting}]
\label{3}
Let $s,s',t$ be vertices in a systolic complex $X$ such
that $|st|=n,|s't|=n'$. Let $(r_k)_{k=0}^n,(r'_k)_{k=0}^{n'}$ be 1--skeleton geodesics such that
$r_k\in \delta_k,r'_k\in \delta'_k$, where $(\delta_k),(\delta'_k)$ are
Euclidean geodesics for $t,s$ and for $t,s'$ respectively. Then for all $0\leq
c\leq 1$ we have $|r_{\lfloor cn\rfloor}r'_{\lfloor cn'\rfloor}|\leq
c|ss'|+C$, where $C$ is a universal constant.
\end{wtw}

The article is organized as follows. It consists of an introductory part (Sections 1--2), the two main parts (Sections 3--6 and Sections 7--13), which can be read independently, and of a concluding Section \ref{final remarks}.

In Section \ref{Systolic complexes} we give a brief introduction to systolic complexes.

In the first part, assuming we have defined Euclidean geodesics satisfying Theorem \ref{2} and Theorem \ref{3}, we define the boundary:
In Section \ref{def} we define the boundary as a set
of equivalence classes of good geodesic rays. Then
we define topology on the compactification
obtained
by adjoining the boundary (Section \ref{top}) and we show its
compactness and finite dimensionality (Section \ref{cmpt}).
Finally, in Section \ref{ma}, we prove Theorem \ref{1}
---the main result of the paper.

In the second part of the article we define Euclidean geodesics and establish Theorem \ref{2} and Theorem \ref{3}:  In Section \ref{Flat surfaces} we recall Elsner's results on minimal surfaces. In Section \ref{Layers} we study \emph{layers}, whose union contains all geodesics between given vertices. We define Euclidean geodesics in Section \ref{Euclidean geodesic}.

In the next two sections we prove Theorem \ref{main_theorem_on_euclidean_geodesics} which is a weak version of Theorem \ref{2} (though with a better constant). Apart from the definitions these sections can be skipped by a hurried reader. We decided to include them since this way of obtaining (the weak version of) Theorem \ref{2} is straightforward in opposition to the strategy in Section \ref{Characteristic discs for Euclidean geodesics}, which is designed to obtain Theorem \ref{3}. In Section \ref{Directed geodesics between simplices of Euclidean geodesics} we study the position of directed geodesics between two simplices of a given Euclidean geodesic with respect to the minimal surface appearing in its construction. Then we verify Theorem \ref{main_theorem_on_euclidean_geodesics} in Section \ref{CAT(0) geometry of characteristic discs} by studying $CAT(0)$ geometry of minimal surfaces.

The last two sections are devoted to the proofs of Theorem \ref{2} and Theorem \ref{3}: In Section \ref{Characteristic discs for Euclidean geodesics} we prove (in a technically cumbersome manner) powerful Proposition \ref{Euclidean_near_CAT(0)} linked with $CAT(0)$ properties of the triangles, whose two sides are Euclidean
geodesics. Proposition \ref{Euclidean_near_CAT(0)} easily implies Theorem \ref{2}, but its main application comes in Section \ref{Contracting property}, where we use it to derive Theorem \ref{3}.

We conclude with announcing some further results for which we do not provide proofs in Section \ref{final remarks}.
\medskip

Acknowledgments. We are grateful to Tadeusz Januszkiewicz and Jacek \'Swi\k{a}tkowski for discussions right from the birth of our ideas and to Mladen Bestvina for encouragement. We thank the Mathematical Sciences Research Institute and the Institut des Hautes \'{E}tudes Scientifiques for the hospitality during the preparation of this article.

\section{Systolic complexes}
\label{Systolic complexes}
In this section we recall (from \cite {JS},\cite{JS2},\cite{SH}) the definition and basic properties of systolic
complexes and groups.

\begin{defin}
\label{full}
A subcomplex $K$ of a~simplicial complex $X$ is called
\emph{full} in $X$ if any simplex of $X$ spanned by vertices of
$K$ is a~simplex of $K$. The \emph{span} of a subcomplex $K\subset X$ is
the smallest full subcomplex of $X$ containing $K$. We denote
it by span$(K)$. A~simplicial complex $X$ is called
\emph{flag} if any set of vertices, which are pairwise connected
by edges of $X$, spans a~simplex in $X$. A simplicial complex $X$
is called \emph{$k$--large}, $\infty\geq k\geq 4$, if $X$ is flag and there
are no embedded cycles of length $<k$, which are full subcomplexes
of $X$ (i.e. $X$ is flag and every simplicial loop of length $<k$
and $\geq 4$ "has a diagonal").
\end{defin}

\begin{defin}
\label{systolic}
A simplicial complex $X$ is called
\emph{systolic} if it is connected, simply connected and links of
all simplices in $X$ are 6--large. A group $\Gamma$ is called
\emph{systolic} if it acts cocompactly and properly  by simplicial
automorphisms on a systolic complex $X$. (\emph{Properly} means
$X$ is locally finite and for each compact subcomplex $K\subset X$
the set of $\gamma\in \Gamma$ such that $\gamma(K)\cap K\neq
\emptyset$ is finite.)
\end{defin}
\medskip

Recall \cite{JS}, Proposition 1.4, that systolic complexes are
themselves 6--large. In particular they are flag. Moreover, we have the following.

\begin{theorem}[\cite{JS}, Theorem 4.1(1)]
\label{con}
Systolic complexes are contractible.
\end{theorem}

\medskip

Now we briefly treat the definitions and facts concerning
convexity.

\begin{defin}
\label{convex}
For every pair of subcomplexes (usually vertices) $A,B$ in
a~simplicial complex $X$ denote by $|A,B|$ ($|ab|$ for vertices $a,b$)
the combinatorial distance between $A^{(0)},B^{(0)}$ in $X^{(1)}$,
the 1--skeleton of $X$ (i.e. the minimal number of edges in a simplicial path connecting both sets).
A subcomplex $K$ of a~simplicial complex
$X$ is called \emph{3--convex} if it is a~full subcomplex of $X$
and for every pair of edges $ab,bc$ such that $a,c\in K, |ac|=2$,
we have $b\in K$. A subcomplex $K$ of a~systolic complex $X$ is
called \emph{convex} if it is connected and links of all simplices
in $K$ are 3--convex subcomplexes of links of those simplices in
$X$.
\end{defin}

In Lemma 7.2 of \cite{JS} authors conclude that convex
subcomplexes of a~systolic complex $X$ are full and
3--convex in $X$, and systolic themselves, hence contractible by Theorem \ref{con}.
The intersection of a family of convex subcomplexes is convex.
For a subcomplex $Y\subset X$, $n\geq 0$, the
\emph{combinatorial ball} $B_n(Y)$ of radius $n$ around $Y$ is
the span of $\{p\in X^{(0)} \colon |p,Y|\leq n\}$. (Similarly
$S_n(Y)= \Span\{p\in X^{(0)} \colon |p,Y|= n\}$.) If $Y$ is
convex (in particular, if $Y$ is a~simplex) then $B_n(Y)$ is also
convex, as proved in \cite{JS}, Corollary 7.5. Combining this with previous remarks we record:

\begin{cor}
\label{balls contractible}
In systolic complexes, balls around simplices are contractible.
\end{cor}

%

Haglund--\'Swi\k{a}tkowski prove the following.

\begin{proposition}[\cite{SH}, Proposition 4.9]
\label{second convexity}
A full subcomplex $Y$ of a systolic complex $X$ is convex if and
only if $Y^{(1)}$ is geodesically convex in $X^{(1)}$ (i.e. if all
geodesics in $X^{(1)}$ joining vertices of $Y$ lie in $Y^{(1)}$).
\end{proposition}

We record:

\begin{cor}
\label{balls convex}
In systolic complexes balls around simplices are geodesically convex.
\end{cor}

\medskip
We will need a~crucial "projection lemma".
The \emph{residue}
of a~simplex $\sigma$ in $X$ is the union of all simplices in $X$,
which contain $\sigma$.

\begin{lemma}[\cite{JS}, Lemma 7.7]
\label{projection lemma}
Let $Y$ be a~convex subcomplex of a~systolic complex $X$ and let
$\sigma$ be a~simplex in $S_1(Y)$. Then the
intersection of the residue of $\sigma$ and of the complex $Y$ is
a~simplex (in particular it is nonempty).
\end{lemma}

\begin{defin}
\label{projection}
The simplex as in Lemma \ref{projection lemma} is called the \emph{projection} of
$\sigma$ onto $Y$.
\end{defin}

The following lemma immediately follows from Definition \ref{projection}.

\begin{lemma}
\label{inclusions}
Let $\sigma\subset\tilde{\sigma}$ be simplices in $S_1(Y)$ for some convex $Y$ and let $\pi,\tilde{\pi}$
be their projections onto $Y$. Then $\tilde{\pi}\subset\pi$.
\end{lemma}

\begin{defin}
\label{directed geodesic}
For a~pair of vertices $v,w,\  |vw|=n$ in a~systolic
complex $X$ we define inductively a~sequence of simplices $\sigma_0=v,
\sigma_1,\ldots,\sigma_n=w$ as follows. Take $\sigma_{i}$ equal to the
projection of $\sigma_{i-1}$ onto $B_{n-i}(w)$ for $i=1,\ldots,n-1,n$. The
sequence $(\sigma_i)_{i=0}^n$ is called the \emph{directed geodesic} from $v$ to $w$
(this notion is introduced and studied in \cite{JS}).

We can extend this construction to any pair
$(\sigma_0,W)$, where $W$ is a convex subcomplex of $X$ and $\sigma_0$ is a
simplex. Namely, if for some $n$ we have $\sigma_0\subset S_n(W)$ then take
$\sigma_{i}$ to be the projection of $\sigma_{i-1}$ onto $B_{n-i}(W)$. If
$\sigma_0$ intersects both $S_n(W)$ and $S_{n-1}(W)$ then take
$\sigma_1=\sigma_0\cap S_{n-1}(W)$ and then proceed as previously. We call the
final $\sigma_n\subset W$ the \emph{projection} of $\sigma_0$ onto $W$. Note
that this coincides with Definition \ref{projection}. Observe that if $\sigma_0\subset W$
then the projection of $\sigma_0$ onto $W$ is equal to $\sigma_0$.
\end{defin}

Finally, recall a powerful observation.

\begin{lemma}[\cite{JS2}, Lemma 4.4]
\label{asphericity of full subcomplexes}
Every full subcomplex of a systolic complex is aspherical.
\end{lemma}

\section{Definition of the boundary}
\label{def}
Let $X$ be a systolic complex.
In this section we give two equivalent definitions of the boundary
of $X$ as a set. We use the notion of Euclidean geodesics which will be introduced in
Section \ref{Euclidean geodesic}, but actually we need only its features given by Theorem \ref{2} and Theorem \ref{3}. Thus, it is enough to read Sections \ref{introduction}--\ref{Systolic complexes} to follow the first part of the article (Sections \ref{def}--\ref{ma}). Let $C$ be a natural number, which is a universal constant satisfying assertions of both Theorem \ref{2} and Theorem \ref{3}.

\begin{rem}
\label{1--skeleton geodesic}
Let $(\delta_i)_{i=0}^n$ be a Euclidean geodesic and let $v_k$ be a vertex in $\delta_k$ for some $0\leq k\leq n$.
Then there exists a 1--skeleton geodesic $(v_i)_{i=0}^n$ such that $v_i\in\delta_i$ for $0\leq i\leq n$. This follows from the fact that $\delta_{i+1}\subset S_1(\delta_i)$, which we use for $k\leq i<n-1$, and from  $\delta_{i}\subset S_1(\delta_{i+1})$, which we use for $1\leq i <k$ (see Section \ref{introduction} or Lemma \ref{properties of euclidean geodesics}(i)).
\end{rem}

\begin{defin}
\label{good}
Let $v,w$ be vertices of a systolic complex $X$.
Let $\gamma=(v_0=v,v_1,v_2, \ldots ,v_n=w)$ be a geodesic in the $1$--skeleton of $X$ between $v$ and $w$ or let $\gamma=(v=v_0,v_1,v_2, \ldots )$
be a $1$--skeleton geodesic ray starting at $v$ (then we set $n=\infty$).
For $0\leq i< j \leq n$, by $(\delta_i^{i,j}=v_i,\delta_{i+1}^{i,j}, \ldots
,\delta_{j}^{i,j}=v_j)$ we denote the Euclidean geodesic between
$v_i$ and $v_j$. We say that $\gamma$ is a \emph{good geodesic} between
$v$ and $w$ or that $\gamma$ is a \emph{good geodesic ray} starting at $v$
if for every $0\leq i< j \leq n$ and every $i\leq k \leq j$ we have
$|v_k,\delta_{k}^{i,j}|\leq C+1$ (the constant $C$ is defined in the beginning of this section).

By ${\cal R}$ we denote the set of all good geodesic rays
in $X$. For a given vertex $O$ of $X$, by ${\cal R}_O$ we denote the set of all good geodesic rays starting at $O$.
\end{defin}

The following two results are immediate corollaries of Theorem
\ref{2} and Theorem \ref{3}.

\begin{cor}
\label{goodex}
For every two vertices $v,w\in X$ there exists a good geodesic between them.
\end{cor}
\dow
Let $(\delta_0=v,\delta_1, \ldots ,\delta_n=w)$ be the Euclidean geodesic
between $v$ and $w$. By Remark \ref{1--skeleton geodesic}, there exists
a $1$--skeleton geodesic
$\gamma=(v_0=v,v_1,v_2, \ldots ,v_n=w)$ with $v_i\in \delta_i$.
We claim that $\gamma$ is a good geodesic.
To justify the claim let $0\leq i<j\leq n$. Let $(\widetilde{\delta}_i, \widetilde{\delta}_{i+1},\ldots, \widetilde{\delta}_j)$
be the Euclidean geodesic between
$v_i$ and $v_j$. By Theorem
\ref{2}, for every $i\leq k \leq j$,
we have
$$
|v_k, \widetilde{\delta}_k|\leq
|\delta_k,\widetilde{\delta}_k|+1 \leq C+1,
$$
which justifies the claim.
\kon

\begin{cor}
\label{contr}
Let $(v_0=O,v_1,v_2, \ldots ,v_n),\ (w_0=O,w_1,w_2, \ldots ,w_m)$
be good geodesics in $X$. Then for all $0\leq c\leq 1$ we have
$|v_{\lfloor cn \rfloor}w_{\lfloor cm\rfloor}|\leq
c|v_nw_m|+D$, where $D=3C+2$.
\end{cor}
\dow
Let $(\delta_i^v),(\delta_i^w)$ be the Euclidean geodesics
between $O$ and, respectively, $v_n$, $w_m$. Fix $0\leq c\leq 1$.
Pick vertices $v_{\lfloor cn\rfloor}'\in \delta^v_{\lfloor cn \rfloor}$ and
$w_{\lfloor cm\rfloor}'\in \delta^w_{\lfloor cm \rfloor}$ which realize the distance to $v_{\lfloor cn\rfloor},w_{\lfloor cm\rfloor}$, respectively.
Find 1--skeleton geodesics $(v_i')_{i=0}^{\lfloor cn\rfloor}$ and $(w_i')_{i=0}^{\lfloor cm\rfloor}$
such that $v_i'\in \delta_i^v$ and $w_i'\in \delta_i^w$.
Their existence is guaranteed by Remark
\ref{1--skeleton geodesic}.
By Theorem \ref{3}, we have
\begin{align*}
|v_{\lfloor cn \rfloor}w_{\lfloor cm\rfloor}|&\leq
|v_{\lfloor cn \rfloor}v_{\lfloor cn \rfloor}'|+
|v_{\lfloor cn\rfloor}'w_{\lfloor cm\rfloor}'|+
|w_{\lfloor cm\rfloor}'w_{\lfloor cm\rfloor}|=
\\
&= |v_{\lfloor cn \rfloor},\delta^v_{\lfloor cn \rfloor}|+
|v_{\lfloor cn\rfloor}'w_{\lfloor cm\rfloor}'|+
|\delta^w_{\lfloor cm\rfloor},w_{\lfloor cm\rfloor}|\leq
\\
&\leq(C+1)+(c|v_nw_m|+C)+(C+1),
\end{align*}
as desired.
\kon
\medskip

The following simple corollary of Corollary \ref{contr} will be useful.

\begin{cor}
\label{contr'}
Let $(v_0=O,v_1,v_2, \ldots ,v_k),(w_0=O,w_1,w_2, \ldots ,w_l)$
be good geodesics in $X$. Then for all $0\leq N\leq \min \lk k,l \rk$ we have
$|v_Nw_N|\leq 2|v_kw_l|+D$.
\end{cor}
\dow
W.l.o.g. we can assume that $k\leq l$. Observe that $l-k\leq |v_kw_l|$. Hence, by Corollary \ref{contr}, we have
$$
|v_Nw_N|\leq |v_kw_k| +D\leq |v_kw_l|+|w_lw_k|+D=
|v_kw_l|+(l-k)+D\leq 2|v_kw_l|+D.
$$
\kon
\medskip

Below we define the central object of the article.

\begin{defin}
\label{def1}
The \emph{(ideal) boundary} of a systolic complex $X$ is the set
$\partial X={\cal R}/ \sim$ of equivalence classes
of good geodesic rays, where
rays
$\eta=(v_0,v_1,v_2, \ldots ),\ \xi=(w_0,w_1,w_2, \ldots )$ are identified if $|v_iw_i|$ is bounded above by a constant independent of $i$
(one can check this happens exactly when the Hausdorff distance between $\eta$ and $\xi$ is finite). For a good geodesic ray $\eta$, we denote
its equivalence class in $\partial X$ by $[\eta]$.
\end{defin}

In order to introduce topology on $\cx=X\cup \partial X$ we give
another definition of the boundary. The two definitions
will turn out to be equivalent in the case of a systolic
complex with a geometric group action.

\begin{defin}
\label{def2}
Let $O$ be a vertex of a systolic complex $X$. Then
the \emph{(ideal) boundary of $X$ with respect to the basepoint vertex $O$}
is the set
$\partial_O X={\cal R}_O/ \sim$ of equivalence classes
of good geodesic rays starting at $O$, where
rays
$\eta=(v_0=O,v_1,v_2, \ldots ),\ \xi=(w_0=O,w_1,w_2, \ldots )$ are identified if $|v_iw_i|$ is bounded above by a constant independent of $i$
(again this happens exactly when the Hausdorff distance between $\eta$ and $\xi$ is finite).
For $\eta \in {\cal R}_O$, we denote
its equivalence class in $\bdo X$ by $[\eta]$ (we hope this ambiguity of the notation will not cause confusion).
\end{defin}

\begin{lem}
\label{repr}
Let $\eta=(v_0=O,v_1,v_2, \ldots ),\ \xi=(w_0=O,w_1,w_2, \ldots )\in {\cal R}_O$. Then $[\eta]=[\xi]$ iff $|v_iw_i|\leq D$ for all $i$.
\end{lem}
\dow
We show that if for some $i$ we have $|v_iw_i|-D\geq 1$, then
$[\xi]\neq [\eta]$. Let $i$ be as above and $R$ be a natural number.
Then, by Corollary \ref{contr}, we have
$$|v_{Ri}w_{Ri}|\geq R(|v_{i}w_{i}|-D)\geq R.$$
Since $R$ can be chosen arbitrarily large,
we get $[\xi]\neq [\eta]$.\kon
\medskip

In the remaining part of this section we prove equivalence
of the above two notions of boundaries in the case
of locally finite complexes.
Assume that $X$ is a locally finite systolic complex.
Let $O\in X$ be a fixed vertex and let
$\eta=(v^0,v^1,v^2, \ldots )$ be a good
geodesic
ray in $X$.
For every $i\geq 0$ we
choose a good geodesic
$\eta^i=(v_0^i=O,v_1^i,v_2^i, \ldots ,v_{n(i)}^i=v^i)$, guaranteed by Corollary \ref{goodex}.
Since $B_1(O)$ is finite, for some vertex $v_1\in S_1(O)$ there are infinitely many $i$ such that $n(i)=|Ov^i|\geq 1$ and $v^i_1=v_1$. Similarly, since all balls are finite, we obtain inductively vertices $v_k\in S_k(O)$ satisfying the following. For each $k$ there are infinitely many $i$ such that $n(i)\geq k$ and for all $j\leq k$ we have $v^i_j=v_j$. For each $k$ denote some such $i$ by $i(k)$.
The following easy facts hold.

\begin{lem}
\label{vv'}
The sequence $(v_0=O,v_1,v_2, \ldots )$ obtained as above
is a good geodesic ray. Moreover, for every $j$ we have
$|v^jv_j|\leq 3|Ov^0|+D$.
\end{lem}
\dow
The first assertion follows from the fact that
for every $k$ the sequence $(v_0=O,v_1,v_2, \ldots ,v_k)$
is a subsequence of the good geodesic
$\eta^{i(k)}$ and hence, by Definition \ref{good},
it is a good geodesic.

Now we prove the second assertion. Let $j\geq 0$.
Consider the case $n(i(j))\leq i(j)$ (the case
$n(i(j))> i(j)$ can be examined analogically).
Then for $k=i(j)-n(i(j))$
we have $|v^kv^{i(j)}|=|v_0^{i(j)}v^{i(j)}|$.
Thus we can
apply Corollary \ref{contr} with $m=n$ to good geodesics
$\eta^{i(j)}$ and $(v^k,v^{k+1}, \ldots ,v^{i(j)})$, which yields the following.
$$
|v^{k+j}v_j|=|v^{k+j}v^{i(j)}_j|\leq |v^kv_0^{i(j)}|+D\leq
(|Ov^0|+k)+D.
$$
Hence
$$ |v^jv_j|\leq k+|v^{k+j}v_j|\leq |Ov^0|+2k+D\leq 3|Ov^0|+D,$$
where the last inequality follows from $k\leq |Ov^0|$, which is the triangle inequality for $v^0, v^{i(j)}$ and $O$.
\kon

\begin{cor}
\label{OO'}
Let $X$ be a locally finite systolic complex and $O,O'$ its vertices.
Then the map $\Phi_{O} \colon \partial X \to \bdo X$ given
by $\Phi_{O}([(v^0,v^1,v^2, \ldots )])=[(v_0=O,v_1,v_2, \ldots )]$ is well defined.
It is a bijection between $\partial X$ and $\bdo X$.
Its restriction $\Phi_{O'O}=\Phi_O|_{\bdoo X}$ is a bijection between
$\bdoo X$ and $\bdo X$.
\end{cor}

\section{Topology on $\cx=X\cup \bdo X$}
\label{top}
Let $X$ be a systolic complex and $O\in X$ be its vertex.
In this section we
define the topology on the set $\overline{X}=X\cup \bdo X$, which extends
the usual topology on the simplicial complex $X$.
The idea is to define the topology through neighborhoods
(not necessarily open) of points in $\cx$. The only problem is to define the neighborhoods of points in the boundary.


For a $1$--skeleton geodesic or a geodesic ray
$\eta=(v_0,v_1,v_2, \ldots )$, we denote by $B_1(\eta)$
the combinatorial ball of radius $1$ around
the subcomplex $\lk v_0, v_1, v_2, \ldots  \rk$.
Let $C$ and $D=3C+2$ be the constants from the previous section.
\begin{defin}
\label{otocz}
Let $\eta=(v_0=O,v_1,v_2, \ldots )$ be a good geodesic ray in $X$ and let $R>D$ (i.e. $R\geq D+1$) and $N\geq 1$ be natural numbers (in fact we could also allow $N=0$, but this would complicate some computations later on). By ${\cal G}_O(\eta,N,R)$ we denote the set
of all good geodesics $(w_0=O,w_1, \ldots ,w_k)$ with $k\geq N$ and good
geodesic rays $(O=w_0,w_1, \ldots )$, such that $|w_Nv_N|\leq R$.

By ${\cal G'}_O(\eta,N,R)$ we denote the set $\lk (w_N,w_{N+1}, \ldots )|\;
(w_0=O,w_1, \ldots ) \in {\cal G}_O(\eta,N,R) \rk$.

A \emph{standard neighborhood} of $[\eta]\in \bdo X\subset \cx$ is the set
$$
U_O(\eta,N,R)=\lk [\xi]|\; \xi \in  {\cal G}_O(\eta,N,R)\cap {\cal R}_O \rk \cup
\bigcup \lk  \mr{int}B_{1}(\xi)|\; \xi \in {\cal G'}_O(\eta,N,R)\rk.
$$

We write ${\cal G}(\eta,N,R)$, ${\cal G'}(\eta,N,R)$ and
$U(\eta,N,R)$ instead of ${\cal G}_O(\eta,N,R)$, ${\cal G'}_O(\eta,N,R)$ and
$U_O(\eta,N,R)$ if it is clear what is the basepoint $O$.
\end{defin}

Before we define the topology, we need the following useful lemmas.
The first one is an immediate consequence of the above definition.

\begin{lem}
\label{U<U'}
Let $\eta,\xi \in {\cal R}_O$ and let $N,N',R>D,R'>D$ be natural
numbers such that $N'\geq N$. If ${\cal G}(\xi,N',R')\subset
{\cal G}(\eta,N,R)$ then $U(\xi,N',R')\subset U(\eta,N,R)$.
\end{lem}

\begin{lem}
\label{RR'}
Let $U(\eta,N,R)$ be a standard neighborhood,
let $\xi\in {\cal R}_O$ be such that $[\xi]=[\eta]$ and
let $R'\geq D$ be a natural number. Then, for $N'\geq (R'+D)N$, we have
$U(\xi,N',R')\subset U(\eta,N,R)$.
\end{lem}
\dow
Denote $\eta=(v_0=O,v_1,v_2, \ldots )$ and $\xi=(w_0=O,w_1,w_2, \ldots )$.

By Lemma \ref{U<U'}, it is enough to show that
for every $\zeta \in {\cal G}(\xi,N',R')$ we have
$\zeta \in {\cal G}(\eta,N,R)$.

Let $\zeta=(z_0=O,z_1,z_2, \ldots )\in {\cal G}(\xi,N',R')$.
By Corollary \ref{contr} and Lemma \ref{repr}, we have
\begin{align*}
|z_Nv_N| & \leq \frac{1}{R'+D}|z_{N'}v_{N'}|+D\leq
\frac{1}{R'+D}(|z_{N'}w_{N'}|+|w_{N'}v_{N'}|)+D\leq
\\
& \leq \frac{1}{R'+D}(R'+D)+D\leq R.
\end{align*}
Thus $\zeta \in {\cal G}(\eta,N,R)$ and the lemma follows.\kon
\medskip

The following defines topology on $\cx$.

\begin{prop}
\label{topo}
Let ${\cal A}$ be the family of subsets $A$ of $\cx=X\cup \bdo X$
satisfying the following. $A\cap X$ is open in $X$ and for every
$x\in A\cap \bdo X$ there is some $\eta\in {\cal R}_O$ such that $[\eta]=x$ and there is a standard neighborhood
$U(\eta,N,R)\subset A$.
Then ${\cal A}$ is a topology on $\cx$.
\end{prop}
\dow
The only thing we have to check is the following. If $A_1,A_2\in {\cal A}$
and $[\eta]\in A_1\cap A_2\cap \bdo X$ then there is a standard
neighborhood $U(\eta,N,R)$ of $[\eta]$ contained in $A_1\cap A_2$.

Since $A_1,A_2\in {\cal A}$, for $i=1,2$, there are standard neighborhoods  $U(\eta_i,N_i,R_i)\subset A_i$
such that $[\eta_i]=[\eta]$.
Thus, by Lemma \ref{RR'}, for any natural $R>D$ there exists
$N \geq N_i$, $i=1,2$, such that
$U(\eta,N,R)\subset U(\eta_1,N_1,R_1)\cap U(\eta_2,N_2,R_2)\subset
A_1\cap A_2$.\kon
\medskip

\begin{rem}
\label{Gromovbd}
It is easy to verify that when $X$ is $\delta$--hyperbolic (in the sense
of Gromov) then our boundary $\bdo X$ (with topology induced from $\cx$) is homeomorphic in a natural way
with the Gromov boundary of $X$.
\end{rem}

We still did not prove that the topology defined in Proposition \ref{topo} is non--trivial. This will follow from the next two lemmas, in which we characterize the intersections with the boundary of the interiors of standard neighborhoods. In particular, we show that
$[\xi]$ is contained in the interior of
$U(\xi,N,R)$.

\begin{lem}
\label{in}
For a set $A\subset \cx$, the intersection ${\mr {int}}A\cap \bdo X$
consists of those points $x\in \bdo X$ for which there exists a representative $\eta$ with
a standard neighborhood $U(\eta, N,R)\subset A$.
\end{lem}
\dow
Let $B$ be the set of those points $x\in \bdo X$ for which there exists a representative $\eta$ of $x$ with
a standard neighborhood $U(\eta, N,R)\subset A$.

It is clear that ${\mr {int}}A\cap \bdo X \subset B$, since ${\mr {int}}A$ is open in the topology defined in Proposition \ref{topo}.
We want now to prove the converse inclusion $B\subset {\mr {int}}A\cap \bdo X$. It is clear that $B\subset A\cap \bdo X$. Thus to prove the lemma we only have to show that $B$ is open in $\bdo X$ (in the topology induced from $\cx$).

Let $x\in B$ and let its representative $\eta$ be such that the standard neighborhood
$U(\eta,N,R')$ is contained in $A$.
By Lemma \ref{RR'}, we can assume
that $R'\geq 2(D+1)$.
Choose natural number
$N'\geq RN$.
We claim that $U(\eta,N',R)\cap \bdo X\subset B$ (i.e. that  equivalence classes of elements in ${\cal G}(\eta,N',R)\cap {\cal R}_O$ lie in $B$).
This implies that $B$ is open in $\bdo X$.

To justify the claim
let $\xi\in {\cal G}(\eta,N',R)\cap {\cal R}_O$.
To prove that $[\xi] \in B$ it is enough to establish $U(\xi,N',R)\subset U(\eta,N,R')$, since the latter is contained in $A$.
By Lemma \ref{U<U'}, it is enough to show that for
every $\zeta\in {\cal G}(\xi,N',R)$,
we have $\zeta \in {\cal G}(\eta,N,R')$.
Let $\zeta=(z_0=O,z_1, \ldots )\in {\cal G}(\xi,N',R)$. Denote $\eta=(v_0=O,v_1, \ldots ), \ \xi=(w_0=O,w_1, \ldots )$.

By Corollary \ref{contr}, we have
\begin{align*}
|z_Nv_N|& \leq |z_Nw_N|+|w_Nv_N|\leq
\\
& \leq \Big(\frac{1}{R}|z_{N'}w_{N'}|+D\Big)
+\Big(\frac{1}{R}|w_{N'}v_{N'}|+D\Big)\leq
\\
& \leq\Big(\frac{1}{R}R+D\Big)+\Big(\frac{1}{R}R+D\Big)=2(D+1)\leq R'.
\end{align*}
Thus $\zeta\in {\cal G}(\eta,N,R')$ and it follows that
$U(\xi,N',R)\subset U(\eta,N,R')$, which justifies the claim.
\kon

\begin{lem}
\label{int}
Let $U(\eta,N,R)$ be a standard neighborhood. Let
$\xi=(w_0=O,w_1,w_2, \ldots ) \in {\cal R}_O$ be such that
$v_N=w_N$, where $\eta=(v_0=O,v_1,v_2, \ldots )$.
Then $[\xi]$ is contained in the interior of
$U(\eta,N,R)$.
\end{lem}
\dow
By Lemma \ref{in}, it is enough to show
that there exists a standard neighborhood
$U(\xi,N',R)$ of $[\xi]$
contained in $U(\eta,N,R)$.
Let $N'\geq RN$.
By Lemma \ref{U<U'}, it is enough to show that
for $(z_0=O,z_1,z_2, \ldots )\in {\cal G}(\xi, N',R)$, we have
$|z_Nv_N|\leq R$.
By Corollary \ref{contr}, we have
\begin{align*}
|z_Nv_N|= |z_Nw_N|\leq
\frac{1}{R}|z_{N'}w_{N'}|+D \leq \frac{1}{R}R+D\leq R ,
\end{align*}
as desired.
\kon
\medskip

Below we give a sufficient condition for two standard neighborhoods to be disjoint.

\begin{lem}
\label{disjoint}
Let $U(\eta,N,R)$ and $U(\xi,N,S)$ be two standard
neighborhoods, with $\eta=(v_0=O,v_1,v_2, \ldots )$ and
$\xi=(w_0=O,w_1,w_2, \ldots )$. If $|v_Nw_N|> R+S+D+2$,
then $U(\eta,N,R)\cap U(\xi,N,S)=\emptyset$.
\end{lem}
\dow
By contradiction. Assume $U(\eta,N,R)\cap U(\xi,N,S)
\neq \emptyset$.
\medskip

\textbf{Case 1:} Let $x\in U(\eta,N,R)\cap U(\xi,N,S)\cap X$.
Then, by Definition \ref{otocz}, there exist $\eta'=(v_0'=O,v_1',v_2', \ldots )\in
{\cal G}(\eta,N,R)$ and $\xi'=(w_0'=O,w_1',w_2', \ldots )
\in {\cal G}(\xi,N,S)$ such that
$x$ belongs to the interior of both some simplex with vertex $v_k'$ and some simplex with vertex $w_l'$,
for some $k,l\geq N$.
Then these simplices coincide and $|v_k'w_l'|\leq 1$.
By Corollary \ref{contr'},
we have
$$|v_Nw_N|\leq|v_Nv_N'|+|v_N'w_N'|+|w_N'w_N|\leq R+(2|v_k'w_l'|+D)+S\leq R+(2+D)+S,$$
contradiction.
\medskip

\textbf{Case 2:} Let $\eta'=(v_0'=O,v_1',v_2', \ldots)\in
{\cal G}(\eta,N,R)$ and
$\xi'=(w_0'=O,w_1',w_2', \ldots ) \in {\cal G}(\xi,N,S)$
be such that $[\eta']=[\xi']$. Then, by Lemma \ref{repr}, we get
$$|v_Nw_N|\leq|v_Nv'_N|+|v'_Nw'_N|+|w'_Nw_N|\leq
R+D+S,$$
contradiction.
\kon

\section{Compactness and finite dimensionality}
\label{cmpt}
Let $X$ be a locally finite systolic complex and let $O\in X$ be its vertex.
In this section we show that $\cx=X\cup \bdo X$ is compact metrizable
and (if $X$ satisfies some additional local finiteness
conditions) finitely dimensional.
We also prove that, for a different vertex $O'$
of $X$, the compactifications $X\cup \bdo X$ and $X\cup \bdoo X$
are homeomorphic.

\begin{prop}
\label{reg}
If $X$ is locally finite then the space
$\cx=X\cup \bdo X$ is second countable and regular.
\end{prop}
\dow
It is clear that $\cx$ is second countable. We show that $\cx$ is regular.

First we show that $\cx$ is Hausdorff. We consider only the case
of two points of the boundary---the other cases are obvious.
Let $[\eta]\neq [\xi]$ be
two boundary points with
$\eta=(v_0=O,v_1,v_2, \ldots )$ and
$\xi=(w_0=O,w_1,w_2, \ldots )$.
Fix a natural number $R>D$ (for example $R=D+1$).
We can find $N$ such that $|v_Nw_N|>2R+D+2$.
Then, by Lemma \ref{int}, we have
$[\eta]\in {\mr {int}}U(\eta,N,R)$ and $[\xi]\in {\mr {int}}U(\xi,N,R)$
 and, by Lemma \ref{disjoint}, we get
${\mr {int}}U(\eta,N,R)\cap {\mr {int}}U(\xi,N,R)\subset
U(\eta,N,R)\cap U(\xi,N,R)=\emptyset$. Thus
we get disjoint non--empty open neighborhoods of $[\eta]$ and $[\xi]$.
\medskip

To show that $\cx$ is regular it is now enough, for every point
$x\in \cx$ and every closed subset $A\subset \cx$ which does not contain $x$,
to find disjoint open sets $U,V$ such that $x\in U$ and $A\subset V$.
Let $x\notin A$ be as above. The  case $x\in X$ is obvious thus we
consider only the case of $x=[\eta]\in \bdo X$, for
$\eta=(v_0=O,v_1,v_2, \ldots )$. Fix some natural $R>D$. Since $\cx \setminus A$ is open,
by definition of the topology (Proposition \ref{topo}) and by Lemma \ref{RR'},
we can find a natural number $N>0$ such that
$U(\eta,N,R')\subset \cx \setminus A$, where $R'\geq 2D+2$.
Let $N'=(R+1)N+1$ and let $U={\mr {int}}U(\eta,N',R)$.
Observe that, by Lemma \ref{int}, we have $x\in U$.
Now we define $V$. For each $y\in A\cap X$, choose an open set $V_y={\mr {int}}B_1(z')$ for some vertex $z'$ in $X$ such that $y\in {\mr {int}}B_1(z')$.
Then we set $V=\bigcup \lk V_y
\; | \; \; y\in A\cap X \rk \cup \bigcup \lk  {\mr {int}}U(\xi,N',R)\; |
\; \; [\xi]\in A\cap \bdo X\rk$.
By Lemma \ref{int}, we have $A\cap \bdo X \subset V$, hence $A\subset V$.
Thus to prove that $U$ and $V$
are as desired we only need to show that $U\cap V=\emptyset$.
\medskip

First we prove that $U\cap {\mr {int}}U(\xi,N',R)=\emptyset$,
for $[\xi]\in A\cap \bdo X$. Let $\xi=(w_0=O,w_1,w_2, \ldots )$.
Then, by Corollary \ref{contr} and by $A\cap U(\eta,N,R')=\emptyset$, we have
\begin{align*}
|v_{N'}w_{N'}|&\geq \frac{N'}{N}(|v_Nw_N|-D)>
(R+1)(R'-D)\geq
\\
&\geq (R+1)(D+2)> 2R+D+2.
\end{align*}
Thus, by Lemma \ref{disjoint},
$U\cap {\mr {int}}U(\xi,N',R)\subset U(\eta,N',R)\cap U(\xi,N',R)=
\emptyset$.
\medskip

Now we show that $U\cap V_y=\emptyset$, for $y\in A\cap X$.
By contradiction, assume $p\in U\cap V_y$.
Since $p\in U$, there exist a vertex $z$ of the simplex containing $p$ in its interior and a good
geodesic $(z_0=O,z_1, \ldots ,z_k=z)\in {\cal G}(\eta,N',R)$,
where $k\geq N'$.
Then, by Corollary \ref{contr}, we have
$$
|v_Nz_N|\leq \frac{N}{N'}|v_{N'}z_{N'}|+D<
\frac{1}{R}R+D\leq D+1.
$$
On the other hand, since $p\in V_y$, there is a vertex
$z'$ such that $\lk y,p \rk \in {\mr {int}}B_1(z')$.
Then
$|zz'|\leq 1$.
Let $(O=z_0',z_1', \ldots ,z_l'=z')$ be a good geodesic.
We have $l\geq N'-1$, hence by Corollary \ref{contr} and Corollary \ref{contr'}, we get
\begin{align*}
|z_{N}z_{N}'|& \leq \frac{N}{N'-1}|z_{N'-1}z_{N'-1}'|+D\leq
\frac{1}{R+1}(2|zz'|+D)+D\leq
\\
&\leq \frac{1}{D+2}(2+D)+D= D+1.
\end{align*}
Summarizing, we have $|v_Nz_N'|\leq |v_Nz_N|+|z_{N}z_{N}'|\leq 2D+2\leq R'$.
It follows that
$(O=z_0',z_1', \ldots ,z_L'=z')\in {\cal G}(\eta,N,R')$ and hence
$y\in U(\eta,N,R')$---contradiction.
\kon

\begin{cor}
\label{metr}
If $X$ is locally finite then the space $\cx=X\cup \bdo X$ is metrizable.
\end{cor}
\dow
This follows from the Urysohn Metrization Theorem---cf.
\cite[Corollary 9.2]{Du}.
\kon

\begin{prop}
\label{cpt}
If $X$ is locally finite then the space
$\cx=X\cup \bdo X$ is compact.
\end{prop}
\dow
By Corollary \ref{metr}, it is enough to
show that every infinite sequence of points in $\cx$
contains a convergent subsequence. Let $(x^1,x^2,x^3, \ldots )$ be a given
sequence of points in $\cx$. If for some $n>0$ there is only finitely
many $x^i$ outside the ball $B_n(O)$ (which is finite), then we can find
a convergent subsequence.
From now on we assume there is no $n$ as above.

For every $i$ we
choose a good geodesic
 or a good geodesic ray $\eta^i=(v_0^i=O,v_1^i,v_2^i, \ldots )$
 the following way.
If $x^i\in X$ then $\eta^i=(v_0^i=O,v_1^i,v_2^i, \ldots ,v_{n(i)}^i)$
is a good geodesic between $O$ and a vertex $v_{n(i)}^i$ lying in a
common simplex with the point $x^i$.
If $x^i\in \bdo X$ then we set
$\eta^i=\zeta$ for an arbitrary $\zeta$ such that $x^i=[\zeta]$ and
we set $n(i)=\infty$.
By our assumptions on $(x^1,x^2,x^3, \ldots )$, for every $n>0$ there exists
an arbitrarily large $i$ such that $n(i)>n$.
Since $S_1(O)$ is finite, for some vertex $v_1\in S_1(O)$ there are infinitely many $i$ such that $n(i)\geq 1$ and $v^i_1=v_1$.
Let $i(1)$ be some such $i$. Similarly, since all spheres are finite, we obtain inductively vertices $v_k\in S_k(O)$ and numbers $i(k)$ satisfying the following. For each $k$ there are infinitely many $i$ such that $n(i)\geq k$ and for all $j\leq k$ we have
$v^i_j=v_j$; we denote some such $i>i(k-1)$ by $i(k)$.

Observe that for every $k$ the sequence $(v_0=O,v_1,v_2, \ldots ,v_k)$
is a subsequence of the good geodesic or the good geodesic ray
$\eta^{i(k)}$ and hence, by Definition \ref{good}, it
is a good geodesic. Thus every subsequence of the infinite
sequence $(v_0=O,v_1,v_2, \ldots )$ is a good geodesic and again, by
Definition \ref{good}, $(v_0=O,v_1,v_2, \ldots )$ is a good geodesic ray.

We claim that
the sequence $(x^{i(k)})_{k=1}^{\infty}$ of points of $\cx$
converges
to $[\eta]\in \bdo X$, where $\eta=(v_0=O,v_1,v_2, \ldots )$.
To prove the claim it is enough to show (since every open set containing
$[\eta]$ contains some $U(\eta, N, R)$, by Lemma \ref{RR'}) that we have $\eta^{i(k)} \in {\cal G}(\eta,N,R)$,
for every $k\geq N$.
This follows from the equality $v_N^{i(k)}=v_N$, which holds for every $k\geq N$.
\kon
\medskip

Observe that by the above proof we get the following.

\begin{cor}
\label{nonemp}
If a locally finite systolic complex is unbounded then its
boundary is non--empty.
\end{cor}

Below we prove that the bijection
$\Phi_{O'O}$ defined in Corollary \ref{OO'} extends to a homeomorphism
of compactifications coming from different basepoints.

\begin{lem}
\label{OO''}
Let $X$ be a locally finite systolic complex and
let $O,O'$ be its vertices. Then the map
$\Phi_{O'O} \colon X\cup \bdoo X \to X\cup \bdo X$ defined as an extension
by identity on $X$ of the map
$\Phi_{O'O} \colon \bdoo X \to \bdo X$ is a homeomorphism.
\end{lem}
\dow
By compactness (Proposition \ref{cpt}) and by Corollary \ref{OO'}, we only have to show that
$\Phi_{O'O}$ is continuous. It is enough to check the continuity
at the boundary points. Let $\xi=(v_0=O,v_1,v_2, \ldots )$ be obtained
from a good geodesic ray $\eta=(v^0=O',v^1,v^2, \ldots )$ as in the definition
of the map $\Phi_{O'O}$. We show that $\Phi_{O'O}$ is continuous at
$[\eta]$. Let $d=|OO'|$, let $R>D$ be a natural number and $R'=R+3D+6d$
and let $U$ be an open neighborhood
of $[\xi]$ in $X\cup \bdo X$.
We have to show that there exists an open neighborhood $V$
of $[\eta]$ in $X\cup \bdoo X$ such that $\Phi_{O'O}(V)\subset U$.
By Lemma \ref{RR'},
there exists
$N$ such that $U_O(\xi,N,R')\subset U$.
Let $V={\mr {int}}U_{O'}(\eta,N+d,R)$. By Lemma \ref{int}, $[\eta]\in V$. We claim that $\Phi_{O'O}(V)\subset U$---this will finish the proof.
\medskip

First we show that for $x\in V\cap X$
we have $\Phi_{O'O}(x)=x\in U$. For such an $x$ choose, by definition
of $U_{O'}(\eta,N+d,R)$, a good geodesic $(w^0=O',w^1,w^2, \ldots ,w^k)\in
{\cal G}_{O'}(\eta,N+d,R)$ such that $x$ belongs to the interior of a simplex with vertex $w^k$, where $k\geq N+d$.
Let $\zeta=(w_0=O,w_1,w_2, \ldots ,w_l=w^k)$ be a good geodesic
guaranteed by Corollary \ref{goodex}.
Then $|l-k|\leq d$, hence $l\geq N$ and $w_N$ is defined. By Lemma \ref{vv'} and Corollary \ref{contr}, we have
\begin{align*}
|w_Nv_N|&\leq
|w_Nw^N|+|w^Nv^N|+|v^Nv_N|\leq
\\
&
\leq (3d+D)+(|w^{N+d}v^{N+d}|+D)+(3d+D)\leq
\\
&\leq R+3D+6d= R'.
\end{align*}
This inequality implies that
$\zeta \in {\cal G}_{O}(\xi,N,R')$ and hence
$x\in U_{O}(\xi,N,R')\subset U$.
\medskip

Now we show that for $[\rho]\in V\cap \bdoo X$ we have
$\Phi_{O'O}([\rho])\in U$.
Let
$\rho=(w^0=O',w^1,w^2, \ldots )\in
{\cal G}_{O'}(\eta,N+d,R)\cap {\cal R}_{O'}$.
Let $\zeta=(w_0=O,w_1,w_2, \ldots )$ be obtained
from $\rho$ as in the definition of $\Phi_{O'O}$.
Then, by Lemma \ref{vv'} and Corollary \ref{contr}, we can perform the same computation as in the previous case to get $|w_Nv_N|\leq R+3D+6d= R'$.
Thus $\Phi_{O'O}([\rho])=[\zeta]\in U_O(\xi,N,R')\subset U$ and we
have completed the proof of
$\Phi_{O'O}(V)\subset U$ and of the whole lemma.\kon
\medskip

Now we proceed to the question of finite dimensionality of $\cx$.
Let us remind that a simplicial complex $X$ is \emph{uniformly locally
finite} if there exists a natural number $L$ such that every
vertex belongs to at most $L$ different simplices.
This happens for example when some group acts geometrically on $X$.

\begin{prop}
\label{findim}
Let $X$ be a uniformly locally finite systolic complex.
Then $\cx=X\cup \bdo X$ is finitely dimensional.
\end{prop}
\dow
Recall that a space $Y$ has \emph {dimension at most $n$} if, for every open cover
${\cal U}$ of $Y$, there exists an open cover ${\cal V}\prec {\cal U}$
($\cal V$ is a \emph{refinement} of $\cal U$, i.e. every element
of ${\cal V}$ is contained in some
element of ${\cal U}$)
such that every point in $Y$ belongs to at most $n+1$ elements of ${\cal V}$
(i.e. the \emph{multiplicity} of $\cal V$ is at most $n+1$).

It is clear that $X$ is finitely dimensional.
It is thus enough to show that there exists a constant $K$ such that
for every open (in $\cx$) cover $\cal U$ of $\bdo X$ there
exists an open cover $\cal V\prec \cal U$ of $\bdo X$ of multiplicity
at most $K$.

Let $R>D$ be a natural number. Then, by uniform local finiteness, there
is a constant $K$ such that every ball of radius at most $2R+D+2$
contains at most $K$ vertices.

Let ${\cal U}$ be an open cover of $\bdo X$ in $\cx$.
We construct an open cover $\cal V\prec \cal U$ of $\bdo X$ in $\cx$
consisting of interiors of standard neighborhoods such that the multiplicity of $\cal V$
is at most $K$.

Let $R'=2R+2D$. By the definition of topology (Proposition \ref{topo})  and by Lemma \ref{RR'},
for every $[\eta]\in \bdo X$ there exists a standard neighborhood
$U(\eta,N_{\eta},R')$ contained in some element of ${\cal U}$. By Lemma \ref{int} we have $[\eta]\in {\mr {int}}U(\eta,N_{\eta},R')$.
By compactness (Proposition \ref{cpt}), among such neighborhoods we can find a finite family
$\lk U(\eta^j,N_{\eta^j},R') \rk_{j=1}^{m}$ such that the family of
smaller standard neighborhoods $\lk U(\eta^j,N_{\eta^j},R) \rk_{j=1}^{m}$
covers $\bdo X$.
Let $N=\max \lk N_{\eta^1},N_{\eta^2}, \ldots ,N_{\eta^m} \rk$.
Let $A$ denote the set of vertices $v$ in $S_N(O)$ for which
there exists a good geodesic ray starting at $O$
and passing through $v$.
For each $v\in A$, pick some such good geodesic ray
$\xi^v=(w^v_0=O,w^v_1,w^v_2, \ldots ,w^v_N=v, \ldots )$.
We claim that the family
${\cal V}=\lk {\mr {int}}U(\xi^v,N,R)|\; v\in A \rk$
is as desired.

First we show that $\cal V$ covers $\bdo X$.
Let $\zeta=(z_0=O,z_1,z_2, \ldots )$ be an arbitrary
good geodesic ray.
Then $z_N=w^{z_N}_{N}$ and thus, by Lemma \ref{int},
$[\zeta]\in {\mr {int}}U(\xi_{z_N},N,R)$.


Now we show that $\cal V\prec \cal U$. To prove
this it is enough to show that for every $v\in A$ there
exists $j\in \lk 1,2, \ldots ,m \rk$ such that
$U(\xi^v,N,R)\subset U(\eta^j,N_{\eta^j},R')$.
Let $v\in A$.
Choose $j$ such that $[\xi^v] \in U(\eta^j,N_{\eta^j},R)$.
By Lemma \ref{U<U'}, to show that
$U(\xi^v,N,R)\subset U(\eta^j,N_{\eta^j},R')$ it is enough
to show that, for every $\zeta \in {\cal G}(\xi^v,N,R)$,
we have $\zeta \in {\cal G}(\eta^j,N_{\eta^j},R')$.
Let $\zeta=(z_0=O,z_1,z_2, \ldots ,z_N, \ldots )\in {\cal G}(\xi^v,N,R)$
and denote
$\eta^j=(v_0^j=O,v_1^j,v_2^j, \ldots )$. By Lemma \ref{repr}, we have
$|w^v_{N_{\eta_j}}v_{N_{\eta_j}}^j|\leq R+D$.
Then, by Corollary \ref{contr}, we have
\begin{align*}
|z_{N_{\eta_j}}v_{N_{\eta_j}}^j|& \leq
|z_{N_{\eta_j}}w^v_{N_{\eta_j}}|+|w^v_{N_{\eta_j}}v_{N_{\eta_j}}^j|\leq
(|z_Nw^v_{N}|+D)+(R+D)\leq
\\
& \leq 2R+2D=R'.
\end{align*}
Thus $\zeta \in {\cal G}(\eta^j,N_{\eta^j},R')$ and it follows that
$\cal V\prec \cal U$.

Finally, we claim that the multiplicity of $\cal V$ is at most $K$.
By Lemma \ref{disjoint}, if $|vv'|> 2R+D+2$ then
${\mr {int}}U(\xi^{v},N,R)\cap {\mr {int}}U(\xi^{v'},N,R)\subset
U(\xi^{v},N,R)\cap U(\xi^{v'},N,R)=\emptyset$.
Thus multiplicity of $\cal V$ is at most the number
of vertices in a ball of radius $2R+D+2$ in $X$, i.e. it
is at most $K$.
\kon

\section{The main result}
\label{ma}
The aim of this section is to prove the main result of the
paper---Theorem \ref{1} (Theorem \ref{main}).

The following result
will be crucial.
\begin{prop}[{\cite[Proposition 2.1]{BeMe}, \cite[Lemma 1.3]{Be}}]
\label{BM}
Let $(Y,Z)$ be a pair of finite--dimensional compact metrizable
spaces with $Z$ nowhere dense in $Y$, and such that
$Y\setminus Z$ is contractible and locally contractible
and the following condition holds:
\begin{itemize}
 \item For every $z\in Z$ and every open neighborhood $U$ of $z$
in $Y$, there exists an open neighborhood $V$ of $z$ contained
in $U$ such that $V\setminus Z\hookrightarrow U\setminus Z$
is null--homotopic.
\end{itemize}
Then $Y$ is an ER and $Z$ is a $Z$--set in $Y$.
\end{prop}

Before proving Theorem \ref{1} we need an important preparatory lemma.

\begin{lem}
\label{BMsys}
Let $[\eta]\in \bdo X$ and let $U(\eta,N,R)$ be a standard
neighborhood of $[\eta]$ in $\cx$. Then
there exists $N'$ such that
$U(\eta,N',R)\subset U(\eta,N,R)$ and the inclusion map
$U(\eta,N',R)\cap X\hookrightarrow U(\eta,N,R)\cap X$ is
null--homotopic.
\end{lem}

\dow
Let $R'=4D+7$. By Lemma \ref{RR'}, there exists $\widetilde N$ such that $U(\eta,\widetilde N,R')\subset U(\eta,N,R)$, so that it is enough
to prove the following. For natural $R\geq D$ there exists $N'$ such that
$U(\eta,N',R)\subset U(\eta,N,R')$ and the inclusion map
$U(\eta,N',R)\cap X\hookrightarrow U(\eta,N,R')\cap X$ is
null--homotopic.
\medskip

Before we start, let us give a rough idea of the proof. Let us restrict to the problem of contracting loops from $U(\eta,N',R)\cap X$ in $U(\eta,N,R')\cap X$ (this turns out to be the most complicated case).
Let $\alpha$ be such a loop. We connect each vertex of $\alpha$ by a good geodesic with $O$, and we are interested in the vertex of this geodesic at certain distance $M$ from $O$, where $N<M<N'$. All vertices constructed in this way lie in a certain ball (see Condition 1 below), which is in turn contained in $U(\eta,N,R')\cap X$ (see Condition 3 below). If we connect these vertices by $1$--skeleton geodesics in the right order we obtain a loop $\alpha_M$, which lies in the ball considered (Corollary \ref{balls convex}) and is contractible inside this ball (Corollary \ref{balls contractible}). So we need to find a free homotopy between $\alpha$ and $\alpha_M$, which we construct via intermediate loops $\alpha_l$. To find that two such consecutive loops are homotopic in $U(\eta,N,R')\cap X$, we need Condition 2. This condition guarantees that all relatively small loops by which consecutive $\alpha_l$ differ can be contracted inside $U(\eta,N,R')\cap X$.
\medskip

Let $M=N+R+1$ and $N'-1\geq (R+D+4)M$. We will show that $N'$ is as desired. Denote $\eta=(v_0=O,v_1,v_2, \ldots )$.
The choice of $M$ and $N'$ guarantees that the following three conditions hold.
\medskip

\textbf{Condition 1.}
Let $\xi=(w_0=O,w_1, \ldots ,w_{k})$ be a good geodesic with $k\geq N'-1$ and $w_k\in \overline{U(\eta,N',R)\cap X}$. Then $w_M\in B_{D+1}(v_M)$.
\medskip

Indeed, let $(z_0=O,z_1,\ldots, z_l)\in {\cal G}(\eta,N',R)$ be such
that $|w_kz_l|\leq 1$ (guaranteed by definition of $U(\eta,N',R)$).
Since $k\geq N'-1$, we have, by Corollary \ref{contr'}, that
\begin{align*}
|w_{N'-1}v_{N'-1}|& \leq |w_{N'-1}z_{N'-1}|+|z_{N'-1}v_{N'-1}|\leq
\\
& \leq
(2|w_kz_l|+D)+(1+|z_{N'}v_{N'}|+1)\leq R+D+4.
\end{align*}
Thus, by Corollary \ref{contr},
we have
\begin{align*}
|w_Mv_M|\leq \frac{M}{N'-1}|w_{N'-1}v_{N'-1}|+D
\leq \frac{1}{R+D+4}|w_{N'-1}v_{N'-1}|+D\leq D+1.
\end{align*}
\medskip

\textbf{Condition 2.}
Let $\xi=(w_0=O,w_1, \ldots ,w_{k})$ be as in Condition 1.
Then, for every $k\geq l\geq M+1$ we have $B_{D+3}(w_l)\subset U(\eta,N,R')\cap X$.
\medskip

To show this observe that, as in the proof of the previous condition,
we have $|w_{N'-1}v_{N'-1}|\leq R+D+4$.
Now, let $z$ be a vertex of $B_{D+3}(w_l)$. Choose a good geodesic
$(z_0=O,z_1,z_2, \ldots ,z_m=z)$ (guaranteed by Corollary \ref{goodex}).
Since $l\geq M+1=N+R+2\geq N+(D+3)$, we have that $m\geq N$ and $z_N$ is defined. Thus,
by Corollary \ref{contr} and Corollary \ref{contr'}, we have
\begin{align*}
|z_Nv_N|& \leq |z_Nw_N|+|w_Nv_N|\leq (2|z_mw_l|+D)+\Big(\frac{N}{N'-1}|w_{N'-1}v_{N'-1}|+D\Big)<
\\
& < (2(D+3)+D)+\Big(\frac{1}{R+D+4}(R+D+4)+D\Big)=4D+7=R'.
\end{align*}
Thus $z\in U(\eta,N,R')\cap X$ and it follows that $B_{D+3}(w_l)\subset U(\eta,N,R')\cap X$.
\medskip

\textbf{Condition 3.}
We have $B_{D+1}(v_M)\subset
U(\eta,N,R')\cap X$.
\medskip

This follows immediately from Condition 2, but we want to record it separately.
\medskip

\textbf{The goal.} First observe that $U(\eta,N',R)\subset U(\eta,N,R')$ by
Lemma \ref{RR'} and the definition of $N'$.
Now we show that the map
$\pi_i(U(\eta,N',R)\cap X)\to \pi_i(U(\eta,N,R')\cap X)$
induced by inclusion is trivial, for every $i=0,1,2, \ldots $.
Let $A$ be the smallest full subcomplex of $X$ containing
$U(\eta,N',R)\cap X$.
Observe that the vertices of $A$ lie in
$\overline{U(\eta,N',R)\cap X}$.
By Condition 2, $A$ is contained in $U(\eta,N,R')\cap X$.
Thus it is enough to show that the map
$\pi_i(A)\to \pi_i(U(\eta,N,R')\cap X)$ induced by the inclusion
is trivial and we may restrict ourselves
only to simplicial spherical cycles.
\medskip

\textbf{Case ($i=0$).} Let $z^1,z^2$ be two vertices of $A$.
We will construct a simplicial path in $U(\eta,N,R')\cap X$
connecting $z^1$ and $z^2$.

Choose (using Corollary \ref{goodex})
good geodesics $(z^j_0=O,z^j_1, \ldots ,z^j_{k(j)}=z^j)$, $j=1,2$.
By Condition 2, $(z^j_M,z^j_{M+1}, \ldots ,z^j_{k(j)}=z^j)$ is contained in
$U(\eta,N,R')$ and, by Condition 1, we have $z^j_M\in B_{D+1}(v_M)$.
Choose a $1$--skeleton
geodesic
$(u_1=z^1_M,u_2, \ldots ,u_l=z^2_M)$.
Since balls are geodesically convex (Corollary \ref{balls convex}), this geodesic
is contained in $B_{D+1}(v_M)$ and hence, by Condition 3, it is
contained in $U(\eta,N,R')\cap X$.

Then the $1$--skeleton path
$(z^1=z^1_{k(1)},z^1_{k(1)-1}, \ldots ,z^1_M=u_1,u_2, \ldots ,u_l=z^2_M,
z^2_{M+1}, \ldots ,z^2_{k(2)}=z^2)$ connects $z^1$ and $z^2$
and is contained in $U(\eta,N,R')\cap X$. Therefore the map
$\pi_0(A)\to \pi_0(U(\eta,N,R')\cap X)$ is trivial.
\medskip

\textbf{Case ($i=1$).} Let $\alpha=(z^0,z^1, \ldots ,z^n=z^0)$ be a $1$--skeleton loop  in $A$. We show that this loop can be contracted
within $U(\eta,N,R')\cap X$.

Choose
good geodesics $(z^j_0=O,z^j_1, \ldots ,z^j_{k(j)}=z^j)$ (guaranteed by Corollary \ref{goodex}), for $j=0,1,2, \ldots ,n-1$.
By $z^j_k$, for $k>k(j)$, we denote $z^j$.
Let $K=\max \lk k(0),k(1), \ldots , k(n-1) \rk$.
Observe that, by Corollary
\ref{contr'}, we have $|z^j_lz^{j+1}_l|\leq D+2$ (we consider $j$ modulo $n$), for every $l=M,M+1,M+2, \ldots ,K$ (we are not interested in smaller $l$).
For these $l$ let $(z^j_l=t^{j,0}_l,t^{j,1}_l, \ldots ,t^{j,p_l(j)}_l=z^{j+1}_l)$ be
arbitrary
$1$--skeleton geodesics. Record that $p_l(j)\leq D+2$.

Thus, for every $l=M+1,M+2,\ldots,K$ and for every $j=0,1, \ldots ,n-1$, we have a $1$--skeleton loop
\begin{align*}
\gamma^j_l=(z^j_l,z^j_{l-1}=t^{j,0}_{l-1},t^{j,1}_{l-1},\ldots, t^{j,p_{l-1}(j)}_{l-1}&=z^{j+1}_{l-1},
\\
z^{j+1}_l&=t^{j,p_l(j)},
t^{j,p_l(j)-1}_l,\ldots,t^{j,0}=z^j_l)
\end{align*}
of length at most
$$1+p_{l-1}(j)+1+p_l(j)\leq 1+(D+2)+1+(D+2)=2D+6.$$
Hence $\gamma^j_l\subset B_{D+3}(z^j_l)$. Since balls are contractible (Corollary \ref{balls contractible}), $\gamma^j_l$ is contractible inside $B_{D+3}(z^j_l)$, which is, by Condition 2, contained in $U(\eta,N,R')$.
Thus, for $M\leq l\leq K$, the loops
\begin{align*}
\alpha_l =(z^0_{l}=t^{0,0}_l,t^{0,1}_l,&\ldots ,t^{0,p_l(0)}_l=z^1_l=t^{1,0}_l,t^{1,1}_l, \ldots, t^{1,p_l(1)}_l=z^2_l,\ldots
\\
&\ldots ,z^{n-1}_l=t^{n-1,0}_l,t^{n-1,1}_l,\ldots, t^{n-1,p_l(n-1)}_{l}=z^n_{l}=z^0_l)
\end{align*}
for consecutive $l$ are freely homotopic in $U(\eta,N,R')$.

Observe that $\alpha=\alpha_K$. On the other hand $\alpha_M\subset B_{D+1}(v_M)$, by Condition 1 and by geodesic convexity of balls (Corollary \ref{balls convex}).
Moreover, since balls are contractible (Corollary \ref{balls contractible}),
$\alpha_M$ can be contracted inside $B_{D+1}(v_M)$, which lies in $U(\eta,N,R')$, by Condition 3. Thus $\alpha$ is contractible in $U(\eta,N,R')$.
It follows that the map
$\pi_1(A)\to \pi_1(U(\eta,N,R)\cap X)$ is trivial.
\medskip

\textbf{Case ($i>1$).} Since $A$ is a full subcomplex of a systolic complex it is,
by Lemma \ref{asphericity of full subcomplexes},
aspherical and thus $\pi_i(A)=0$ and the map in question is obviously trivial.
\kon

\begin{tw}[Theorem \ref{1}]
\label{main}
Let a group $G$ act geometrically by simplicial automorphisms
on a systolic complex $X$.
Then $\cx=X\cup \bdo X$, where $O$ is a vertex of $X$,
 is a compactification of $X$ satisfying the
following:
\begin{enumerate}
 \item $\cx$ is a Euclidean retract (ER),
 \item $\bdo X$ is a $Z$--set in $\cx$,
 \item for every compact set $K\subset X$,
$(gK)_{g\in G}$ is a null sequence,
 \item the action of $G$ on $X$ extends to an action,
by homeomorphisms, of $G$ on $\cx$.
\end{enumerate}
\end{tw}
\dow
{\sl (1. and 2.)} By Corollary \ref{metr}, Proposition \ref{cpt}, and Proposition \ref{findim}, we have that
$\cx=X\cup \bdo X$ is a finitely dimensional metrizable compact space.

Since $X$ is a simplicial complex, it is locally contractible and, by Theorem \ref{con}, it is contractible since it is
a systolic complex. By the definition of the topology
on $\cx$ (c.f. Proposition \ref{topo}), it is clear that $\bdo X$ is nowhere dense in $\cx$.
Thus we are in a position to apply Proposition \ref{BM}.
Let $x\in \bdo X$ and let $U$ be its open neighborhood in $\cx$.

By definition of the topology (Proposition \ref{topo})
we can find a standard neighborhood $U(\eta,N,R)\subset U$, where
$[\eta]=x$.
By Lemma \ref{BMsys}, there exists a standard
neighborhood $U(\eta,N',R)\subset U(\eta,N,R)\subset U$
(with $[\eta] \in {\mr {int}}(U(\eta,N',R)$, by Lemma \ref{int}) such that the map
${\mr {int}}(U(\eta,N',R)\cap X)\hookrightarrow U(\eta,N',R)\cap X
\hookrightarrow U(\eta,N,R)\cap X\hookrightarrow U\cap X$
is null--homotopic.
Thus $\cx$ is an ER and $\bdo X$ is a $Z$--set in $\cx$.
\medskip

{\sl (3.)}
Let ${\cal U}$ be an open cover of $\cx$ and let $K\subset X$ be a compact set. We will show that all but finitely many translates $gK$, for $g\in G$,
are ${\cal U}$--small.

Let $R>D$ be such that $K\subset B_R(z)$, for some vertex $z$.
As in the proof of Proposition \ref{findim}, we can find
a natural number $N$,
a finite set of vertices $A\subset S_N(O)$ and a collection
of good geodesic rays $\lk \xi^v \; |\; \; v\in A \rk$ with $\xi^v$ passing through $v$ such that
the following holds. The
family ${\cal V}=\lk {\mr {int}}U(\xi^v,N,R)\;| \; \; v\in A \rk$
covers $\bdo X$
and the family ${\cal V}'=\lk U(\xi^v,N,4R)\;| \; \; v\in A \rk$
is a refinement of ${\cal U}$.
Thus we can find an open cover ${\cal W}={\cal V}\cup {\cal W'}$
of $\cx$ such that every $W\in {\cal W'}$ is contained
in $X$. By compactness---Proposition \ref{cpt}---there is a finite
subfamily of ${\cal W}$ covering $\cx$. It follows that there
exists a natural number $N'>N$ such that
$\cx \setminus B_{N'}(O)\subset \bigcup {\cal V}$.
By properness of the action there exists a cofinite subset
$H\subset G$ such that
$gK\subset B_R(gz)\subset X\setminus B_{N'}(O)$, for $g\in H$.

We claim that, for every $g\in H$, we have
$gK\subset B_R(gz)\subset U(\xi^v,N,4R)\cap X$, for some $v\in A$.
Assertion {\sl (3.)} follows then from the claim.
Let $g\in H$.
Since $\cx \setminus B_{N'}(O)\subset \bigcup {\cal V}$,
there exists $v\in A$ such that $gz\in {\mr {int}}
U(\xi^v,N,R)$.
We show that $B_R(gz)\subset U(\xi^v,N,4R)$.
Let $x\in B_R(gz)$ and let
$\zeta=(z_0'=O,z_1', \ldots ,z_l')$ be a good
geodesic (which exists by Corollary \ref{goodex}) such that $z_l'\in B_R(gz)$ is a vertex of the simplex containing $x$ in its interior.
Since $gz\in U(\xi^v,N,R)$ there exists a good geodesic
$(z_0=O,z_1,z_2, \ldots ,z_k=gz)$, such that
$|z_Nv|\leq R$.
We have $l,k\geq N'$ and $|z_l'z_k|\leq R$.
Hence, by Corollary \ref{contr'}, we have
\begin{align*}
|z_N'v|& \leq |z_N'z_N|+|z_Nv|\leq
(2|z_l'z_k|+D)+|z_Nv|\leq
\\
& \leq(2R+D)+R< 4R.
\end{align*}
Thus $\zeta \in {\cal G}(\xi^v,N,4R)$ and hence
$x\in U(\xi^v,N,4R)$.
It follows that $B_R(gz)\subset U(\xi^v,N,4R)$.
Since $g\in H$ was arbitrary we have that elements of
$(gK)_{g\in H}$ are ${\cal V}'$--small and thus they
are ${\cal U}$--small.
\medskip

{\sl (4.)}
For $g\in G$ we define a map $g\circ \colon X\cup \bdo X \to X\cup \partial_{gO}X$ as follows. For $x\in X$ let $g\circ x=gx$ and for
$x=[(v_0=O,v_1,v_2, \ldots )]\in \bdo X$ let $g\circ x=[(gv_0=gO,gv_1,gv_2, \ldots )]$. This is obviously a well defined homeomorphism.

We extend the action of $G$ on $X$ to $X\cup \bdo X$ by the formula
$g\cdot x= \Phi_{gOO} (g\circ x)$, for $x\in \bdo X$.
By Lemma \ref{OO''}, the map $g\ \cdot \ \colon X\cup \bdo X \to
X\cup \bdo X$ is a homeomorphism.
To see that $(gh)\cdot x=g\cdot (h\cdot x)$, for $x\in \bdo X$, pick some representative $\eta=(v_0=O,v_1,\ldots)$ of $x$.
We need to show that $$\Phi_{ghOO}(gh\circ [\eta])=\Phi_{gOO}(g\circ \Phi_{hOO}(h\circ [\eta])).$$
Recall that, by Lemma \ref{vv'}, mappings $\Phi_{gOO},\Phi_{hOO}$ and $\Phi_{ghOO}$
displace representative rays by a finite Hausdorff distance.
Hence $\Phi_{ghOO}(gh\circ [\eta])$ is the class of rays starting at $O$ at a finite Hausdorff distance from $(ghv_0=ghO, ghv_1,\ldots)$. On the other hand, $\Phi_{hOO}(h\circ [\eta])$ is the class of rays starting at $O$ at a finite Hausdorff distance from $(hv_0=hO, hv_1,\ldots)$, hence $g\circ \Phi_{hOO}(h\circ [\eta])$ as well as $\Phi_{gOO}(g\circ \Phi_{hOO}(h\circ [\eta]))$ is the class of rays (starting at, respectively, $gO$ and $O$) at a finite Hausdorff distance from
$(ghv_0=ghO, ghv_1,\ldots)$.
This proves the desired equality.


Hence we get an extension of the action of $G$ on $X$ to an action
on $\cx$ by homeomorphisms.
\kon


\section{Flat surfaces}
\label{Flat surfaces}
With this section we start the second part of the article, in which we define Euclidean geodesics, establish Theorem \ref{2} and Theorem \ref{3}. Before we define Euclidean geodesics, we first need to study, as mentioned in Section \ref{introduction}, the minimal surface spanned on a pair of directed geodesics connecting given vertices. The tools for this are minimal surfaces (Section \ref{Flat surfaces}) and layers (Section \ref{Layers}).

In this section we recall some definitions and facts concerning flat
minimal surfaces in systolic complexes proved by Elsner
\cite{E},\cite{E3}.

\begin{defin}
\label{def_of_disc}
The \emph{flat systolic plane} is a systolic 2--complex obtained
by equilaterally triangulating Euclidean plane. We denote it by
$\mathbb{E}_\Delta^2$. A \emph{systolic disc} is a systolic
triangulation of a 2--disc and a \emph{flat disc} is any systolic
disc $\Delta$, which can be embedded into $\mathbb{E}_\Delta^2$,
such that $\Delta^{(1)}$ is embedded isometrically into
1--skeleton of $\mathbb{E}_\Delta^2$. A systolic disc $\Delta$ is
called \emph{wide} if $\partial\Delta$ is a full subcomplex of
$\Delta$. For any vertex $v\in \Delta^{(0)}$ the \emph{defect} at $v$
(denoted by def$(v)$) is $6-t(v)$ for $v\notin
\partial \Delta^{(0)}$, and $3-t(v)$ for $v \in
\partial\Delta^{(0)}$, where $t(v)$ is the number of triangles in
$\Delta$ containing $v$. It is clear that internal vertices of a
systolic disc have nonpositive defects.
\end{defin}

We will need the following easy and well known fact.

\begin{lemma}[Gauss-Bonnet Lemma]
\label{Gauss-Bonnet}
If $\Delta$ is any triangulation of a 2--disc, then
$$
\sum_{v\in \Delta^{(0)}} \mathrm{def} (v)=6
$$
\end{lemma}

Flat systolic discs can be characterized as follows.

\begin{lemma}[\cite{E}, Lemma 2.5]
\label{characterization_of_flat}
A systolic disc $D$ is flat if and only if it satisfies the
following three conditions: \item{(i)}   $D$ has no internal
vertices of defect $<0$ \item{(ii)}  $D$ has no boundary vertices
of defect $<-1$ \item{(iii)} any segment in $\partial D$
connecting vertices with defect $<0$ contains a vertex of defect
$>0$.
\end{lemma}

Now we recall another handful of definitions.

\begin{defin}
\label{def_of_surface}
Let $X$ be a systolic complex. Any simplicial map $S\colon \Delta
\rightarrow X$, where $\Delta$ is a triangulation of a 2--disc, is
called a \emph{surface}. We say that $S$ is \emph{spanned} on a loop $\gamma$, if $S|_{\partial
\Delta}=\gamma$. A loop $\gamma$ is \emph{triangulable}, if there exists a
surface $S$ spanned on $\gamma$, such that all the vertices of $\Delta$ are in
$\partial \Delta$. A surface $S$ is \emph{systolic, flat} or
\emph{wide} if the disc $\Delta$ satisfies the corresponding
property. If $S$ is injective on $\partial\Delta$ and minimal (the
smallest number of triangles in $\Delta$) among surfaces with the
given image of $\partial \Delta$, then $S$ is called
$\emph{minimal}$. A geodesic in $\Delta^{(1)}$ is called
$\emph{neat}$ if it stays out of $\partial\Delta$ except possibly
at its endpoints. A surface $S$ is called \emph{almost geodesic} if it
maps neat geodesics in $\Delta^{(1)}$ isometrically into
$X^{(1)}$.
\end{defin}

The following is part of the main theorem of \cite{E}.

\begin{theorem}[\cite{E}, Theorem 3.1]
\label{main_Elsner_theorem}
Let $X$ be a systolic complex. If $S$ is a wide flat minimal
surface in $X$ then $S$ is almost geodesic.
\end{theorem}

We will also use the following handy fact, whose proof can be
extracted from \cite{E3}. In case $\gamma$ has length 2 it follows
immediately from 6--largeness.

\begin{proposition}[\cite{E3}, Proposition 3.10]
\label{stable_geodesics_in characteristic image}
Let $X$ be a systolic complex and $S\colon \Delta\rightarrow X$ a wide flat minimal
surface. Let $\gamma$ be a neat 1--skeleton geodesic in $\Delta\subset \mathbb{E}^2_{\Delta}$,
which is contained in a straight line. Then for any 1--skeleton geodesic $\overline{\gamma}$ in $X$ with the same endpoints as $S(\gamma)$ there is another minimal surface $S'\colon \Delta\rightarrow X$ such that $S'(\gamma)=\overline{\gamma}$ and $S=S'$ on the vertices of $\Delta$ outside $\gamma$.
\end{proposition}

\section{Layers}
\label{Layers}
In this section we introduce and study the notion of layers for a pair of convex subcomplexes of a systolic complex.
If those subcomplexes are vertices $v,w$, then the \emph{layer $k$} is the span of all vertices, in 1--skeleton geodesics $vw$,
at distance $k$ from $v$ (c.f. Definition \ref{layer}). In particular, simplices of the directed geodesics between $v$ and $w$ (c.f. Definition \ref{directed geodesic}), as well as the simplices of Euclidean geodesics (which we construct in Section \ref{Euclidean geodesic})
lie in appropriate layers.

On the other hand, layers in systolic complexes seem to be interesting on their own.

\begin{defin}
\label{layer}
Let $V,W$ be convex subcomplexes of a systolic complex $X$ and
$n=|V,W|$. For $i=0,1,\ldots,n$ we define the \emph{layer $i$}
between $V$ and $W$ as the subcomplex of $X$ equal to $B_i(V)\cap
B_{n-i}(W)$. We will denote it by $L_i(V,W)$ (or shortly $L_i$ when $V,W$ are understood).
\end{defin}

\begin{rem}
\label{layers convex}
$L_i$ are convex, since they are intersections of convex $B_i(V),B_{n-i}(W)$ (see remarks after Definition \ref{convex}).
\end{rem}

\begin{lemma}
\label{layers close}
\item(i) $L_i=S_i(V)\cap S_{n-i}(W)$, for $0\leq i\leq n$.
\item(ii) $L_{j}\subset S_{j-i}(L_i)$, for $0\leq i<j\leq n$. In particular $L_{i+1}\subset S_1(L_i)$, for $0\leq  i< n$.
\end{lemma}
\proof
(i) W.l.o.g we only need to prove that $L_i\subset S_i(V)$. Take a vertex $x\in L_i$. Then we have $|x,V|\leq i$ and $|x,W|\leq n-i$, while $|V,W|=n$. Thus by the triangle inequality we have $|x,V|=i$, as desired.

(ii) By (i) we have that $B_{j-i-1}(L_i)\cap L_j=\emptyset$, thus we only need to prove that $L_j\subset B_{j-i}(L_i)$. Let $x$ be a vertex in $L_j$. Since, by (i), we have $x\in S_j(V)$, there is a vertex  $y\in B_i(V)$ at distance $j-i$ from $x$. Since $x\in B_{n-j}(W)$, we have $y\in B_{n-i}(W)$. Thus $y\in L_i$ and $x\in B_{j-i}(L_i)$.
\qed
\medskip

Now we study the properties of layers.

\begin{lemma}
\label{layers_are_infinity_large}
For $0<i<n$ we have that $L_i$ is $\infty$--large.
\end{lemma}
\proof Suppose the layer $L_i$ is not $\infty$--large. Then there exists an embedded cycle $\Gamma$ in
 $L_i$ (denote its consecutive vertices by $p_1,p_2,\dots,
p_k,p_1,k\geq 4$) which is a full subcomplex of $X$.

Denote $D_1=\text{span} \{B_{i-1}(V),\Gamma\},D_2=\text{span}
\{B_{n-i-1}(W),\Gamma \}$. We have that $D_1\cap D_2=\Gamma$. Notice that $D_1\cup D_2$
is a full subcomplex of $X$, because there are no edges in $X$ between vertices
in $B_{i-1}(V)$ and vertices in $B_{n-i-1}(W)$.

Observe that $\Gamma$ is
contractible in $D_1$ (and similarly in $D_2$). Indeed, by Lemma \ref{layers close}(i) we have that $\Gamma\subset S_i(V)$. Thus we can project the edges of $\Gamma$ onto $B_{i-1}(V)$ (c.f. Definition \ref{projection}). If we choose a vertex in each of these projections, we get, by Lemma \ref{projection lemma}, that these vertices form a loop. This loop is homotopic to $\Gamma$ in $D_1$. Moreover, since $B_{i-1}(V)$ is contractible (by remarks after Definition \ref{convex}) it follows that $\Gamma$ is contractible in $D_1$ (and similarly in $D_2$), as desired.
The simplicial
sphere $S$ formed of these two contractions is contractible in $D_1\cup D_2$ as
full subcomplexes of $X$ are aspherical (Lemma \ref{asphericity of full subcomplexes}).

Now use Meyer-Vietoris sequence of the pair $D_1,D_2$. Since $[\Gamma]$ is the
image of $[S]=0$ under $H_2(D_1\cup D_2)\rightarrow H_1(D_1\cap D_2)$ it
follows that the cycle $\Gamma$ is homological to zero in itself. This is a
contradiction.\qed

\begin{lemma}
\label{no_trapezoid_general}
Let $\sigma_1, \sigma_2, \sigma_3$ be maximal simplices in the layer $L_i$ for some $0\leq i\leq n$
and $\tau_1=\sigma_1\cap
\sigma_2,\ \tau_2=\sigma_2\cap \sigma_3$. Then
$\tau_1\cap\tau_2=\emptyset$ or $\tau_1\subset\tau_2$ or
$\tau_2\subset\tau_1$.
\end{lemma}

\proof W.l.o.g. assume that $i\neq 0$. Suppose the
lemma is false. Then there exist vertices $p_1\in \tau_1\backslash\tau_2,\
p_2\in \tau_2\backslash\tau_1, \ r\in \tau_1\cap\tau_2$. By Lemma \ref{layers close}(ii) we have that $\sigma_1,\sigma_3\subset S_1(L_{i-1})$. Denote by $q_1, q_2$
some vertices in the projections (c.f. Definition \ref{projection}) of $\sigma_1, \sigma_3$
onto $L_{i-1}$. We have $|q_1q_2|\leq 1$, because both $q_1$ and $q_2$ are neighbors of
$r$ and the projection of $r\in L_i\subset S_1(L_{i-1})$ (c.f. Lemma \ref{layers close}(ii)) onto $L_{i-1}$ is a simplex (Lemma \ref{projection lemma}). Now we will argue that we can assume that $q_1p_2$ is an edge. In
case $q_1\neq q_2$ consider the 4--cycle $q_1q_2p_2p_1q_1$. It must have a
diagonal. We can then assume w.l.o.g. that $q_1p_2$ is an edge. In case $q_1=q_2$ we also have that $q_1p_2$ is an edge.
In both cases it follows that $p_2$ belongs to the simplex which is the projection of $q_1\in L_{i-1}\subset S_1(L_i)$ (c.f. Lemma \ref{layers close}(ii)) onto $L_i$. This simplex also contains $\sigma_1$. But $p_2\notin\sigma_1$, which
contradicts the maximality of $\sigma_1$. \qed

\begin{cor}
\label{no_trapezoid}
Let $T$ be the following simplicial complex: the trapezoid build
of three triangles $p_1rs_1, p_1rp_2,p_2rs_2$. Then there is no
isometric embedding of $T^{(1)}$ into $L_i^{(1)}$, for $0\leq i\leq n$.
\end{cor}

\proof Extend the images of those three triangles to maximal simplices
$\sigma_1, \sigma_2, \sigma_3$ and apply Lemma \ref{no_trapezoid_general}.\qed

\begin{cor}
\label{corollary_to_no_trapezoid}
Let $0<i<n$. Let $|p_0r_0|\leq 1, |p_dr_d|\leq 1$ for vertices $p_0,r_0,p_d,r_d\in L_i$ such
that $|p_0p_d|=|r_0r_d|=d\geq 2$ and $|p_0r_d|\geq d, |r_0p_d|\geq d$. Then, for
any 1--skeleton geodesics $(p_i),(r_i),0\leq i\leq d$ connecting $p_0$ with
$p_d$ and $r_0$ with $r_d$, respectively, and for any $0\leq i,j\leq d$ such that
$|i-j|\leq 1$, we have that $|p_ir_j|\leq 1$ (i.e. $p_ir_j$ is an edge or $p_i=r_j$).
\end{cor}

\proof We will prove the corollary by induction on $d$. First observe that since $L_i$ is $\infty$--large (Lemma \ref{layers_are_infinity_large}), the loop $p_0p_1\ldots p_dr_d\ldots
r_1r_0p_0$ is triangulable and there exists a diagonal cutting off a triangle.
There are only four possibilities for this diagonal and we can w.l.o.g suppose
this diagonal is $p_0r_1$. Now since $p_0\in S_d(r_dp_d)$ and both $p_1$ and
$r_1$ lie in the projection of $p_0$ onto $B_{d-1}(r_dp_d)$, then by Lemma \ref{projection lemma} either $p_1r_1$ is
an edge or $p_1=r_1$.

Now we start the induction. If $d=2$ and the loop $p_1r_1r_2p_2p_1$ is embedded, then it has
a diagonal. The rest of the required inequalities follows from applying twice
Corollary \ref{no_trapezoid}.

Suppose that for $d-1$ the
corollary is already proved. Then applying it to the pair $p_1r_1,p_dr_d$
yields all the required inequalities except for the estimate on $|r_0p_1|$. But this
follows from Corollary \ref{no_trapezoid} applied to the trapezoid $r_0p_0r_1p_1p_2$. \qed

\begin{cor}
\label{corollary_2_with_leq_to_no_trapezoid}
If $pr,p'r'$ are edges in $L_i$, for some $0<i<n$, such that $|pp'|=|rr'|=d\geq 2$ and $|pr'|\leq d, |p'r|\leq d$, then
$|pr'|=|p'r|=d$.
\end{cor}
\proof
By contradiction.

\textbf{Case $|pr'|=|p'r|=d-1$.} If $d>2$ (if $d=2$ there is a diagonal in the square $pr'p'rp$) then
Corollary \ref{corollary_to_no_trapezoid} applied to $d-1$ in place of $d$,\ $p_0=p\ ,p_{d-1}=r',\ r_0=r,\ r_{d-1}=p'$ gives $|pp'|=|rr'|=d-1$, contradiction.

\textbf{Case $|pr'|=d-1,\ |p'r|=d$.} Again apply Corollary \ref{corollary_to_no_trapezoid}, with $p_0=p,\ r_0=r,\ p_d=r_d=p',\ p_{d-1}=r'$, getting $|rr'|=d-1$, contradiction.
\qed
\medskip

Below we present another important property of layers. Since it will not be needed in the article, we do not include the proof. Denote $L=\Span(L_i\cup L_{i+1})$ for some $1\leq i<n-1$.

\begin{lemma}
\label{union of two layers is infinity-large}
$L$ is $\infty$--large.
\end{lemma}

We end with a simple, but useful observation.

\begin{lemma}
\label{difference1}
For any edges $vw,xy$ such that $v,x\in L_i, \ w,y\in L_{i+1}$, where $0\leq i< n$, we have that $||vx|-|wy||\leq 1$.
\end{lemma}
\proof
By contradiction. Suppose, w.l.o.g., that $|wy|=2+|vx|$. Hence $v$ lies on a 1--skeleton geodesic $wy$. Thus, by convexity of layers (Remark \ref{layers convex}) and by Proposition
\ref{second convexity}, we have that $v$ lies in $L_{i+1}$, which is, by Lemma \ref{layers close}, disjoint with $L_i$, contradiction.
\qed

\section{Euclidean geodesics}
\label{Euclidean geodesic}

In this section we define, for a pair of simplices $\sigma,\tau$ as below,
a sequence of simplices in the layers between  $\sigma$ and $\tau$, which can be considered as a "Euclidean" geodesic
joining $\sigma$ and $\tau$. Unlike the directed geodesics defined by Januszkiewicz and \'Swi\k{a}tkowski (see Definition \ref{directed geodesic}), Euclidean geodesics are symmetric with respect to
$\sigma$ and $\tau$.

The definition requires a lenghty preparation. Roughly speaking, we start by spanning a minimal surface between directed geodesics from $\sigma$ to $\tau$ and from $\tau$ to $\sigma$. We observe that this surface is
flat whenever the two directed geodesics are far apart (we call the corresponding layers \emph{thick}). Next we show that this "piecewise" flat surface is in some sense unique. This occupies the first part of the section, up to Definition \ref{characteristic image}. Then we look at the geodesics in the Euclidean metric in the flat pieces and use them to define \emph{Euclidean geodesics} in systolic complexes, c.f. Definition \ref{euclidean geodesic}. Finally, we establish some of their basic properties.

The setting, which we fix for Sections \ref{Euclidean geodesic}---\ref{Contracting property}
is the following. Let $\sigma, \tau$ be simplices of a systolic complex $X$, such that for some natural $n\geq 0$ we have $\sigma \subset S_n(\tau), \tau \subset S_n(\sigma)$. Let
$\sigma_0\subset\sigma, \sigma_1,\ldots, \sigma_n\subset\tau$ and $\tau_n\subset\tau, \tau_{n-1}, \dots,\tau_0\subset\sigma$ be
sequences of simplices in $X$, such that for $0\leq k<n$ we have that $\sigma_k, \sigma_{k+1}$ span a simplex and
$\tau_k,\tau_{k+1}$ span a simplex. In particular, $\sigma_k,\tau_k$ lie in the layer $k$ between $\sigma$ and $\tau$ (c.f. Definition \ref{layer}).

Note that if $\sigma_0=\sigma, \sigma_1,\ldots, \sigma_n\subset\tau$ is the directed geodesic
from $\sigma$ to $\tau$ and $\tau_n=\tau, \tau_{n-1}, \dots,\tau_0\subset\sigma$ is
the directed geodesic from $\tau$ to $\sigma$ (c.f. Definition \ref{directed geodesic}), then the above condition is satisfied. This special choice of $(\sigma_k),(\tau_k)$ will be very important later and we will frequently distinguish it.

\begin{defin}
\label{thickness}
For $0\leq i\leq n$ the \emph{thickness} of the layer $i$ for $(\sigma_k),(\tau_k)$ is
the maximal distance between vertices in $\sigma_i$ and in $\tau_i$.
If the layer $i$ for $(\sigma_k),(\tau_k)$ has thickness
$\leq 1$ we say that the layer $i$ for
$(\sigma_k),(\tau_k)$ is \emph{thin}, otherwise we say that
the layer $i$ for $(\sigma_k),(\tau_k)$ is \emph{thick}. If $(\sigma_k),(\tau_k)$ are directed geodesics from $\sigma$ to $\tau$ and from $\tau$ to $\sigma$, respectively, then
we skip "for $(\sigma_k),(\tau_k)$" for simplicity.
\end{defin}

\noindent \textbf{Caution.} Perhaps, to avoid confusion, we should not have used the word "layer" in the above definition, since we are in fact only checking the position of $\sigma_i$ w.r.t. $\tau_i$. Even if the layer $i$ between $\sigma$ and $\tau$ is large, it can happen that the thickness of the layer $i$ for $(\sigma_k),(\tau_k)$ is small.
However, we decided that this terminology suits well our approach, in which we will be mostly interested in the part of the layer $i$ between $\sigma$ and $\tau$, which lies between $\sigma_i$ and $\tau_i$.
\medskip

\begin{defin}
\label{thick_interval}
A pair $(i,j)$, where $0\leq i<j\leq n$ is called a \emph{thick interval} (for $(\sigma_k),(\tau_k)$) if
the layers $i$ and $j$ (for $(\sigma_k),(\tau_k)$) are thin, $i+1<j$, and for every $l$, such that
$i<l<j$, the layer $l$ (for $(\sigma_k),(\tau_k)$) is thick. We say that the thick interval $(i,j)$ \emph{contains} $l$ if $i<l<j$.
\end{defin}

\begin{lemma}
\label{thickness_varies_by_1}
\item (i) The thickness of consecutive layers varies at most by 1.
\item(ii) If $(i,j)$ is a thick interval (for $(\sigma_k),(\tau_k)$), then $\sigma_i, \tau_i$ are disjoint.
\end{lemma}
\proof
Both parts follow immediately from Lemma \ref{difference1}.
\qed

\begin{defin}
\label{characteristic disc and surface}
Let $(i,j)$ be a thick interval (for $(\sigma_k),(\tau_k)$).
Let vertices
$s_k\in \sigma_k,\ t_k\in \tau_k$ be such that for each $i\leq k\leq j$ the
distance $|s_kt_k|$ is maximal (i.e. $s_k,t_k$ realize the thickness of the layer $k$).
By Lemma \ref{thickness_varies_by_1}(ii) the sequence $s_i,s_{i+1},\ldots, s_j, t_j,t_{j-1},\ldots t_i,s_i$ is
an embedded loop, thus we can consider a minimal surface $S\colon \Delta\rightarrow X$
spanned on this loop (c.f. Definition \ref{def_of_surface}). We say that $S$ is a
\emph{characteristic
surface} (for the thick interval $(i,j)$) and $\Delta$ is a \emph{characteristic disc}.
\end{defin}

\begin{lemma}
\label{realizing_thickness_in_pairs_implies_realizing_as_a_pair}
For $s_k,s'_k\in \sigma_k, \ t_k,t'_k\in \tau_k$, if distances $|s_kt'_k|,|s'_kt_k|$ equal the thickness of the layer $k$ then also $|s_kt_k|$ equals the thickness of the layer $k$, i.e. if vertices $s_k\in\sigma_k, t_k\in\tau_k$ realize the thickness in some pairs, then they also realize the thickness as a pair.
\end{lemma}
\proof
Immediate from definition of thickness and Corollary \ref{corollary_2_with_leq_to_no_trapezoid}.
\qed
\medskip

The lemma below summarizes the geometry of characteristic discs, which we need to introduce the concept of a Euclidean geodesic. The special features of characteristic discs, in the case that $(\sigma_k),(\tau_k)$ are directed geodesics, will be given in Lemma \ref{properties_of_special_characteristic discs} at the end of this section.

Let $S\colon \Delta\rightarrow X$ be a characteristic surface. Denote by $v_k,w_k\in\Delta$ the preimages of $s_k,t_k$ in $X$, respectively. This notation will be fixed for the entire article. Let us point out that we use numbers $i,\ldots,j$ to number the layers in $\Delta$ (c.f. Definition \ref{layer}) between $v_iw_i$ and $v_jw_j$, instead of $0,\ldots j-i$, for the sake of clarity.

\begin{lemma}
\label{properties_of_characteristic discs}
\item (i) $\Delta$ (and thus the characteristic surface $S$) is wide and flat,
\item (ii) if we embed $\Delta\subset \mathbb{E}_{\Delta}^2$, then the edges $v_iw_i$ and $v_jw_j$ are parallel and consecutive layers between them are contained in consecutive straight lines (treated as subcomplexes
    of $\mathbb{E}_{\Delta}^2$) parallel to the lines containing $v_iw_i$ and $v_jw_j$.
\end{lemma}
\proof
(i) To prove wideness it is enough to show that any nonconsecutive vertices of the boundary loop are at distance
$\geq 2$. Since the layers $k$, where $i< k< j$, are thick (for $(\sigma_k),(\tau_k)$), the only possibility for this to fail is that (w.l.o.g.) $|s_kt_{k+1}|=1$ for some $i< k< j$. If this happens, then
both $s_k$ and $t_k$ lie in the projection of $t_{k+1}$ onto the layer $k$ between $\sigma$ and $\tau$ (the projection is defined by Lemma \ref{layers close}(ii)), hence
they are neighbors (Lemma \ref{projection lemma}), which contradicts $|s_kt_k|\geq 2$. Thus a characteristic disc is wide.

Before proving flatness, we need the following general observation. If $\Gamma$ is a 1--skeleton geodesic, which is in the boundary of a triangulation of a disc, then the sum of the defects at the vertices in the interior of $\Gamma$ is $\leq 1$. Moreover, all the defects at these vertices are $\leq 1$ and each two vertices of positive defect are separated by a vertex of negative defect.

To prove flatness, compute possible defects at the boundary vertices of
$\Delta$. By wideness, they are $\leq 1$ at $v_i,v_j,w_i,w_j$. Moreover, their sum over the interior vertices
of each of the 1--skeleton geodesics $(v_k)_{k=i}^j,(w_k)_{k=i}^j$ is $\leq 1$ (they are 1--skeleton geodesics, since their images are). Thus Gauss--Bonnet Lemma \ref{Gauss-Bonnet}
implies that the defects of the interior vertices are equal to zero, the sums of the defects over the vertices
$(v_k)_{k=i+1}^{j-1},(w_k)_{k=i+1}^{j-1}$ equal 1 each
and the defects at $v_i,v_j,w_i,w_j$ are equal to $1$.

We now want to say more about the defects at $(v_k)_{k=i+1}^{j-1}$. Up to now we know that their sum is 1, they equal $1,0,-1$ or $-2$ and each two vertices of positive defect are separated by a vertex of negative defect (since $(v_k)_{k=i}^{j}$ is a 1--skeleton geodesic). This implies
that the defects equal alternatingly $1,-1,1-1,\ldots, 1$ with possible $0$'s between them. The same holds for
the defects at $(w_k)_{k=i+1}^{j-1}$. Thus, by Lemma \ref{characterization_of_flat} (characterization of flatness), the characteristic disc $\Delta$ is flat, i.e. we have an embedding $\Delta\subset \mathbb{E}_{\Delta}^2$ isometric on the 1--skeleton.

(ii) By the computation of defects in the proof of (i) we get that the edges $v_iw_i$ and $v_jw_j$ are parallel in $\mathbb{E}_{\Delta}^2$. We also get that $v_k,w_k$, for $i\leq k\leq j$ are at combinatorial distances $k-i,j-k$ from the lines containing the edges $v_iw_i,v_jw_j$. Hence $v_k,w_k$ lie on the appropriate line parallel to $v_iw_i$ and the vertices of $\Delta$ split into families lying on geodesics $v_kw_k$. By convexity of layers, Remark \ref{layers convex}, (or by direct observation) these geodesics are equal to the layers.
\qed
\medskip

When speaking about the layers in $\Delta$ between $v_iw_i$ and $v_jw_j$, we will often skip "between $v_iw_i$ and $v_jw_j$".

\begin{rem}
\label{layers preserved 0}
Denote the layer $k$ in $\Delta$ (between $v_iw_i$ and $v_jw_j$) by $L_k$. Then $S(L_k)$ is contained in the layer $k$ in $X$ between $\sigma$ and $\tau$. This follows from
\begin{align*}
S(L_k)&\subset S(B_{k-i}(v_iw_i))\cap S(B_{j-k}(v_jw_j))\subset
\\
&\subset
 B_{k-i}(S(v_iw_i))\cap B_{j-k}(S(v_jw_j))\subset
 \\
&\subset B_{k-i}(\sigma_i\tau_i)\cap B_{j-k}(\sigma_j\tau_j)\subset
\\
&\subset B_{k-i}(B_i(\sigma))\cap B_{j-k}(B_{n-j}(\tau))
= B_k(\sigma)\cap B_{n-k}(\tau).
\end{align*}
\end{rem}
\medskip

The next lemma summarizes some uniqueness properties of characteristic surfaces for a fixed thick interval $(i,j)$.

\begin{lemma}
\label{properties of characteristic surfaces}
\item(i) A characteristic surface is almost geodesic. In particular, it is an isometric embedding on the 1--skeleton
of a subcomplex spanned by any pair of consecutive layers between $v_iw_i$ and $v_jw_j$ in $\Delta$.
\item (ii) A characteristic disc $\Delta\subset \mathbb{E}_{\Delta}^2$ does not depend (up to isometry) on the choice
of $s_k,t_k$ and the choice of a characteristic surface.
\item If we have two characteristic surfaces $S_1\colon\Delta_1\rightarrow X, S_2\colon\Delta_2\rightarrow X$, then after identifying the characteristic discs $\Delta_1=\Delta_2$ (which is possible by (ii)) we have that
\item (iii) for any vertices $x,y\in \Delta_1=\Delta_2$ at distance $1$, $S_1(x)$ and $S_2(y)$ are also at distance $1$, i.e. for any two characteristic surfaces $S_1,S_2$ we can substitute an image of a vertex of the first surface with the corresponding image in the second and get another characteristic surface,
    \item (iv) for any vertex $x\in \Delta_1=\Delta_2$, $S_1(x)$ and $S_2(x)$ are at distance $\leq 1$.
\end{lemma}
\proof
(i) This follows from Elsner's Theorem \ref{main_Elsner_theorem}, since, by Lemma \ref{properties_of_characteristic discs}(i), a characteristic disc is
flat and wide. The second part follows from the fact that any two vertices in a same or consecutive layers in $\Delta\subset \mathbb{E}_{\Delta}^2$
can be connected by a neat geodesic, which can be verified by direct observation.
\medskip

(ii) Observe that, by Lemma \ref{properties_of_characteristic discs}(ii), the isometry class of $\Delta$ is determined by the distances $|v_kw_k|$,
$|v_kw_{k+1}|$, for $i\leq k\leq j-1$, which are equal, by (i), to $|s_kt_k|, |s_kt_{k+1}|$, respectively. The value $|s_kt_k|$ equals the thickness of the layer $k$, so it does not depend on the choices. To prove the same for $|s_kt_{k+1}|$,
consider two characteristic surfaces constructed for choices  $s_l,s'_l\in \sigma_l,\ t_l,t'_l\in \tau_l$, where $l=k,k+1$.
We will prove that $|s_kt_{k+1}|=|s_kt_{k+1}'|=|s_k't_{k+1}'|$. We restrict ourselves to proving the first equality (the second is proved analogically).
By Lemma \ref{realizing_thickness_in_pairs_implies_realizing_as_a_pair}
we have that $|s_{k+1}t_{k+1}'|$ is the thickness of the layer $k+1$.
Thus there is a characteristic surface spanned on a loop passing through $s_k, t_k, s_{k+1}, t_{k+1}'$. Hence, by (i), the distance $|s_kt_{k+1}'|$ is determined by $|s_kt_k|$ and $|s_{k+1}t_k|$, thus it is the same as
$|s_kt_{k+1}|$, as desired.
\medskip

(iii) If $x$ and $y$ are both boundary vertices, then this is obvious. Otherwise, w.l.o.g.
assume that $x$ is an interior vertex of $\Delta$. Suppose $x$ lies in the layer $k$ (we denote it by $L_k$) in $\Delta$ between $v_iw_i$ and $v_jw_j$. Denote the thickness of the layer $k$ for $(\sigma_t),(\tau_t)$ by $d$.

First consider the case that $y\in L_k$. By Remark \ref{layers preserved 0}
we have that $S_1(L_k)$ and $S_2(L_k)$ lie in the layer $k$ in $X$ between $\sigma$ and $\tau$. By Lemma \ref{realizing_thickness_in_pairs_implies_realizing_as_a_pair} we
have that $|S_2(v_k)S_1(w_k)|=|S_1(v_k)S_2(w_k)|=d$.
Hence Corollary \ref{corollary_to_no_trapezoid} applied to $S_1(L_k)$ and $S_2(L_k)$ gives $|S_1(x)S_2(y)|=1$, as desired.

Now, w.l.o.g., consider the remaining case that $y$ is in the layer $k-1$ (denoted by $L_{k-1}$) in $\Delta$ between $v_iw_i$ and $v_jw_j$. Denote by $y',x''$ the common neighbors of $x,y$ in $L_{k-1},L_k$, respectively, and by $x'$ the neighbor of $x$ in $L_k$ different from $x''$.  Then, from the previous case, we have that $S_1(x)S_2(x')S_2(y')S_2(y)S_2(x'')S_1(x)$ is a loop of length 5, hence it is triangulable. By (i), all $|S_2(x')S_2(x'')|,\ |S_2(x')S_2(y)|,\ |S_2(x'')S_2(y')|$ equal 2, hence we obtain $|S_1(x)S_2(y)|=1$, as desired.

Observe that this proof actually implies Proposition \ref{stable_geodesics_in
characteristic image} in the case that $\gamma\subset v_kw_k$ for some $k$.
\medskip

(iv) For boundary vertices this is obvious. For an interior vertex $x$, let $x',x''$ be its neighbors in a common layer in $\Delta$ between $v_iw_i, v_jw_j$. Then, by (iii), we have that $S_1(x)S_2(x')S_2(y)S_2(x'')S_1(x)$ is a loop of length 4. Moreover, by (i), we have that $|S_2(x')S_2(x'')|=2$. Thus $|S_1(x)S_2(y)|\leq 1$, as desired.
\qed
\medskip

As a corollary, the following definition is allowed.

\begin{defin}
\label{characteristic image} Let $\rho$ be a simplex of the
characteristic disc $\Delta$ for some thick interval $(i,j)$ (for
$(\sigma_k),(\tau_k)$). Its \emph{characteristic image} is a simplex
in $X$, denoted by $\mathcal{S}(\rho)$, which is the span of the
images of $\rho$ under all possible characteristic surfaces. Note
that $\mathcal{S}(\rho)$ is a simplex by Lemma \ref{properties of
characteristic surfaces}(iii,iv), and if $\rho\subset \rho'$, then $\mathcal{S}(\rho)\subset \mathcal{S}(\rho')$, i.e. $\mathcal{S}$ respects inlucions. The \emph{characteristic image} of a subcomplex of $\Delta$
is the union of the characteristic images of all its simplices. We call this assignment the
\emph{characteristic mapping}.

If $\overline{v}$ is a vertex in $\mathcal{S}(\Delta)$, we denote by
$\mathcal{S}^{-1}(\overline{v})$ the vertex $v\in \Delta$
such that $\mathcal{S}(v)\ni\overline{v}$. We claim that this vertex is unique. Indeed, characteristic images of different layers in $\Delta$ between $v_iw_i, v_jw_j$ are disjoint since, by Remark \ref{layers preserved 0}, they lie in different layers in $X$ between $\sigma,\tau$, disjoint by Lemma \ref{layers close}. Moreover, by Lemma \ref{properties of characteristic surfaces}(i,iii), we have that $S_1(v)\neq S_2(v')$ for any characteristic surfaces $S_1,S_2$ and any vertices $v\neq v'$ in a common layer in $\Delta$. This justifies the claim. If $\overline{\rho}$ is a simplex
in $\mathcal{S}(\Delta)$, we denote by $\mathcal{S}^{-1}(\overline{\rho})$
the span of the union of $\mathcal{S}^{-1}(\overline{v})$ over all $\overline{v}\in\overline{\rho}$. We have that $\mathcal{S}^{-1}(\overline{\rho})$ is a simplex, by Remark \ref{layers preserved 0}, Lemma \ref{layers close}, and Lemma \ref{properties of
characteristic surfaces}(i,iii).
If $Y$ is a subcomplex
of $\mathcal{S}(\Delta)$, we denote by $\mathcal{S}^{-1}(Y)$
the union of $\mathcal{S}^{-1}(\overline{\rho})$ over all $\overline{\rho}\subset Y$.
\end{defin}

Having established the uniqueness properties of characteristic
surfaces, we start to exploit the $CAT(0)$ structure of the
corresponding characteristic discs. From now on, up to the end of
Section \ref{Contracting property}, unless stated otherwise, assume
that $(\sigma_k), (\tau_k)$ are the directed geodesics between
$\sigma,\tau$.

\begin{defin}
\label{euslidean diagonal} Let $(i,j)$ be a thick interval and
let $\Delta\subset \mathbb{E}_{\Delta}^2$ be its characteristic disc. We
will define a sequence of simplices $\rho_k\in \Delta$, where
$i<k<j$, which will be called the \emph{Euclidean diagonal} of the
characteristic disc $\Delta$.

Let $v_k', w_k'$ be points (barycenters of edges) on the straight
line segments $v_kw_k$ at distance $\frac{1}{2}$ from $v_k,w_k$,
respectively. In particular $v'_i=w'_i, \ v'_j=w'_j$. Consider the closed
polygonal domain $\Delta' \subset \Delta$ enclosed by the piecewise linear loop with consecutive vertices
$v'_i,v'_{i+1},\ldots ,v'_j=w'_j,w'_{j-1},\ldots, w'_i=v'_i$.
Note that, since $\Delta'$ is simply--connected, it is $CAT(0)$ with the Euclidean path metric induced from $\mathbb{E}_{\Delta}^2$
identified with $\mathbb{E}^2$. We
call $\Delta'$ a \emph{modified characteristic disc}. Let
$\gamma'$ be the $CAT(0)$ geodesic joining $v'_i=w'_i$ to $v'_j=w'_j$
in $\Delta'$. We call $\gamma'$ a \emph{CAT(0) diagonal} of
$\Delta$. For each $i<k<j$, among the vertices of $\Delta$ lying in
the interior of the 1--skeleton geodesic $v_kw_k$ find the ones nearest to
$\gamma'\cap v_kw_k$. For each $k$ this is either a single vertex or
two vertices spanning an edge (if $\gamma'$ goes through its
barycenter and $v_k,w_k$ are not some of its vertices). We put
$\rho_k$ equal to this vertex or this edge, accordingly.
\end{defin}

At first sight it might seem strange that in the above definition we pass
to $\Delta'$ and take the geodesic $\gamma'$ there instead of doing
it in $\Delta$ itself. However, this construction allows us to
exclude $v_k,w_k$ from being in $\rho_k$, which a careful reader
will find to be a necessary condition for the arguments of the
combinatorial Proposition
\ref{main_proposition_to_euclidean_geodesics} to be valid.

Here are some basic properties of the Euclidean diagonals.

\begin{lemma}
\label{properties_of_euclidean_diagonal}
\item (i) Each pair of consecutive $\rho_k, \rho_{k+1}$, for $i< k<j-1$, spans a simplex.
\item (ii) $\rho_{i+1}, v_i, w_i$ span a simplex and $\rho_{j-1}, v_j, w_j$ span a simplex.
\end{lemma}
\proof
Part (ii) is obvious, since we excluded $v_k,w_k$ from being in $\rho_k$.
To prove (i), consider $\Delta'\subset \Delta \subset \mathbb{E}_{\Delta}^2$ oriented in such a way that
$v_kw_k$ are horizontal, this is possible by Lemma \ref{properties_of_characteristic discs}(ii). Moreover, Lemma  \ref{properties_of_characteristic discs}(ii) yields that the boundary of $\Delta'$ consists of line segments at angle $30^\circ$ from the vertical direction. Let $\gamma'$ be as in Definition \ref{euslidean diagonal}. It is a broken line with vertices at the boundary of $\Delta'$.

We claim that any line segment of $\gamma'$
is at angle $< 30^\circ$ from the vertical direction. First we prove that this angle is $\leq 30^\circ$. Otherwise, let $p$ be an endpoint of such a line segment. Obviously $p$ is different from the endpoints of $\gamma'$.
The interior angle at $p$ between the segment of $\gamma'$ and any of the boundary line segments of $\Delta'$ is $< 180^\circ$, which contradicts the fact that $p$ is an interior vertex of a geodesic $\gamma'$.
Thus we proved that  any line segment of $\gamma'$
is at angle $\leq 30^\circ$ from the vertical direction.

If for some line segment of $\gamma'$ this angle equals $30^\circ$, then by the previous considerations the whole $\gamma'$ is in fact a straight line at angle $30^\circ$ from
the vertical. This implies that the defects at all vertices in $(v_k)_{k=i+1}^{j-1}$ or
all vertices in $(w_k)_{k=i+1}^{j-1}$ are zero. Contradiction.
\medskip

Now part (i) follows from the following observation, whose proof is easy and is left for the reader. Consider two consecutive horizontal lines $\alpha_1,\alpha_2$ in $\mathbb{E}_{\Delta}^2$. Let $\beta$ be some straight line segment joining points $p\in \alpha_1, r\in \alpha_2$ at
angle $< 30^\circ$ from the vertical direction. Then there exist two 2--simplices $abc, bcd$ in $\mathbb{E}_{\Delta}^2$ such that $ab\subset \alpha_1, cd\subset \alpha_2$ and  $p\in ab, r\in cd$. Moreover, it cannot happen simultaneously that $|pa|\leq |pb|$ and $|rd|\leq |rc|$.
\qed

\medskip
Thus we can finally introduce the main definition of this section.

\begin{defin}
\label{euclidean geodesic}
We define a sequence of simplices $\delta_k$, where $0\leq k\leq n$, which is called the \emph{Euclidean geodesic} between $\sigma, \tau$, as follows. For each $k$, if the layer $k$ is thin, then we take
$\delta_k$ as the span of $\sigma_k$ and $\tau_k$.

If the layer $k$ is thick, consider the
thick interval $(i,j)$ which contains $k$. Let $\rho_k$ be an appropriate simplex of the Euclidean diagonal of the
characteristic disc $\Delta$ for $(i,j)$ (c.f. Definition \ref{euslidean diagonal}). We take $\delta_k = \mathcal{S}(\rho_k)$ (c.f. Definition \ref{characteristic image}).
\end{defin}

\begin{rem}
\label{delta v_i}
In the above setting, we have $\sigma_i=\mathcal{S}(v_i), \ \tau_i=\mathcal{S}(w_i)$, by Lemma \ref{thickness_varies_by_1}(ii).
Hence $\delta_i=\Span\{\sigma_i,\tau_i\}=\mathcal{S}(v_iw_i)$.
\end{rem}

\begin{rem}
\label{symmetry of euclidean geodesic}
By the symmetry of the construction, the Euclidean geodesic between
$\sigma$ and $\tau$ becomes the Euclidean geodesic between
$\tau$ and $\sigma$ if we take the simplices of the sequence in the opposite order.
\end{rem}

Here is the justification for using the name "geodesic" in Definition \ref{euclidean geodesic}.

\begin{lemma}
\label{properties of euclidean geodesics}
\item (i) For any $0\leq k< l \leq n$ we have that $\delta_k\subset S_{l-k}(\delta_{l}),\delta_{l}\subset S_{l-k}(\delta_k)$.
\item (ii) For any $0\leq k\leq n-1$ if the layer $k$ or the layer $k+1$ is thick, then $\delta_k$ and  $\delta_{k+1}$ span a simplex.
\item (iii) For any $0\leq l<m\leq n$ such that there exists $l\leq k\leq m $ such that the layer $k$ is thick, and for any vertices $x\in\delta_m, y\in \delta_l$, we have
$|xy|=m-l$.
\end{lemma}
\proof
Assertion (ii) follows from Lemma \ref{properties_of_euclidean_diagonal}(i,ii), Remark \ref{delta v_i} and Lemma \ref{properties of characteristic surfaces}(iii,iv).

To prove assertion (i), say the first inclusion, observe that for any $0\leq k<n$ we have $\Span (\sigma_k\cup \tau_k)\subset B_1(\Span (\sigma_{k+1}\cup \tau_{k+1}))$. Hence, assertion (ii) gives already, for any $0\leq k< l \leq n$, that $\delta_k\subset B_{l-k}(\delta_{l})$. Then $\delta_k\subset S_{l-k}(\delta_{l})$ follows from Remark \ref{layers preserved 0} and Lemma \ref{layers close}(ii).

To prove part (iii),  assume that $l<k<m$ (other cases are easier). Take any vertex $z\in \delta_k$. Then, by (i), there are vertices $x'\in \delta_{k-1},\ y'\in \delta_{k+1}$ such that $|xx'|=(k-1)-l,\ |yy'|=m-(k+1)$. By (ii) (and (i)), we have $|zx'|=|zy'|=1$. Hence $|xy|\leq m-l$ and by (i) we have $|xy|= m-l$, as desired.
\qed
\medskip

Now we state an extra property of characteristic discs in the case that $(\sigma_k)$ (but $(\tau_k)$ not necessarily) is the directed geodesic. This property was not necessary for Definition \ref{euclidean geodesic}, but will become indispensable in the next section.

\begin{lemma}
\label{properties_of_special_characteristic discs}
\item (i) If the defect at some $v_k$, where $i+1<k<j-1$, equals $-1$, then the defect at $v_{k+1}$ equals 1.
\item (ii) The defect at $v_{i+1}$
equals 1.
\end{lemma}
\proof
(i) Proof by contradiction. Suppose the defect at some $v_k$, where $i+1<k<j-1$, equals $-1$, and the defect at $v_{k+1}$ equals 0. Denote by $x$ the vertex next to $v_{k+1}$ on the 1--skeleton geodesic $v_{k+1}w_{k+1}$ and by $y$ the vertex next to $v_{k}$ on the 1--skeleton geodesic $v_{k}w_{k}$. We aim to prove that, for any characteristic surface $S$,
$S(x)$ belongs to $\sigma_{k+1}$. Suppose for a moment we have already proved this. Then, since by Lemma \ref{properties of characteristic surfaces}(i) we have $|S(x)S(v_{k+2})|=2$ and at the same time $S(v_{k+2})\in \sigma_{k+2}$,  we get a contradiction.

Now we prove that $S(x)\in \sigma_{k+1}$. By Remark \ref{layers preserved 0}, $S(x)$ lies in the layer $k+1$ between $\sigma$ and $\tau$. Now by definition of projection (c.f. Definition \ref{projection}) we need to prove that $S(x)$ is a neighbor of each $\bar{z}\in \sigma_k$. Case $\bar{z}=S(y)$ is obvious, so suppose $\bar{z}\neq S(y)$. Since, by definition of thickness, $|\bar{z}S(w_k)|\leq |S(v_k)S(w_k)|$, we have by Lemma \ref{corollary_to_no_trapezoid} (applied to $r_0=S(v_k),r_1=S(y),r_d=p_d=S(w_k)$ and to $p_0=\bar{z}$ in case of $|\bar{z}S(w_k)|= |S(v_k)S(w_k)|$ or to $p_0=S(v_k),p_1=\bar{z}$ in case of $|\bar{z}S(w_k)|< |S(v_k)S(w_k)|$)
that $|\bar{z}S(y)|=1$. Considering the loop $\bar{z}S(y)S(x)S(v_{k+1})\bar{z}$, since
$|S(y)S(v_{k+1})|=|yv_{k+1}|=2$ (Lemma \ref{properties of characteristic surfaces}(i)), we get $|\bar{z}S(x)|=1$, as desired.
\medskip

(ii) By contradiction. Denote by $x$ the vertex between $v_{i+1}$ and $w_{i+1}$ on the 1--skeleton geodesic $v_{i+1}w_{i+1}$.
Since $\sigma_i=\mathcal{S}(v_i)$ (see Remark \ref{delta v_i}), we have by Remark \ref{layers preserved 0} and Lemma \ref{properties of characteristic surfaces}(iii) that $S(x)$ belongs to $\sigma_{i+1}$. By Lemma \ref{properties of characteristic surfaces}(i) we have $|S(x)S(v_{i+2})|=2$. At the same time $S(v_{i+2})\in \sigma_{i+2}$, contradiction.
\qed
\medskip

We will repeat some steps of this proof later on in the proof of Lemma \ref{position of
directed_geodesic_simplices_in_characteristic_image}. We decided, for clarity, not to interwind these two proofs.

As a consequence of Lemma \ref{properties_of_special_characteristic discs}, we get the following lemma, whose proof, similar to the proof of Lemma \ref{properties_of_euclidean_diagonal}, we omit. Here we assume that both $(\sigma_k),(\tau_k)$ are directed geodesics.

\begin{lemma}
\label{transversality}
If $j-i>2$ then the $CAT(0)$ diagonal $\gamma'$ in $\Delta$ crosses each line orthogonal to the layers transversally.
\end{lemma}

\section{Directed geodesics between simplices of Euclidean geodesics}
\label{Directed geodesics between simplices of Euclidean geodesics}

In this section we start to prove a weak version of Theorem \ref{2}, which concerns one of the main properties of Euclidean geodesics. Roughly speaking, the theorem says that pieces
of Euclidean geodesics are coarsely also Euclidean geodesics.

We keep the notation from the previous section. The simplices
$(\sigma_k),(\tau_k)$ are in this section the directed geodesics
between $\sigma,\tau$.

\begin{theorem}[weak version of Theorem \ref{2}]
\label{main_theorem_on_euclidean_geodesics}
Let $\sigma,\tau$ be simplices of a systolic complex $X$, such that for some natural $n$ we have
$\sigma \subset S_n(\tau), \tau \subset S_n(\sigma)$ (as required in the definition of the Euclidean geodesic).
Let $(\delta_k)_{k=0}^n$ be the Euclidean geodesic between $\sigma$ and $\tau$. Take some
$0\leq l<m\leq n$ and consider the simplices $\tilde{\delta}_l=\delta_l,\tilde{\delta}_{l+1},\ldots, \tilde{\delta}_m=\delta_m$ of the Euclidean geodesic between
$\delta_l$ and $\delta_m$ (we can define it by Lemma \ref{properties of euclidean geodesics}(i)). Then for each $l\leq k\leq m$ we have $|\delta_k,\tilde{\delta}_k|\leq 3$.
\end{theorem}

The proof of Theorem \ref{main_theorem_on_euclidean_geodesics}
splits into two steps. The first step is to prove that directed
geodesics between $\delta_l$ and $\delta_m$ stay close to the
union of characteristic images of all characteristic discs (for
$(\sigma_k),(\tau_k)$). This is the content of Proposition
\ref{main_proposition_to_euclidean_geodesics}, whose proof occupies
the rest of this section.

The second step is to check that characteristic images for the directed geodesics between $\delta_l$ and $\delta_m$ also stay close to the
union of characteristic images for $(\sigma_k),(\tau_k)$. Properties of layers actually imply that characteristic discs of the former are embedded into characteristic discs of the latter, modulo small neighborhood of the boundary. So everything boils down to the fact that Theorem \ref{main_theorem_on_euclidean_geodesics} is valid for $CAT(0)$ subspaces of the Euclidean plane.  We carry out this program in the next
section. We also indicate there an argument, how to promote Theorem \ref{main_theorem_on_euclidean_geodesics} to Theorem \ref{2}, with a reasonable constant $C$.

A complete alternative proof of Theorem \ref{2}, with a worse constant $C$, is obtained as a consequence of Proposition \ref{Euclidean_near_CAT(0)}. We present it at the end of Section \ref{Characteristic discs for Euclidean geodesics}.
We advise the reader to have a look at the proof of Theorem \ref{main_theorem_on_euclidean_geodesics} via Proposition \ref{main_proposition_to_euclidean_geodesics}. This proof is straightforward
and allows us to introduce gradually some concepts needed later. However, to save time, one can skip the remaining part of Section \ref{Directed geodesics between simplices of Euclidean geodesics}, go over the definitions in Section \ref{CAT(0) geometry of characteristic discs} and then go directly to Section \ref{Characteristic discs for Euclidean geodesics}.
\medskip

For each thick layer $l\leq k\leq m$ contained in a thick interval
$(i,j)$ (for $(\sigma_t),(\tau_t)$; from now on we often skip "for $(\sigma_t),(\tau_t)$"), denote by $\alpha_k$ the appropriate simplex (in the
corresponding characteristic disc $\Delta$) of the directed geodesic
from $\rho_l$, if $i<l$, or $v_i$ otherwise, to $\rho_m$, if $m<j$,
or $v_j$ otherwise. The simplices $(\tilde{\sigma}_k)_{k=l}^{m}$
of the directed geodesic from $\delta_l$ to $\delta_m$ satisfy the
following.

\begin{proposition}
\label{main_proposition_to_euclidean_geodesics} Let $l\leq k\leq m$.
\item (i) If the layer $k$ is thin, then $\sigma_k$ contains or is
contained in $\tilde{\sigma}_k$,
\item (ii) if the layer $k$ is thick, then $\mathcal{S}(\alpha_k)$ contains or is contained in
$\tilde{\sigma}_k$.
\end{proposition}

Before we give the proof of Proposition \ref{main_proposition_to_euclidean_geodesics}, we need to establish some necessary lemmas. The first one describes
the position of $\sigma_k$ with respect to the characteristic image.
Like in Lemma \ref{properties_of_special_characteristic discs}, here $(\tau_k)$ does not need to be the directed geodesic.

\begin{lemma}
\label{position of
directed_geodesic_simplices_in_characteristic_image} For a thick
layer $k$ let $x_k$ be the vertex, which is a neighbor of $v_k$ on
the 1--skeleton geodesic $v_kw_k$ in the characteristic disc for the
thick interval containing $k$. If the defect at $v_k$ equals 1, then
$\sigma_k=\mathcal{S}(v_kx_k)$. Otherwise $\sigma_k=
\mathcal{S}(v_k)$.
\end{lemma}
\proof First of all $\sigma_k\subset\mathcal{S}(v_kx_k)$ follows
from the definition of thickness and Proposition
\ref{stable_geodesics_in characteristic image} (one could also verify this by hand, similarly like in the proofs of Lemma \ref{properties_of_characteristic discs}(iii) and Lemma \ref{properties_of_special_characteristic discs}(i)). Suppose the defect
at $v_k$ is $\neq 1$. Hence $|v_{k-1}x_k|=2$, by Lemma \ref{properties_of_special_characteristic discs}(i,ii). The inclusion
$\mathcal{S}(v_k)\subset\sigma_k$ is obvious and the converse inclusion
follows from $\sigma_k\subset\mathcal{S}(v_kx_k)$ and from Lemma
\ref{properties of characteristic surfaces}(i).

Now suppose the defect at $v_k$ equals $1$. If the layer $k-1$ is thick, then the defect at $v_{k-1}$ is $\neq 1$
and we apply what we have just proved to get
$\mathcal{S}(v_{k-1})=\sigma_{k-1}$. If the layer $k-1$ is thin we get
immediately that $\mathcal{S}(v_{k-1})=\sigma_{k-1}$ (Remark \ref{delta v_i}). In both cases
using Remark \ref{layers preserved 0}, Lemma \ref{properties of characteristic surfaces}(iii), and the definition of projection
we get $\mathcal{S}(v_kx_k)\subset \sigma_k$, as desired. \qed
\medskip

As a corollary we get the following technical lemma.

\begin{lemma} \label{agreeing_of_projections} Suppose $k<m$ do not
satisfy $i\leq k<m<j$ for any thick interval $(i,j)$ or if they
violate this then  $|v_{k+1},\rho_m|=m-(k+1)$. Then the projection
of $\sigma_k$ onto $B_{m-(k+1)}(\delta_m)$ equals $\sigma_{k+1}$.
\end{lemma}
\proof
To justify speaking about the projection
of $\sigma_k$ onto $B_{m-(k+1)}(\delta_m)$ we must show that $\sigma_k\subset
S_{m-k}(\delta_m)$. The simplex $\sigma_k$ is outside
$B_{m-k-1}(\delta_m)$ by Remark \ref{layers preserved 0} and Lemma \ref{layers close}. Thus we only need to check that $\sigma_k\subset
B_{m-k}(\delta_m)$.

To verify this, we prove that $\sigma_{k+1}\subset
B_{m-(k+1)}(\delta_m)$. If the layer $k+1$ is thin then this follows
from Lemma \ref{properties of euclidean geodesics}(i).
If the layer $k+1$ is thick, then denote by $(i,j)$ the thick interval containing
$k+1$. By Lemma \ref{position of
directed_geodesic_simplices_in_characteristic_image} we have
$\sigma_{k+1}\subset \mathcal{S}(v_{k+1}x_{k+1})$ ($x_{k+1}$ as in
Lemma \ref{position of
directed_geodesic_simplices_in_characteristic_image}). Thus it is
enough to establish the inclusion
$\mathcal{S}(v_{k+1}x_{k+1})\subset B_{m-(k+1)}(\delta)$. If $m<j$,
then this follows from our assumptions. If $j\leq m$, then from Remark \ref{delta v_i} and Lemma \ref{properties
of euclidean geodesics}(i) we have
\begin{align*}
\mathcal{S}(v_{k+1}x_{k+1})\subset \mathcal{S}(B_{j-(k+1)}(v_j))\subset B_{j-(k+1)}(\mathcal{S}(v_j))&\subset
\\
\subset B_{j-(k+1)}(\delta_j)&\subset B_{m-(k+1)}(\delta_m),
\end{align*}
as desired.

Hence the projection of $\sigma_k$ onto $B_{m-(k+1)}(\delta_m)$ is defined. Denote it by $\pi$. Since $B_{m-(k+1)}(\delta_m)\subset
B_{n-(k+1)}(\tau)$, we have $\pi\subset \sigma_{k+1}$. For the converse
inclusion we need $\sigma_{k+1}\subset
B_{m-(k+1)}(\delta_m)$, which we have just proved.
\qed
\medskip

The next lemma is valid for any $(\sigma_k),(\tau_k)$, not necessarily directed geodesics.

\begin{lemma}
\label{projection_in_geodesic} Let $e$ be an edge in the layer $k$ of
$\Delta$ (between $v_iw_i,v_jw_j$), such that $e$ has three neighboring vertices in the layer
$k+1$. Let $\overline{x}$ be a~vertex in the residue of $S(e)$ (for
some characteristic surface $S$) in the layer $k+1$ between
$\sigma,\tau$ in $X$. Then $\overline{x}\in \mathcal{S}(x)$, where
$x$ is the vertex in the layer $k+1$ of $\Delta$ in the residue of $e$.
\end{lemma}
\proof Denote by $y_1,y_2$ the neighbors of $e$ in the layer $k+1$ of
$\Delta$ different from $x$, and let
$\overline{y}_1=S(y_1),\overline{y}_2=S(y_2)$. We claim that
$\overline{y}_1,\overline{y}_2$ are neighbors of $\overline{x}$. Indeed, let
$z_1$ be the vertex in $e$, which is a neighbor of $y_1$. Let
$\overline{z}_1=S(z_1)\subset S(e)$. Observe that both
$\overline{y}_1,\overline{x}$ lie in the projection of
$\overline{z}_1$ onto $B_{n-(k+1)}(\tau)$ (by Remark \ref{layers
preserved 0}), hence, by Lemma \ref{projection lemma}, they are
neighbors, as desired. Analogically, $\overline{y}_2,\overline{x}$
are neighbors. Thus, by the easy case of Proposition
\ref{stable_geodesics_in characteristic image}, $\overline{x}\in
\mathcal{S}(x)$, as required. \qed

\medskip
The following lemma describes the behavior of the simplices
$\alpha_k$ appearing in the statement of Proposition
\ref{main_proposition_to_euclidean_geodesics}. The proof of Lemma
\ref{directed geodesics in characteristic discs} requires Lemma
\ref{properties_of_special_characteristic discs}(i,ii), apart from
this it is straightforward and we skip it. For the same reason we will usually not invoke it in the proof of Proposition \ref{main_proposition_to_euclidean_geodesics}.

\begin{lemma}
\label{directed geodesics in characteristic discs}
Let $\Delta$ be a characteristic disc for some thick interval $(i,j)$.
Suppose for some $i\leq l<m\leq j$ we have simplices  $\alpha, \alpha'$ in the layers $l,m$ respectively between $v_iw_i,v_jw_j$ in $\Delta$. Suppose that $\alpha\subset S_{m-l}(\alpha')$ and $\alpha'\subset S_{m-l}(\alpha)$. Moreover, assume that $\alpha$ is an interior vertex of $\Delta$ or an edge disjoint with the boundary or $\alpha=v_i$. Assume that $\alpha'$ is an interior vertex or an edge disjoint with boundary or $\alpha'=v_j$. Let $(\alpha_k)_{k=l}^m$ be the directed geodesic in $\Delta$ joining $\alpha$ to $\alpha'$ (in particular $\alpha_l=\alpha, \alpha_m\subset \alpha'$). Then:

\item (i) If $\alpha_k$ is an edge, then $\alpha_{k+1}$ is
the unique vertex, which is in the residue of $\alpha_{k}$ in the layer
$k$.
\item (ii) If $\alpha_k=v_k$ and the defect at $v_k$ equals $0$, then $\alpha_{k+1}=v_{k+1}$.
\item (iii) If $\alpha_k$ is a vertex with two neighbors in the layer $k+1$, both at
distance $m-(k+1)$ from $\alpha'$, then
    $\alpha_{k+1}$ is an edge spanned by these two vertices.
\item (iv) If  $\alpha_k$ is a vertex with
two neighbors in the layer $k+1$, but only one of them at distance
$m-(k+1)$ from $\alpha'$, then $\alpha_{k+1}$ is this special
vertex.
\item Moreover, $\alpha_k$ never equals $w_k$.
If $\alpha_k$ is an edge containing $w_k$ then the defect at $w_k$
is $-1$.
If $\alpha_k=v_k$, then the defect at $v_k$ is
not equal to 1, except possibly for the cases $k=i,j$.
\end{lemma}

Now we are ready for the following.
\medskip\par\noindent\textbf{Proof of Proposition \ref{main_proposition_to_euclidean_geodesics}.}\ignorespaces
\ We will prove by induction on $k$, for $l\leq k\leq m$, the following
statement, which, by Lemma \ref{directed geodesics in characteristic
discs} and Lemma \ref{position of
directed_geodesic_simplices_in_characteristic_image}, implies the
proposition.
\medskip

\medskip\par\noindent\textbf{Induction hypothesis.}\ignorespaces
\ (1) If the layer $k$ is thick and $\alpha_k$ is an edge disjoint with
the boundary or meeting the boundary at a vertex of defect $\neq 1$,
then $\mathcal{S}(\alpha_k)$ is contained in $\tilde{\sigma}_k$,
\\ (2) if the layer $k$ is thick and $\alpha_k$ is a non--boundary vertex,
then $\mathcal{S}(\alpha_k)$ contains $\tilde{\sigma}_k$,
\\ (3) if the layer $k$ is thick and $\alpha_k$ is a boundary vertex or an edge intersecting the
boundary at a vertex of defect 1, or the layer $k$ is thin, then
$\sigma_k$ contains or is contained in $\tilde{\sigma}_k$.
\\

For $k=l$ the hypothesis is obvious. Suppose it is already proved
for some $l\leq k \leq m-1$. We would like to prove it for $k+1$. First
suppose that the layer $k$ is thick and $\alpha_k$ is an edge disjoint
with the boundary or meeting the boundary at a vertex of defect
$\neq 1$ (case (1)). Then $\alpha_{k+1}$ is a vertex. If it is a
boundary vertex, then $v_k\in
\alpha_k$. By the induction hypothesis, since the defect at $v_k$ is
not 1, $\mathcal{S}(\alpha_k)\subset \tilde{\sigma}_k$, moreover, by
Lemma \ref{position of
directed_geodesic_simplices_in_characteristic_image} we have
$\sigma_k\subset \mathcal{S}(\alpha_k)$, hence $\sigma_k\subset
\tilde{\sigma}_k$. Hence, by Lemma \ref{inclusions},
$\tilde{\sigma}_{k+1}$ is contained in the projection of $\sigma_k$
onto $B_{m-(k+1)}(\delta_m)$, which in this case equals
$\sigma_{k+1}$ by Lemma \ref{agreeing_of_projections}. Thus
$\tilde{\sigma}_{k+1}\subset\sigma_{k+1}$, as desired.

Now, still assuming that the layer $k$ is thick and that $\alpha_k$ is
an edge disjoint with the boundary or meeting the boundary at a
vertex of defect $\neq 1$, suppose that $\alpha_{k+1}$ is not a
boundary vertex. Let $\overline{x}$ be any vertex in
$\tilde{\sigma}_{k+1}$. Our goal is to prove that $\overline{x}\in
\mathcal{S}(\alpha_{k+1})$. By induction hypothesis we know that
$\mathcal{S}(\alpha_k)\subset \tilde{\sigma}_k$. Since
$\overline{x}$ lies in the layer $k+1$ between $\sigma, \tau$, by Lemma
\ref{layers preserved 0}, we can apply Lemma
\ref{projection_in_geodesic} with $e=\alpha_k$. Hence we get
$\overline{x}\in \mathcal{S}(\alpha_{k+1})$, as desired.

Thus we have completed the induction step in case (1), i.e. for
the layer $k$ thick and $\alpha_k$ an edge disjoint with the boundary or
meeting the boundary at a vertex of defect $\neq 1$.
\medskip

Now suppose that the layer $k$ is thick and $\alpha_k$ is a
non--boundary vertex (case (2)). Then it has two neighbors in the layer
$k+1$ of $\Delta$, suppose first that both of them are at distance
$m-(k+1)$ from $\rho_m$ (we put $\rho_m=v_j$ if $m\geq j$). Then
$\alpha_{k+1}$ is the edge spanned by those two vertices. If it
intersects the boundary, the defect at the boundary vertex is not 1.
Thus we must show that $\mathcal{S}(\alpha_{k+1})$ is contained in
$\tilde{\sigma}_{k+1}$. But by induction hypothesis we know that
$\mathcal{S}(\alpha_k)$ contains $\tilde{\sigma}_k$. Thus, by Lemma
\ref{inclusions}, it is enough to observe that
$\mathcal{S}(\alpha_{k+1})\subset B_{m-(k+1)}(\delta_m)$. This
follows from $\alpha_{k+1}\subset B_{m-(k+1)}(\rho_m)$.

If one of the two neighbors of $\alpha_k$ in the layer $k+1$ is not at
distance $m-(k+1)$ from $\rho_m$, then $\alpha_{k+1}$ is the second
neighbor, it is a non--boundary vertex (unless $k+1=j$, which will
be considered in a moment) and $m<j$. Thus we must show that
$\mathcal{S}(\alpha_{k+1})$ contains $\tilde{\sigma}_{k+1}$. Let
$\bar{z}$ be a vertex in $\tilde{\sigma}_{k+1}$. Then $\bar{z}$ lies
on a 1--skeleton geodesic $\overline{\gamma}$ of length $m-k$ from some vertex of
$\tilde{\sigma_k}\subset \mathcal{S}(\alpha_k)$ to some vertex
$\bar{x}\in \delta_m=\mathcal{S}(\rho_m)$. We claim that if $\rho_m$
is an edge, then the vertex $x=\mathcal{S}^{-1}(\bar{x})\in \Delta$
is the vertex closer to $v_m$ then the other vertex of $\rho_m$.
Indeed, let $y\in \rho_m$ be the vertex closer to $w_m$. Since
$|\alpha_ky|>m-k$ and this distance is realized by a neat geodesic,
hence by Lemma \ref{properties of characteristic surfaces}(i) we
have $|\mathcal{S}(\alpha_k),\mathcal{S}(y)|>m-k$. This proves the
claim. Thus we can apply Lemma \ref{stable_geodesics_in
characteristic image} to $\gamma=\alpha_kx$ and
$\bar{z}\in\overline{\gamma}$, and get $\bar{z}\in
\mathcal{S}(\alpha_{k+1})$, as desired.

Now we come back to the case $k+1=j$ and $\alpha_k$
 a non--boundary vertex. By induction hypothesis we have
 $\tilde{\sigma}_k\subset \mathcal{S}(\alpha_k)$. By Lemma \ref{properties of euclidean
 geodesics}(i) we have that $\sigma_{k+1}=\mathcal{S}(v_{k+1})$ (Remark \ref{delta v_i}) lies in
 $B_{m-(k+1)}(\delta_m)$. Hence, by Lemma \ref{inclusions}, we have that $\tilde{\sigma}_{k+1}$ contains
 $\sigma_{k+1}$, as desired.

Thus we have completed the induction step in case (2), i.e. for
the layer $k$ thick and $\alpha_k$ a non--boundary vertex.
\medskip

Now consider the case that the layer $k$ is thick and $\alpha_k$ is a
boundary vertex or the layer $k$ is thin, but the layer $k+1$ is thick (in
this case put $i=k$). In both cases $\alpha_k=v_k$. First consider
the case that the defect at $v_k$ is $-1$ or the layer $k$ is thin. If
the hypothesis of Lemma \ref{agreeing_of_projections} are not
satisfied, then we can finish as in the previous case (no matter
what is the direction of the inclusion given by the induction
hypothesis) getting $\tilde{\sigma}_{k+1}\subset
\mathcal{S}(\alpha_{k+1})$. Otherwise, $\alpha_{k+1}$ is the edge
spanned by two neighbors of $v_k$ in the layer $k+1$. By Lemma
\ref{properties_of_special_characteristic discs}(i,ii) the defect at
$v_{k+1}$ equals 1. Hence we want to prove that $\sigma_{k+1}$
either contains or is contained in $\tilde{\sigma}_{k+1}$. We know,
by the induction hypothesis, that $\sigma_{k}$ contains or is
contained in $\tilde{\sigma}_{k}$, hence it is enough to use Lemma
\ref{inclusions} and Lemma \ref{agreeing_of_projections}.

Now assume that either the layer $k$ is thick and $\alpha_k$ is a
boundary vertex of defect 0 or an edge intersecting the boundary at
a vertex of defect 1, or the layer $k$ is thin and the layer $k+1$ is also
thin. Similarly, as before, we have that that $\sigma_{k}$ contains
or is contained in $\tilde{\sigma}_{k}$ and we want to prove that
$\sigma_{k+1}$ contains or is contained in $\tilde{\sigma}_{k+1}$.
This follows from Lemma \ref{inclusions} and Lemma
\ref{agreeing_of_projections}.

Thus we have exhausted all the possibilities for case (3) and
completed the induction step. \qed

\section{Euclidean geodesics between simplices of Euclidean geodesics}
\label{CAT(0) geometry of characteristic discs} In this section we
complete the proof of Theorem
\ref{main_theorem_on_euclidean_geodesics}.
Its first ingredient is
Proposition \ref{main_proposition_to_euclidean_geodesics}, proved in
section \ref{Directed geodesics between simplices of Euclidean
geodesics}. The second ingredient is easy 2--dimensional
Euclidean geometry, which we present as a series of lemmas in this
section. Throughout the section, we will be treating characteristic
discs simultaneously as simplicial complexes and $CAT(0)$ metric
spaces.

We start with extending in various ways the notion of a
characteristic disc and surface.


\begin{defin}
\label{generalized characteristic disc} A \emph{generalized
characteristic disc} $\Delta$ for an interval $(i,j)$, where $i<j$,
is a closed $CAT(0)$ (i.e. simply connected) subspace of $\mathbb{E}^2$ with the following
properties. Its boundary is a piecewise linear loop with vertices $v_i,\ldots ,v_j,w_j,\ldots ,w_i,v_i$ (possibly $v_k=w_k$) ,
such that for $i\leq k\leq j$ the straight line segments (or points) $v_kw_k$ are
contained in consecutive parallel lines at distance
$\frac{\sqrt{3}}{2}$. We also require, if $\mathbb{E}^2$ is oriented so that
$v_kw_k$ are horizontal, that $v_k$ lies to the left of $w_k$, or $v_k=w_k$.

A \emph{restriction} of a generalized
characteristic disc to the interval $(l,m)$, where $i\leq l<m\leq
j$, is the generalized characteristic disc enclosed by the loop
$v_l\ldots v_mw_m\ldots w_lv_l$. We will denote it by $\Delta|_l^m$.
If a generalized characteristic disc comes from equipping a systolic
$2$--complex with the standard piecewise Euclidean metric, then we
call it a \emph{simplicial generalized characteristic disc}.
\end{defin}

\begin{rem}
\label{characteristic discs are generalized characteristic discs}
Characteristic discs (resp. modified characteristic discs, c.f.
Definition \ref{euslidean diagonal}) with the standard piecewise
Euclidean metric are simplicial generalized characteristic discs
(resp. generalized characteristic discs).
\end{rem}

\begin{defin}
\label{generalized characteristic discs for incomplete boundary}
Suppose we have simplices $(\sigma_k),(\tau_k)$ in the layer $k$ between
$\sigma,\tau$ (not necessarily the simplices of the directed
geodesics) defined only for $0\leq i\leq k\leq j\leq n$, where
$i<j$, and for $i\leq k<j$ we have that $\sigma_k, \sigma_{k+1}$
span a simplex and $\tau_k,\tau_{k+1}$ span a simplex. Suppose that
for $i\leq k\leq j$ the maximal distance between vertices in
$\sigma_k$ and in $\tau_k$ is $\geq 2$. Then we define a \emph{
partial characteristic disc} and a \emph{partial characteristic
surface} in the following way.

Extend $(\sigma_k),(\tau_k)$ to all $0\leq k \leq n$ so that for
$\sigma_k, \sigma_{k+1}$ and $\tau_k,\tau_{k+1}$ span simplices for
$0\leq k<n$, and $\sigma_0,\tau_0\subset \sigma,
\sigma_n,\tau_n\subset \tau$. (This is possible, since, by example, we may issue directed geodesics from $\sigma_i,\tau_i$ to $\sigma$ and from $\sigma_j,\tau_j$ to $\tau$.) Obviously, $\sigma_k, \tau_k$, lie in the layer $k$ between $\sigma, \tau$ for all $0\leq k\leq n$.
Let $(i_{ext},j_{ext})$ be the thick
interval for extended $(\sigma_k),(\tau_k)$ containing $(i,j)$. Let
$S\colon\Delta\rightarrow X$ be a characteristic surface for
$(i_{ext},j_{ext})$. Then we call $\Delta_{res}=\Delta|_{i}^{j}$ a
partial characteristic disc (which is a simplicial generalized
characteristic disc) and $S_{res}=S|_{\Delta_{res}}$ a partial
characteristic surface.
\end{defin}

\textbf{Caution.} A characteristic surface $S\colon \Delta \rightarrow X$, where $\Delta$ is a characteristic disc for a thick interval $(i,j)$ for $(\sigma_k)_{k=0}^n, (\tau_k)_{k=0}^n$ (as in Definition \ref{characteristic disc and surface}) is not a partial characteristic surface for $(\sigma_k)_{k=i}^j, (\tau_k)_{k=i}^j$. This is because the layers $i,j$ are thin.
But if $i+1<j-1$, then already $S$ restricted to $\Delta|_{i+1}^{j-1}$ is a partial characteristic surface.
\medskip

Next we show that partial characteristic surfaces satisfy most of
the properties of characteristic surfaces. Fix an interval $(i,j)$ and simplices $(\sigma_k)_{k=i}^j, (\tau_k)_{k=i}^j$ as in Definition \ref{generalized characteristic discs for incomplete boundary}. Let $S_{res}\colon
\Delta|_{res}\rightarrow X$ be a partial characteristic surface, as
above.

\begin{lemma}
\label{properties of generalized characteristic discs}
\item (i) $\Delta_{res}$ (and thus $S_{res}$) is flat,
\item (ii) if we embed $\Delta_{res}\subset \mathbb{E}_{\Delta}^2$,
then $v_iw_i$ and $v_jw_j$ are parallel and the consecutive layers
between them are contained in consecutive straight lines parallel
to $v_iw_i$ and $v_jw_j$.
\item(iii) $S_{res}$ is an isometric embedding on 1--skeleton
of a subcomplex spanned by any pair of consecutive layers between
$v_iw_i$ and $v_jw_j$ in $\Delta_{res}$.
\item (iv) $\Delta_{res}\subset \mathbb{E}_{\Delta}^2$ does not depend on the choice
of $\sigma_k,\tau_k$ for $k<i$ and $k>j$, the choice of $s_k,t_k$
for $0\leq k\leq n$, and the choice of $S$.
\item If we have two partial characteristic surfaces
$S_1\colon\Delta_1\rightarrow X, S_2\colon\Delta_2\rightarrow X$,
then after identifying partial characteristic discs $\Delta_1=\Delta_2$
(which is possible by (ii)) we have that
\item (v) for any vertices $x,y\in \Delta_1=\Delta_2$ at distance $1$, $S_1(x)$ and $S_2(y)$ are
also at distance $1$,
\item (vi) for any vertex $x\in \Delta_1=\Delta_2$, $S_1(x)$ and $S_2(x)$ are at distance $\leq 1$.
\item(vii) $S(v_kw_k)$ lies in the layer $k$ between $\sigma$ and $\tau$.
\end{lemma}

\proof Assertions (i) and (ii) follow immediately from Lemma
\ref{properties_of_characteristic discs}(i,ii). Assertion (iii)
follows from Lemma \ref{properties of characteristic surfaces}(i).
To prove (iv) notice that $\Delta_{res}=\Delta|_i^j$ is determined by the
distances $|s_kt_k|$ for $i\leq k\leq j$ and $|s_kt_{k+1}|$ for
$i\leq k<j$, by (iii). Hence, if we fix $s_k$ and $t_k$ for $i\leq k\leq j$, then $\Delta_{res}$ does not depend on the extension of $(\sigma_k)_{k=i}^j,(\tau_k)_{k=i}^j$. On the other hand, if we fix such an extension, then $|s_kt_k|, |s_kt_{k+1}|$ do not depend on the choice of $s_k, t_k$, by Lemma \ref{properties of characteristic surfaces}(ii).

It is a bit awkward to try to obtain assertion (v) as a consequence of Lemma \ref{properties of characteristic surfaces}(iii). Let us say, instead, that assertion (v) follows immediately from the proof of Lemma \ref{properties of characteristic surfaces}(iii). Similarly, assertion (vi) follows from the proof of Lemma \ref{properties of characteristic surfaces}(iv).
Assertion (vii) follows directly from Remark \ref{layers preserved 0}.\qed

\begin{defin}
\label{def_of_generalized_S} We define the \emph{partial
characteristic image} $\mathcal{S(\rho)}$ of a simplex $\rho$ in the
partial characteristic disc as the span of  $S(\rho)$ over all
partial characteristic surfaces $S$. By Lemma \ref{properties of
generalized characteristic discs}(v,vi), $\mathcal{S(\rho)}$ is a
simplex. We call this assignment the \emph{partial
characteristic mapping}. Like in Definition \ref{characteristic image} we can consider also the assignment $\mathcal{S}^{-1}$.
\end{defin}

\begin{defin}
\label{def_of_close} Let $\Delta$ be a generalized characteristic
disc and $\gamma$, $\gamma'$ be two paths connecting some points on
$v_iw_i$ to points on $v_jw_j$ such that intersections of $\gamma$,
$\gamma'$ with $v_kw_k$ are unique for each $i\leq k\leq j$. We say
that $\gamma, \gamma'$ are \emph{$d$--close} if they intersect
$v_kw_k$ in points at distance $\leq d$ for each $i\leq k\leq j$.
\end{defin}

The following lemma describes the possible displacements of $CAT(0)$
geodesics in characteristic discs when perturbing the boundary and
the endpoints.

\begin{lemma}
\label{geodesics stable under perturbing boubdary} Let
$\Delta'\subset \Delta$ be two generalized characteristic discs for
$(i,j)$ such that for each $i\leq k \leq j$ we have $v'_kw_k'\subset
v_kw_k$ (and the order is $v_kv'_kw'_kw_k$) and $|v_kv'_k|\leq d,\
|w_kw'_k|\leq d$. Then for any points $x\in v_iw_i, \ y\in v_jw_j,\
 x'\in v'_iw'_i, \ y'\in v'_jw'_j$ such that $|xx'|\leq d,\ |yy'|\leq d$,
the $CAT(0)$ geodesics from $x$ to $y$ in $\Delta$ and from $x'$ to
$y'$ in $\Delta'$ are $d$--close in $\Delta$.
\end{lemma}
\proof Denote by $\gamma, \gamma'$ the geodesics from $x$ to $y$ in
$\Delta$ and from $x'$ to $y'$ in $\Delta'$ respectively. Denote by
$N_d(\gamma)$ the set of points in $\Delta$ at distance $\leq d$
from $\gamma$ in the direction parallel to $v_kw_k$ (i.e. the intersection with $\Delta$ of the union of translates of $\gamma$ by a distance $\leq d$ in the direction parallel to $v_kw_k$), and by
$N'_d(\gamma)$ the intersection $N_d(\gamma)\cap \Delta'$.

Observe that $N'_d(\gamma)$ is connected, since for each $k$ the set
$v'_kw'_k\cap N'_d(\gamma)$ is nonempty and the intersection of
$N'_d(\gamma)$ with each of the parallelograms
$v'_kw'_kw'_{k+1}v'_{k+1}$ is an intersection of two parallelograms,
hence convex and connected. We claim that $N_d(\gamma)$ is convex in
$\Delta$. To establish this, we need to study the interior angle at
vertices of $\partial N_d(\gamma)$ outside $\partial \Delta$. The
only possibility for angle $>180^{\circ}$ is at the horizontal
translates of break points of $\gamma$. But since $\gamma$ is a
$CAT(0)$ geodesic, then each of its break points lies on the boundary
of $\Delta$, and the translate, for which possibly the angle is
$>180^{\circ}$, lies outside $\Delta$. Thus the claim follows. Hence
(by connectedness) $N'_d(\gamma)$ is convex in $\Delta'$. Thus
$\gamma'\subset N'_d(\gamma)$ and we are done. \qed

\medskip
Let us prepare the setting for the next lemma. It will help us deal
with the data given by Proposition
\ref{main_proposition_to_euclidean_geodesics}, which is, roughly
speaking, a pair of surfaces spanned on nearby pairs of geodesics.
To be more precise, let
$\hat{\sigma}_k,\hat{\tau}_k,\tilde{\sigma}_k,\tilde{\tau}_k$ be
simplices in the layers $i\leq k\leq j$ between $\sigma, \tau$
satisfying conditions of Definition \ref{generalized characteristic
discs for incomplete boundary}. Moreover, assume that for each $i\leq k\leq j$ we have that
$\hat{\sigma}_k\subset\tilde{\sigma}_k$ or $\tilde{\sigma}_k\subset\hat{\sigma}_k$, and $\hat{\tau}_k\subset\tilde{\tau}_k$ or $\tilde{\tau}_k\subset\hat{\tau}_k$.
Let $\hat{\Delta}, \tilde{\Delta}$
be associated partial characteristic discs, unique by
\ref{properties of generalized characteristic discs}(iv).
Denote the boundary vertices of $\hat{\Delta}$ (resp. $\tilde{\Delta}$) by $\hat{v}_k,\hat{w}_k$ (resp. $\tilde{v}_k,\tilde{w}_k$), its characteristic mapping by $\hat{\mathcal{S}}$ (resp. $\tilde{\mathcal{S}}$).

\begin{lemma}
\label{common_characteristic_disc} There exists a simplicial
generalized characteristic disc $\overline{\Delta}$ for $(i,j)$
and embeddings (thought of as inclusions, for simplicity)
$\overline{\Delta}\subset \hat{\Delta},\overline{\Delta}\subset
\tilde{\Delta}$ such that the distances $|\overline{v}_k\hat{v}_k|,$
$|\overline{w}_k\hat{w}_k|$ in $\hat{\Delta}$ and distances
$|\overline{v}_k\tilde{v}_k|,|\overline{w}_k\tilde{w}_k|$ in
$\tilde{\Delta}$ are all $\leq 1$ for $i\leq k\leq j$. Moreover,
$|\overline{v}_k\overline{w}_k|\geq 1$ for $i\leq k\leq j$.
\end{lemma}
\proof For each $i\leq k\leq j$, let $\sigma^{max}_k$ be the greater
among $\hat{\sigma_k},\tilde{\sigma_k}$ and let $\sigma^{min}_k$ be
the smaller, let $\tau^{max}_k$ be the greater among
$\hat{\tau_k},\tilde{\tau_k}$ and let $\tau^{min}_k$ be the smaller.
Pick vertices $x_k\in \sigma^{max}_k, y_k\in \tau^{max}_k$ so that
the distance $|x_ky_k|$ is maximal. If possible, choose them from
$\sigma^{min}_k, \tau^{min}_k$ (if it is possible for $x_k,y_k$
independently, then it is possible for both of them at the same
time, by Lemma
\ref{realizing_thickness_in_pairs_implies_realizing_as_a_pair}).
Pick a 1--skeleton geodesic $\phi_k$ connecting $x_k$ to $y_k$
intersecting $\sigma^{min}_k, \tau^{min}_k$ (this is possible by
Corollary \ref{corollary_to_no_trapezoid}).
If $x_k\in \sigma^{min}_k$, then put $\overline{s}_k=x_k$, otherwise
let $\overline{s}_k$ be the neighbor of $x_k$ on $\phi_k$.
Analogically, if $y_k\in \tau^{min}_k$, then put
$\overline{t}_k=y_k$, otherwise let $\overline{t}_k$ be the neighbor vertex
of $y_k$ on $\phi_k$. Thus $\overline{s}_k\in
\sigma^{min}_k,\overline{t}_k\in \tau^{min}_k$. Let
$\overline{\Delta}$ be the partial characteristic disc for
$(\overline{s}_k),(\overline{t}_k)$ for $i\leq k\leq j$. Denote its boundary vertices by $\overline{v}_k,\overline{w}_k$.

The embedding, say $\overline{\Delta}\subset \hat{\Delta}$, is
defined as follows. By Proposition \ref{stable_geodesics_in
characteristic image} there exists a characteristic surface
$\overline{S}\colon \overline{\Delta} \rightarrow X$ such that
$\overline{S}(\overline{v}_k\overline{w}_k)\subset \phi_k$.
Moreover, again by Proposition \ref{stable_geodesics_in
characteristic image}, the sub--geodesic
$\overline{s}_k\overline{t}_k$ of $\phi_k$ lies in
$\hat{\mathcal{S}}(\hat{\Delta})$. Hence we can define the desired
mapping as the composition $\hat{\mathcal{S}}^{-1}\circ
\overline{S}$. To check that this is an embedding it is enough to
check that it preserves the the layers (Lemma \ref{properties of generalized characteristic discs}(vii)) and is isometric on the layers (Lemma \ref{properties of generalized characteristic discs}(iii)).

To prove the last assertion fix $k$ and assume w.l.o.g. that
$\sigma^{min}_k=\hat{\sigma_k}$. Then
$|\overline{v}_k\overline{w}_k|\geq |\hat{v}_k\hat{w}_k|-1\geq 1$, as
desired. \qed

\medskip
Now let us prepare the statement of our final lemma. One can view it
as a simple case of Theorem
\ref{main_theorem_on_euclidean_geodesics}, case of $X$ being flat.

Let $\Delta$ be a characteristic disc for a thick interval $(i,j)$
for the directed geodesics $(\sigma_k),(\tau_k)$ between
$\sigma,\tau$ and let $\gamma'$ be its $CAT(0)$ diagonal, c.f. Definition
\ref{euslidean diagonal}. Let $(\rho_k)_{k=i+1}^{j-1}$ be simplices
of the Euclidean diagonal in $\Delta$ (Definition \ref{euslidean
diagonal}). Fix $i\leq l<m\leq j$. If $i< l<m< j$ then let
$(\alpha_k)_{k=l}^m,(\beta_k)_{k=m}^l$ be directed geodesics in
$\Delta$ from $\rho_l$ to $\rho_m$ and from $\rho_m$ to $\rho_l$
respectively. If $l=i$ then put $\rho_i=v_i$ in the definition of
$(\alpha_k)_{k=l}^m$ and $\rho_i=w_i$ in the definition of
$(\beta_k)_{k=m}^l$. If $m=j$ then put $\rho_j=w_j$ in the
definition of $(\beta_k)_{k=m}^l$ and $\rho_j=v_j$ in the definition
of $(\alpha_k)_{k=l}^m$. For all other purposes we will put
$\rho_i=v_iw_i,\rho_j=v_jw_j$.

Let $\bigcup\hat{\Delta}$ be the subcomplex of $\Delta$ which is the
span of the union of $\conv\{\alpha_k,\beta_k\}$ over all $l\leq
k\leq m$. Note that $\bigcup\hat{\Delta}$ is a simplicial
generalized characteristic disc. Denote the vertices of its boundary loop by $(\hat{v}_k)$ and $(\hat{w}_k)$.
Denote by $\hat{\gamma}$ the $CAT(0)$
geodesic joining in $\bigcup\hat{\Delta}$ the barycenters of
$\rho_l$ and $\rho_m$ (which lie in in $\bigcup\hat{\Delta}$).

\begin{lemma}
\label{geodesic of a subdisc agrees} $\gamma'$ restricted to
$\Delta|_l^m$ and $\hat{\gamma}$ are $\frac{1}{2}$--close in
$\Delta|_l^m$ .
\end{lemma}

\proof Let us denote by $\bigcup\hat{\Delta}_0$ the generalized
characteristic disc obtained from $\bigcup\hat{\Delta}$ by removing
the following triangles: For any boundary vertex of defect 1 in
the layers $\neq l,m$, say $\hat{v}_k$, cut off a triangle along the
segment $\hat{v}_{k-1}\hat{v}_{k+1}$. For any boundary vertex of defect 2 (which is possible in
the layers $l,m$), say $\hat{v}_l$, cut off a triangle joining
$\hat{v}_{l+1}$ to the barycenter of $\hat{v}_l\hat{w}_l$.

We claim that $\bigcup\hat{\Delta}_0$ is convex in $\Delta$ (treated
as $CAT(0)$ spaces). This means that at all vertices of $\partial
\bigcup\hat{\Delta}_0$ outside $\partial \Delta$, the interior angle
of $\bigcup\hat{\Delta}_0$ is $\leq 180^\circ$. We skip the
proof, which is an easy consequence of Lemma \ref{directed geodesics in characteristic discs}.

Let $\hat{\gamma}_0$ be the $CAT(0)$ geodesic in
$\bigcup\hat{\Delta}_0$ joining the barycenter $\hat{x}$ of $\rho_l$
with the barycenter $\hat{y}$ of $\rho_m$ (observe that
$\hat{x},\hat{y}\in\bigcup\hat{\Delta}_0$). Since
$\bigcup\hat{\Delta}_0\subset \Delta$ is convex, $\hat{\gamma}_0$
agrees with the $CAT(0)$ geodesic in $\Delta$ joining
$\hat{x},\hat{y}$.

Now we apply Lemma \ref{geodesics stable under perturbing
boubdary} to $\Delta'|_l^m\subset \Delta|_l^m$ (c.f. Definition
\ref{euslidean diagonal} for the definition of $\Delta'$), and
geodesics $\hat{\gamma}_0$ in $\Delta|_l^m$ and $\gamma'$ restricted
to $\Delta'|_l^m$. Observe that endpoints $\hat{x},\hat{y}$ of
$\hat{\gamma}_0$ are at distance $\leq \frac{1}{2}$ from
$\gamma'\cap v_lw_l,\gamma'\cap v_mw_m$ by the definition of
$\rho_l, \rho_m$. Hence, by Lemma \ref{geodesics stable under
perturbing boubdary}, we have that $\hat{\gamma}_0$ is
$\frac{1}{2}$--close to $\gamma'$ restricted to $\Delta|_l^m$.

Now observe that since $\bigcup\hat{\Delta}_0$ is also convex in
$\bigcup\hat{\Delta}$, we have $\hat{\gamma}_0=\hat{\gamma}$ and we
are done. \qed
\medskip

Finally, we can proceed with the following.

\medskip\par\noindent\textbf{Proof of Theorem \ref{main_theorem_on_euclidean_geodesics}.}\ignorespaces
\ First suppose that the layer $k$ for $(\sigma_t), (\tau_t)$ is thin.
Then, by Proposition \ref{main_proposition_to_euclidean_geodesics}(i),
$\tilde{\sigma}_k$ contains or is contained in $\sigma_k$ and
$\tilde{\tau}_k$ contains or is contained in $\tau_k$. Hence the
thickness of the layer $k$ for $(\tilde{\sigma}_t), (\tilde{\tau}_t)$ is
$\leq 3$ and thus $\tilde{\sigma_k}\subset B_1(\tilde{\delta}_k)$ or
$\tilde{\tau_k}\subset B_1(\tilde{\delta}_k)$, hence
$|\tilde{\delta}_k,\delta_k|\leq 1$.

Now suppose that the layer $k$ for $(\sigma_t), (\tau_t)$ is thick and
suppose it is contained in a thick interval $(i,j)$ with a
characteristic disc $\Delta$. Put $\rho_l=v_iw_i$ if $l\leq i$
and $\rho_m=v_jw_j$ if $m\geq j$. We will use the notation
introduced before Lemma \ref{geodesic of a subdisc agrees}. First
suppose that the layer $k$ for $(\tilde{\sigma}_t), (\tilde{\tau}_t)$ is
thin. Then, by Proposition
\ref{main_proposition_to_euclidean_geodesics}(ii), the maximal distance
between vertices in $\mathcal{S}(\alpha_k)$ and $\mathcal{S}(\beta_k)$, hence in $\alpha_k$ and $\beta_k$ is $\leq 3$. Since
$\hat{\gamma}\cap v_kw_k$ lies in $\conv\{\alpha_k,\beta_k\}$, Lemma
\ref{geodesic of a subdisc agrees} implies that $\gamma'\cap v_kw_k$
is at distance $\leq \frac{1}{2}$ from $\conv\{\alpha_k,\beta_k\}$. Hence
$\alpha_k\subset B_1(\rho_k)$ or $\beta_k\subset B_1(\rho_k)$. Thus
$\tilde{\delta}_k,\delta_k$ are at distance $\leq 1$.

Now suppose that the layer $k$ for $(\tilde{\sigma}_t),
(\tilde{\tau}_t)$ is thick. Let $\tilde{\Delta}$ be the
characteristic disc for the thick interval $(\tilde{i},\tilde{j})$
containing $k$ for $(\tilde{\sigma}_t), (\tilde{\tau}_t)$. If the layer
$k$ for $(\alpha_t),(\beta_t)$ (between $\rho_l,\rho_m$ in $\Delta$)
is thin, then the thickness of the layer $k$ for $(\tilde{\sigma}_t),
(\tilde{\tau}_t)$ is $\leq 3$, by Proposition
\ref{main_proposition_to_euclidean_geodesics}(ii). Hence
$\tilde{\sigma}_k\subset B_1(\tilde{\delta}_k)$ or
$\tilde{\tau}_k\subset B_1(\tilde{\delta}_k)$. By Lemma
\ref{geodesic of a subdisc agrees} we have $|\rho_k,\alpha_k|\leq 1$
and $|\rho_k,\beta_k|\leq 1$, hence altogether
$|\tilde{\delta}_k,\delta_k|\leq 2$.
\medskip

So suppose that the layer $k$ for $(\alpha_t),(\beta_t)$ in $\Delta$ is
thick, let $\hat{i},\hat{j}$ be the thick interval for
$(\alpha_t),(\beta_t)$ containing $k$ and $\hat{\Delta}$ the
corresponding characteristic disc. Observe that
$\hat{\Delta}=\bigcup\hat{\Delta}|_{\hat{i}}^{\hat{j}}$. Let
$i_{max}$ be the maximum of $\hat{i},\tilde{i}$ and $j_{min}$ be the
minimum of $\hat{j},\tilde{j}$. Obviously $i_{max}<k<j_{min}$.
Assume $i_{max}+1<j_{min}-1$, in the case of equality the argument
is similar and we omit it.

By Proposition \ref{main_proposition_to_euclidean_geodesics}(ii) we can
apply Lemma \ref{common_characteristic_disc} to $\tilde{\Delta}$ and
$\hat{\Delta}$ restricted to $(i_{max}+1, j_{min}-1)$. Denote by
$\overline{\Delta}$ the simplicial generalized characteristic disc
for $(i_{max}+1, j_{min}-1)$ guaranteed by Lemma
\ref{common_characteristic_disc}. Denote by  $\overline{\Delta}'$
the generalized characteristic disc obtained from
$\overline{\Delta}$ by removing horizontal (the direction of
$\overline{v}_t\overline{w}_t$) $\frac{1}{2}$--neighborhood of the
boundary, which is allowed since $|\overline{v}\overline{w}|\geq 1$
by Lemma \ref{common_characteristic_disc}. Let $\tilde{\Delta}'$ be
the modified characteristic in $\tilde{\Delta}$ and
$\tilde{\gamma}'$ the $CAT(0)$ diagonal of $\tilde{\Delta}$ (c.f.
Definition \ref{euslidean diagonal}). Define a generalized
characteristic disc $\hat{\Delta}'$ and a $CAT(0)$ geodesic
$\hat{\gamma}'$ in $\hat{\Delta}'$ in the following way. For each
$\hat{i}\leq t\leq \hat{j}$ denote by $\hat{v}'_t, \hat{w}'_t$
points on $\hat{v}_t\hat{w}_t$ at distance $\frac{1}{2}$ from
$\hat{v}_t,\hat{w}_t$, respectively, if $\hat{v}_t\neq \hat{w}_t$.
Otherwise, put $\hat{v}'_t=\hat{v}_t,\hat{w}'_t=\hat{w}_t$. Let
$\hat{\Delta}'$ be the generalized characteristic disc enclosed by
the loop
$\hat{v}'_{\hat{i}}\ldots\hat{v}'_{\hat{j}}\hat{w}'_{\hat{j}}\ldots\hat{w}'_{\hat{i}}\hat{v}'_{\hat{i}}$.
Let $\hat{\gamma}'$ be the $CAT(0)$ geodesic in $\hat{\Delta}'$
joining $\hat{v}'_l= \hat{w}'_l$ and $\hat{v}'_m= \hat{w}'_m$. By
Lemma \ref{common_characteristic_disc} we have inclusions of
$\overline{\Delta}'$ into $\hat{\Delta}'|_{i_{max}+1}^{j_{min}-1},
\tilde{\Delta}'|_{i_{max}+1}^{j_{min}-1}$ with distances
$|\overline{v}'_t\hat{v}'_t|,|\overline{w}'_t\hat{w}'_t|$ in
$\hat{\Delta}'$, and distances
$|\overline{v}'_t\tilde{v}'_t|,|\overline{w}'_t\tilde{w}'_t|$ in
$\tilde{\Delta}'$ all $\leq 1$ for ${i_{max}+1}\leq t\leq
{j_{min}-1}$.

Now we will choose a special point $\overline{x}\in
\overline{v}'_{i_{max}+1}\overline{w}'_{i_{max}+1}$.
W.l.o.g. assume that $i_{max}=\tilde{i}$, hence
$|\tilde{v}_{i_{max}+1}\tilde{w}_{i_{max}+1}|=2$. Choose any
$\hat{x}$ in $\overline{v}'_{i_{max}+1}\overline{w}'_{i_{max}+1}$ at
distance $\leq 1$ from $\hat{\gamma}'$, which is possible, since
$|\overline{v}'_{i_{max}+1}\hat{v}'_{i_{max}+1}|\leq 1$ and
$|\overline{w}'_{i_{max}+1}\hat{w}'_{i_{max}+1}|\leq 1$. Since
$|\tilde{v}'_{i_{max}+1}\tilde{w}'_{i_{max}+1}|\leq 1$,
$\overline{x}$ is also at distance $\leq 1$ from $\tilde{\gamma}'$.
Choose $\overline{y}$ in
$\overline{v}'_{j_{min}-1}\overline{w}'_{j_{min}-1}$ in an analogous
way.

By this construction the endpoints of  $\tilde{\gamma}'$ and
$\hat{\gamma}'$ restricted to $(i_{max}+1,j_{min}-1)$ are at
distance $\leq 1$ from $\bar{x},\bar{y}$ in
$\tilde{\Delta}'|_{i_{max}+1}^{j_{min}-1},\hat{\Delta}'|_{i_{max}+1}^{j_{min}-1}$,
respectively. Thus, using twice Lemma \ref{geodesics stable under
perturbing boubdary}, we get that $\tilde{\gamma}'$ and
$\hat{\gamma}'$ restricted to $(i_{max}+1,j_{min}-1)$ are 1--close
to the $CAT(0)$ geodesic $\bar{x}\bar{y}$ in $\overline{\Delta}'$ (in
$\tilde{\Delta}'|_{i_{max}+1}^{j_{min}-1},\hat{\Delta}'|_{i_{max}+1}^{j_{min}-1}$
respectively).

By Lemma \ref{geodesic of a subdisc agrees}, $\gamma'$ and
$\hat{\gamma}$ are $\frac{1}{2}$--close in $\Delta|_l^m$. By Lemma
\ref{geodesics stable under perturbing boubdary}, $\hat{\gamma}'$
and $\hat{\gamma}$ are $\frac{1}{2}$--close in
$\hat{\Delta}|_{\hat{i}}^{\hat{j}}$. Putting those four estimates
together we get that $\delta_k, \tilde{\delta}_k$ are at distance
$\leq 3$, as desired. \qed
\medskip

We end this section by indicating, how Theorem \ref{main_theorem_on_euclidean_geodesics} can be promoted to Theorem \ref{2}, with a reasonable constant $C$. The difference in statements comes from substituting $\delta_l,\delta_m$ with $x\in\delta_l,\ y\in\delta_m$ such that $|xy|=m-l$. As a first step, we check that Proposition \ref{main_proposition_to_euclidean_geodesics} implies that the directed geodesics between $x$ and $y$ lie near the union of characteristic images of characteristic discs for $(\sigma_k),(\tau_k)$. This follows from the fact that directed geodesics in systolic complexes satisfy the so called fellow traveler property with a good constant, see \cite{JS}. The second step is to reprove Lemma \ref{common_characteristic_disc} allowing $\hat{\sigma}_k$ and $\tilde{\sigma}_k$ (and similarly $\hat{\tau}_k$ and $\tilde{\tau}_k$) to be farther apart, at distance bounded by the above fellow traveler constant. Then some minor changes in the proof of Theorem \ref{main_theorem_on_euclidean_geodesics} yield Theorem \ref{2}.

We will give a different complete proof of Theorem \ref{2} (though with a worse constant) in the next section.

\section{Characteristic discs spanned on Euclidean geodesics}
\label{Characteristic discs for Euclidean geodesics}

In this section we prove the following crucial proposition, which,
roughly speaking, says that in a characteristic disc spanned on a
Euclidean geodesic and an arbitrary other geodesic, the boundary
segment corresponding to the Euclidean geodesic is coarsely a $CAT(0)$
geodesic. We introduce the following notation, which will be fixed
for the whole section.

\medskip
Let $\sigma,\tau$ be simplices in a systolic complex $X$ satisfying as before $\sigma\subset
S_n(\tau), \tau\subset S_n(\sigma)$ and suppose
$(p_k)_{k=0}^n,(r_k)_{k=0}^n$ are 1--skeleton geodesics with
endpoints in $\sigma$ and $\tau$ such that $r_k\in \delta_k$, where
$(\delta_k)_{k=0}^n$ is the Euclidean geodesic between $\sigma$ and
$\tau$. Let $0\leq i_{pr}<j_{pr}\leq n$ be a thick interval for
$(p_k),(r_k)$ and let $\Delta_{pr}, \mathcal{S}_{pr}$ be the
corresponding characteristic disc and mapping. Let $\gamma_{pr}$ be
the $CAT(0)$ geodesic in $\Delta_{pr}$ joining the barycenters of the
unique edges in the layers $i_{pr},j_{pr}$.

\begin{proposition}
\label{Euclidean_near_CAT(0)}
$\gamma_{pr}$ is $99$--close to the boundary path $\mathcal{S}^{-1}_{pr}\big((r_k)\big)$.
\end{proposition}

This proposition has fundamental consequences. One of them is Theorem \ref{3}, which says roughly this:
in a "Euclidean geodesic triangle", the distance between the
midpoints of two sides is, up to an additive constant, smaller than half of the
length of the third side. We study this in the next section.

The second consequence of Proposition \ref{Euclidean_near_CAT(0)} is an alternative proof of the following.
\begin{theorem}[Theorem \ref{2}]
\label{corollary_to_main_theorem_on_euclidean_geodesics}
Let $\sigma,\tau$ be simplices of a systolic complex $X$, such that for some natural $n$ we have
$\sigma \subset S_n(\tau), \tau \subset S_n(\sigma)$.
Let $(\delta_k)_{k=0}^n$ be the Euclidean geodesic between $\sigma$ and $\tau$. Take some $0\leq l<m\leq n$ and let $(r_k)_{k=l}^m$ be a 1--skeleton geodesic such that $r_k\in \delta_k$
for $l\leq k\leq m$.
Consider the simplices $\tilde{\delta}_l=r_l,\tilde{\delta}_{l+1},\ldots, \tilde{\delta}_m=r_m$
of the Euclidean geodesic between
vertices $r_l$ and $r_m$. Then for each $l\leq k\leq m$ we have $|\delta_k,\tilde{\delta}_k|\leq C$, where $C$ is a universal constant.
\end{theorem}

\proof
Extend $(r_k)_{k=l}^m$ to a 1--skeleton geodesic $(r_k)_{k=0}^n$ between $\sigma,\tau$ so that $r_k\in \delta_k$ (this is possible by Lemma \ref{properties of euclidean geodesics}(i)). Let $(\tilde{r}_k)_{k=l}^m$ be any 1--skeleton geodesic between $r_l$ and $r_m$ such that $\tilde{r}_k\in \tilde{\delta}_k$.
Put additionally $\tilde{r}_k=r_k$ for $0\leq k< l$ and for $m< k\leq n$. Let $\Delta_{r\tilde{r}}$ be the characteristic disc for some thick interval for $(\tilde{r}_k)_{k=0}^n, (r_k)_{k=0}^n$ and
let $\gamma_{r\tilde{r}}$ be the $CAT(0)$ geodesic joining the barycenters of its outermost edges.
Let $\mathcal{S}_{r\tilde{r}}$ be the corresponding characteristic mapping.

Notice that $\Delta_{r\tilde{r}}$ is also a characteristic disc for $(r_k)_{k=l}^m,(\tilde{r}_k)_{k=l}^m$ between
$r_l$ and $r_m$.
Applying twice Proposition \ref{Euclidean_near_CAT(0)} we obtain that $\gamma_{r\tilde{r}}$ is $99$--close to both
$\mathcal{S}^{-1}_{r\tilde{r}}\big((r_k)\big)$ and $\mathcal{S}^{-1}_{r\tilde{r}}\big((\tilde{r}_k)\big)$. This proves that for all $l\leq k\leq m$ we have $|r_k\tilde{r}_k|\leq 198$, hence $|\delta_k,\tilde{\delta_k}|\leq 198$. Thus any $C\geq 198$ satisfies the assertion of the theorem. \qed
\medskip

The proof of Proposition \ref{Euclidean_near_CAT(0)} is rather
technical. This is the reason we decided to present the
straightforward proof of Theorem
\ref{main_theorem_on_euclidean_geodesics} (the weak version of Theorem \ref{2}) via Proposition
\ref{main_proposition_to_euclidean_geodesics}. Before we get into
technical details of the proof, split into various lemmas, we
present an outline, which hopefully helps to keep track of the main
ideas.

\medskip\par\noindent\textbf{Outline of the proof of Proposition \ref{Euclidean_near_CAT(0)}.}\ignorespaces
\ We are dealing with configurations of four geodesics between $\sigma$
and $\tau$: the directed geodesics, which we denote by
$(\sigma_k)_{k=0}^n,(\tau_k)_{k=0}^n$, as in the previous sections,
$(r_k)_{k=0}^n$, which goes along the Euclidean geodesic $\delta_k$, and the
fourth arbitrary 1--skeleton geodesic $(p_k)_{k=0}^n$.
For the layer $k$ thick (for $(\sigma_k),(\tau_k)$) we have that $\delta_k=\mathcal{S}(\rho_k)$, where $\rho_k$ is the simplex of the
Euclidean diagonal in appropriate characteristic disc $\Delta$ for $(\sigma_k),(\tau_k)$.
Hence we need to
find out, what is the possible position of $(p_k)$ w.r.t.
$\mathcal{S}(\Delta)$. It turns out that in each layer there are 1--skeleton geodesics between simplices
$\sigma_k,\tau_k$ and $p_k$, which form a very thin triangle
(Lemma \ref{tripods}). The intersection with $\mathcal{S}(\Delta)$ of the center simplex of this triangle will be later denoted by $\overline{\chi}_k$.

In Lemma \ref{preparation_to_triplanes} we study, how do
$\overline{\chi}_k$ vary with $k$. Assume for simplicity that $p_k$ stay away from $\mathcal{S}(\Delta)$.
Then it turns out that first (i.e. for small $k$) $\overline{\chi}_k$ follow $\mathcal{S}(w_k)$, next the
barycenters of $\overline{\chi}_k$ lie in the characteristic image
of a vertical line in $\Delta$ and last $\overline{\chi}_k$ follow $\mathcal{S}(v_k)$. The $CAT(0)$ diagonal $\gamma'$ of $\Delta$
crosses this line at most once. Thus we can divide each "thick"
interval (an interval with all layers thick, in opposition to the thick interval with thin endpoint layers) for $(\sigma_k),(\tau_k)$ into three subintervals: the
"initial" one, for which
${\chi}_k=\mathcal{S}^{-1}(\overline{\chi}_k)$ is far to the right
from $\rho_k$ or near $w_k\in \partial \Delta$, the "middle" one, for
which ${\chi}_k$ is near
$\rho_k$, and the "final" one, for which
${\chi}_k$ is far to the left from
$\rho_k$ or near $v_k\in \partial \Delta$, see Lemma
\ref{triplanes}. Moreover, in the "initial" (resp. "final") interval we can distinguish a "pre--initial" (resp. "post--final") interval in which $\gamma'$ stays away from $w'_k\in \partial\Delta'$ (resp. $v'_k\in \partial\Delta'$), where $\Delta'$ is the modified characteristic disc. This distinction is done in the main body of the proof of Proposition \ref{Euclidean_near_CAT(0)}.
The vertices $\mathcal{S}^{-1}_{pr}(r_k)$ in
$\Delta_{pr}$, for $k$ in one of these intervals, are positioned as follows.
The vertices of the "middle" interval together with the vertices of the other ones outside the "pre--initial" and "post--final" intervals form a coarse vertical line (this is a consequence of Lemma \ref{joining 0}),
while the vertices of the "pre--initial" and "post--final" intervals form also coarse $CAT(0)$
geodesics,
fortunately forming with the coarse vertical line angles
$\geq 180^{\circ}$ at the endpoints. This proves Proposition \ref{Euclidean_near_CAT(0)} in the simple case of
a single "thick" interval for $(\sigma_k),(\tau_k)$.

In the complex case, the question is, how may the various "thick" intervals and thin layers
for $(\sigma_k),(\tau_k)$ alternate. We define roughly the following notions.
A "thin" interval is an interval of not very thick layers.
A "proper thin" interval is a "thin" interval with thin layers at the beginning and at the end.
A "very thick" interval is an interval containing a layer that is very thick.
In Lemma \ref{thinvertical} we
prove that the vertices $\mathcal{S}^{-1}_{pr}(r_k)$, for $k$ in
a "thin" interval, form a coarse vertical
line.
In Corollary \ref{joining} we prove that if at the beginning of
a thin layer there is an adjoined "thick" interval, then this
"thick" interval has the "final" subinterval constructed above "thin".
Similarly, if at the end of a thin layer there is an adjoined
"thick" interval, then this thick interval has the "initial"
subinterval "thin".
The last piece of the puzzle is an assertion in Lemma
\ref{triplanes}, that for a "very thick" interval, either its "initial" or "final" subinterval is non--"thin".

The way to put these pieces together is the following. We take a
maximal "proper thin" interval. The "very thick" interval
adjoined at the beginning of this "proper thin" interval must have
either the "initial" or the "final" subinterval non--"thin" (Lemma \ref{triplanes}), but the possibility
of the "final" subinterval non--"thin" is excluded (Corollary \ref{joining}). Thus its "initial"
subinterval is non--"thin" and this excludes the possibility that some
thin layer (hence any layer) is adjoined at the beginning of this "very thick" interval (Corollary \ref{joining}). We can apply analogous considerations to the
"very thick" interval adjoined at the end of the "proper thin"
interval. Altogether, we have the following configuration: the "proper thin"
interval, with a "very thick" interval with "thin"
"final" subinterval adjoined at the beginning, and with a "very thick"
interval with "thin" "initial" subinterval adjoined at the
end. Moreover, in the first of the "very thick" intervals we distinguish the "pre--initial" interval and in the second one we distinguish the "post--final" interval.
The vertices
$\mathcal{S}^{-1}_{pr}(r_k)$, for $k$ outside the "pre--initial" and "post--final" intervals, form a coarse vertical line (Lemma \ref{joining 0} and Lemma \ref{thinvertical}), and the ones for $k$ in the "pre--initial" and "post--final" intervals form
also
coarse $CAT(0)$ geodesics forming with
the coarse vertical line angles $\geq 180^{\circ}$ at the endpoints. This ends the outline of
the proof of Proposition \ref{Euclidean_near_CAT(0)}.
\medskip

The following lemma treats configurations of three vertices in
a~layer. Denote the layers between $\sigma,\tau$ by $L_k$.

\begin{lemma}
\label{tripods} Suppose $p,s,t$ are three vertices in $L_k$. Then
either there exists a~vertex such that there are 1--skeleton
geodesics $ps,pt,st$ passing through this vertex or there exists
a~triangle (i.e. a 2--simplex) such that there are 1--skeleton geodesics $ps,pt,st$
passing through the edges of this triangle.
\end{lemma}

\proof
Let $p'$ be a~vertex farthest from $p$ lying both on some 1--skeleton geodesic $ps$ and some 1--skeleton geodesic $pt$.
Let $s'$ be a~vertex farthest from $s$ lying both on some 1--skeleton geodesic $sp'$ and some 1--skeleton geodesic $st$.
Finally let $t'$ be a~vertex farthest from $t$ lying both on some 1--skeleton geodesic $tp'$ and some
1--skeleton geodesic $ts'$. If two of the vertices $p',s',t'$ coincide, then all three coincide and the lemma follows immediately. Suppose now that
those three vertices are distinct.

From the choice of $p',s',t'$ it follows that any loop $\Gamma$ obtained by concatenating
some 1--skeleton geodesics $p's',s't',t'p'$ is embedded in $L_k$. Since $L_k$ is convex (Remark \ref{layers convex}), it is contractible (see remarks after Definition \ref{convex}), hence $\Gamma$ is contractible in $L_k$ (we could also invoke Lemma \ref{layers_are_infinity_large}).
Consider a~surface $T\colon D \rightarrow L_k$ of minimal area
spanned on such a~geodesic triangle $\Gamma$ (we allow the geodesics to
vary). By minimality of area the defects at interior vertices of $D$
and at interior vertices of the boundary geodesics are
non--positive. Since by Gauss--Bonnet Lemma \ref{Gauss-Bonnet} the
total sum of defects equals 6, we get that all mentioned
vertices have defects 0 and the vertices of the geodesic triangle
$D$ have defect 2. Hence $D$ is a~subcomplex of
$\mathbb{E}^2_{\Delta}$ which is a Euclidean equilateral triangle.
Denote the length of the side of this triangle by $d>0$. If $d\geq
2$ then let $u$ be the vertex such that $T(u)=p'$, let $u_1,u_2$ be
its neighbors in $D$, let $u_3$ be the common neighbor of $u_1,u_2$
in $D$ different from $u$ and let $u_4$ be the neighbor of $u_1$
different from previously mentioned vertices. By Lemma
\ref{no_trapezoid} applied to the trapezoid
$T(u)T(u_1)T(u_2)T(u_3)T(u_4)$ either we have an edge $T(u)T(u_3)$
or $T(u_2)T(u_4)$. In the first case the vertex $T(u_3)$ turns out
to lie on some 1--skeleton geodesics $sp,tp$ contradicting the
choice of $p'$. In the second case the vertex $T(u_2)$ turns out to
lie some 1--skeleton geodesics $sp,tp$, also giving a~contradiction.
Hence $d=1$ and the lemma follows. \qed
\medskip

In the next lemma we analyze the possible position of $(p_k)$ w.r.t. the
partial characteristic image $\mathcal{S}(\Delta)$ of a partial
characteristic disc $\Delta$ for $(i,j)$ for $(\sigma_k),(\tau_k)$. This
means that we assume that the layers $i\leq k\leq j$ are thick, c.f.
Definition \ref{generalized characteristic discs for incomplete
boundary}. In the language of the outline of the proof of Proposition
\ref{Euclidean_near_CAT(0)} this is the "thick" interval.
The boundary vertices of $\Delta$ are, as always, denoted
by $(v_k),(w_k)$.

For each $i\leq k\leq j$ let $s_k\in \sigma_k,t_k\in \tau_k$ be
chosen as in the previous sections to maximize the distance
$|s_kt_k|$.
Moreover, among those, choose $s_k,t_k$ to maximize the
distances $|p_ks_k|,|p_kt_k|$ (it is possible to do this
independently by Lemma
\ref{realizing_thickness_in_pairs_implies_realizing_as_a_pair}).
For
each $k$ perform in $L_k$ the construction of $s'_k,t'_k,p'_k$ as in the proof of
Lemma \ref{tripods} and denote $\overline{\chi}_k=s'_kt'_k$, which
is an edge or a~vertex in some 1--skeleton geodesic $s_kt_k$. Denote
$\chi_k=\mathcal{S}^{-1}(\overline{\chi}_k)$. Observe that $\chi_k$
does not depend on the choice of $s_k,t_k,s'_k,t'_k,p'_k$, since it is
determined by the distances $|s_kt_k|,|s_kp_k|,|t_kp_k|$. Lemmas \ref{preparation_to_triplanes} --- \ref{triplanes}
are devoted to studying the position of $\chi_k$ w.r.t. $\rho_k$ (the simplices of the Euclidean diagonal).

We refer to the path $(v_k)$ as one \emph{boundary component} of
$\Delta$, and to the path $(w_k)$ as the other boundary component.


Finally, note that in the lemma below we actually do not have to
assume that $(\sigma_k),(\tau_k)$ are directed geodesics.

\begin{lemma}
\label{preparation_to_triplanes} In the above setting, assume that
for all $i\leq k\leq j$ we have $p_k\neq p'_k$ (this does not depend on the choice of $p'_k$). Then for $i\leq k<j$,
\item (i) if $\chi_k,\chi_{k+1}$ are both edges, then they both intersect the same boundary component,
\item (ii) if one of $\chi_k,\chi_{k+1}$, say $\chi_k$, is an edge, and the second is a~vertex, then either
      $\chi_{k}, \chi_{k+1}$ span a~simplex, or they intersect the same boundary component,
\item (iii) if $\chi_k,\chi_{k+1}$ are both vertices, then they both lie on the same boundary component.
\item If we remove the assumption that $p_k\neq p'_k$, then in case (i) we only have that $\chi_k\subset S_1(\chi_{k+1})$
and $\chi_{k+1}\subset S_1(\chi_{k})$, case (ii) remains unchanged,
and in case (iii) we only have that $\chi_k,\chi_{k+1}$ span an edge.
\end{lemma}

\proof First let us prove the last assertion. We need to prove (up
to interchanging $k$ with $k+1$) that for a~vertex $u_0\in\chi_k$
either there exists a~neighbor of $u_0$ in $\chi_{k+1}$, or
$\chi_k,\chi_{k+1}$ intersect the same boundary component. Suppose
the first part of this alternative does not hold. Then, up to
interchanging $v_k$ with $w_k$, we have the following configuration (which it will take some time to describe, since we need to name all the vertices that come into play):

We have $u_0\neq w_k$, and we denote by $u_1$ the vertex following
$u_0$ on 1--skeleton~geodesic in $\Delta$ from $u_0$ to $w_k$, and
by $u_2$ the vertex following $u_1$ if $u_1\neq w_k$. In the layer $k+1$
we denote by $z_1\neq w_{k+1}$ the vertex in the residue of $u_0u_1$
and by $z_2$ the vertex following $z_1$ on 1--skeleton geodesic
$z_1w_{k+1}$. The configuration is the following: $\chi_{k+1}$ lies on the 1--skeleton geodesic
$z_2w_{k+1}$.

Fix some 1--skeleton geodesics $s_{l}\ldots s'_{l},t_{l}\ldots t'_{l},p_{l}\ldots p'_{l}$
for $l=k,k+1$. Consider a partial characteristic surface $S\colon \Delta\rightarrow X$ such
that for $l=k,k+1$ we have that $S(v_lw_l)$ (where $v_lw_l$ is the 1--skeleton geodesic in $\Delta$) contains $s_{l}\ldots s'_{l}$ and $t'_{l}\ldots t_{l}$ (this is possible by Proposition \ref{stable_geodesics_in characteristic image}).
Then $S(z_2)\in s_{k+1}\ldots s'_{k+1}\subset s_{k+1}\ldots s'_{k+1}p'_{k+1}\ldots p_{k+1}$ (where possibly $s'_{k+1}=p'_{k+1}$).
By Proposition \ref{stable_geodesics_in characteristic image} applied
to the partial characteristic surface for $p_k,p_{k+1},s_k,s_{k+1}$,
there is a neighbor of $S(z_2)$ on $s_k\ldots s'_kp'_k\ldots p_k$ (where possibly $s'_k=p'_k$).
Denote this neighbor by $\overline{x}$. Since $S(u_0)\in \overline{\chi}_k$, we have that $\overline{x}\neq S(u_1), \ \overline{x}\neq S(u_2)$. Moreover, since the vertices in the
1--skeleton geodesic $v_ku_0$ are not neighbors of $z_2$, we have by
Lemma \ref{properties of generalized characteristic discs}(iii) that
$\overline{x}\notin s_k\ldots s'_k$. So $\overline{x}\in p'_k\ldots
p_k$. But by Lemma \ref{projection lemma} the vertices
$\overline{x},S(u_1)$, together with $S(u_2)$, if defined, span
a~simplex. On the other hand, $S(u_1)$, and $S(u_2)$ if defined, lie
on the~1--skeleton geodesic $p_k\ldots p'_kt'_k\ldots t_k$ passing through $\overline{x}$.
Since $\overline{x}, S(u_1)$, and $S(u_2)$, if defined, are different vertices,
this is only possible if $\overline{x}=p'_k,
S(u_0)=s'_k, S(u_1)=t'_k$ and $u_1=w_k$, i.e.
$u_2$ is not defined. Then $\chi_k$ is an edge, $\chi_{k+1}$ is
a~vertex, and they intersect the same boundary component, which is
the second possibility of the alternative. Thus we have proved the last assertion of the lemma. In
particular, we have proved assertion (ii).
\medskip

Now we will be proving assertions (i,iii) and we may already assume
that $p_k\neq p'_k$ for $i\leq k\leq j$.

First we prove (i), by contradiction. Suppose that
$\chi_k, \chi_{k+1}$ are both edges, and
w.l.o.g. suppose that $\chi_k$ does not intersect the
boundary. This implies that $s'_k\neq s_k, t'_k\neq t_k$. Let
$\overline{z}$ be a~vertex in the projection (c.f. Definition \ref{projection}) of the triangle
$s'_kt'_kp'_k$ onto the layer $L_{k+1}$. By Lemma
\ref{projection_in_geodesic} applied thrice we get that
$\overline{z}$ lies on 1--skeleton geodesics between all pairs of
vertices from $\{s_{k+1},t_{k+1},p_{k+1}\}$, thus
$\overline{\chi}_{k+1}$ is a~vertex. Contradiction.

Now we prove (iii), by contradiction. Suppose that $\chi_k,
\chi_{k+1}$ are both vertices and one of them is non--boundary. Then in the layers $k,k+1$
of $\Delta$ there are vertices, which are common neighbors of
$\chi_k,\chi_{k+1}$, denote them by $u$ (in the layer $k$) and by $z$
(in the layer $k+1$). Moreover, either $\chi_k\neq v_k$ and $\chi_{k+1}\neq v_{k+1}$, or
$\chi_k\neq w_k$ and $\chi_{k+1}\neq w_{k+1}$. Assume w.l.o.g that the latter holds.
Consider the partial characteristic disc $\Delta_{pt}$ for
$p_k,p_{k+1},t_k,t_{k+1}$ (we are allowed to do this, since $|p_kt_k|= |p_k\overline{\chi}_k|+|\overline{\chi}_kt_k|\geq 2$ and similarly $|p_{k+1}t_{k+1}|\geq 2$)
and the corresponding partial characteristic mapping $\mathcal{S}_{pt}$. Let $x$ be
the common neighbor of
$\mathcal{S}_{pt}^{-1}(\overline{\chi}_k),\mathcal{S}^{-1}_{pt}(\overline{\chi}_{k+1})$ in $\Delta_{pt}$
lying on $\mathcal{S}_{pt}^{-1}(p_k\overline{\chi}_k)$ or $\mathcal{S}_{pt}^{-1}(p_{k+1}\overline{\chi}_{k+1})$. Assume,
w.l.o.g., that $\mathcal{S}_{pt}(x)\subset L_k$. Since vertices in $\mathcal{S}_{pt}(x),\mathcal{S}(u)\subset L_k$
are neighbors of $\overline{\chi}_{k+1}\in L_{k+1}$, we have by Lemma
\ref{projection lemma} that $\mathcal{S}_{pt}(x)$ and $\mathcal{S}(u)$ span
a simplex. On the other hand, $\overline{\chi}_k$ lies by definition
on some 1--skeleton geodesic $p_ks_k$.
By Proposition
\ref{stable_geodesics_in characteristic image}, its segments $p_k\overline{\chi}_k$ and $\overline{\chi}_ks_k$ intersect $\mathcal{S}_{pt}(x)$ and $\mathcal{S}(u)$, respectively (outside $\overline{\chi}_k$). Hence $\overline{\chi}_k$
separates vertices from
$\mathcal{S}_{pt}(x)$ and $\mathcal{S}(u)$ on a 1--skeleton geodesic $p_ks_k$.
Contradiction. Thus we have proved assertion (iii) and hence the whole lemma. \qed
\medskip

Let us introduce the following
language.

\begin{defin}
We will refer to the \emph{horizontal coordinates} of points in
various characteristic discs. Namely, we view a characteristic disc
as a~$CAT(0)$ subspace of $\mathbb{E}^2$. There we consider cartesian
coordinates such that the layers are contained in horizontal lines. We
also specify that the horizontal coordinate increases (from the left to
the right) in the direction from $v_k$ to $w_k$. We denote the
horizontal coordinate of a point $z$ by $z^x$. If $\lambda$ is a
vertical line in $\Delta$, then its horizontal coordinate is denoted
by $\lambda^x$.
\end{defin}

We will need the following technical lemma, which helps to compare
the horizontal coordinates of the preimages of vertices of $X$ in
various characteristic discs.

\begin{lemma}
\label{comparing characteristic discs} Suppose $\Delta^1,\Delta^2$
are partial characteristic discs (and $\mathcal{S}^1, \mathcal{S}^2$ resp. characteristic mappings) for the interval $(i,j)$ for some
sequences of simplices
$(\sigma_k^1),(\tau_k^1),(\sigma_k^2),(\tau_k^2)$ in the layers $L_k$ between
$\sigma,\tau$. Suppose $(p_k)_{k=i}^j,(\tilde{p}_k)_{k=i}^j$ are 1--skeleton geodesics such
that for $i\leq k\leq j$ we have that $p_k,\tilde{p}_k\in L_k$ and, for $l=1,2$, we have $p_k,\tilde{p}_k\in
\mathcal{S}^l(\Delta^l)$  but
$(\mathcal{S}^l)^{-1}(p_k)\neq (\mathcal{S}^l)^{-1}(\tilde{p}_k)$.
Then, if we vary $i\leq k\leq j$, the differences between $((\mathcal{S}^1)^{-1}(p_k))^x$ and between $((\mathcal{S}^2)^{-1}(p_k))^x$ agree.

\end{lemma}
\proof Apply Lemma \ref{properties of generalized characteristic
discs}(iii). \qed
\medskip

The following notions will help us formulate neatly the upcoming lemma.

\begin{defin}
Let $\Delta$ be a simplicial generalized characteristic disc for
$(i,j)$ such that $|v_kw_k|\geq 2$ for $i\leq k\leq j$. Let
$\chi,\rho$ be some simplices in the layer $k$ of $\Delta$, and $c\in
\mathbb{Z}_+$. We say that $\chi$ is
\item
\emph{$\partial$--left} if either $v_k\in \chi$ or $\chi$ is a
neighbor vertex of $v_k$, which has defect 1 in case $k\neq i,j$ or
defect 2 in case $k=i$ or $k=j$,
\item
\emph{$\partial$--right} if either $w_k\in \chi$ or $\chi$ is a
neighbor vertex of $w_k$, which has defect 1 in case $k\neq i,j$ or
defect 2 in case $k=i$ or $k=j$,
\item
\emph{$c$--left} from $\rho$ if $|\chi,\rho|\geq c$ and $\chi$
lies on $v_k\rho$,
\item
\emph{$c$--right} from $\rho$ if $|\chi,\rho|\geq c$ and $\chi$
lies on $\rho w_k$.
\end{defin}

In all that follows, $c$ is a positive integer. When all the pieces of
the proof of Proposition \ref{Euclidean_near_CAT(0)} are put
together, we assign $c=5$. But before this happens, we use the
variable $c$, in order to help keeping track of the role of the
constant in the various lemmas.

\begin{lemma}
\label{triplanes} Assume that for some $i<j$ and each $i\leq k\leq
j$ the layer $k$ is thick for $(\sigma_k),(\tau_k)$, and
$|p_k,\delta_k|\geq c+4$.
Then there exist $i\leq l\leq m\leq j$ such that
\item (i) for $i\leq k<l$ we have that $\chi_k$ is $\partial$--right or $c$--right from
$\rho_k$,
\item (ii) among $l\leq k\leq m$ the differences between $(\mathcal{S}^{-1}_{pr}(r_k))^x$ are $\leq c+1$,
\item (iii) for $m< k\leq j$ we have that $\chi_k$ is $\partial$--left or $c$--left from
$\rho_k$.
\item Moreover, if the maximal thickness of the layers (for $(\sigma_k),(\tau_k)$)
from $i$ to $j$ is $\geq 2c+4$ and the layers $i-1,j+1$ are thin, then
there are $l,m$ as above such that either  $m<j$ and $v_{j+1}^x-v_{m+1}^x\geq c$ (in the
characteristic disc for the thick interval $(i-1,j+1)$) or $l>i$ and
$w_{l-1}^x-w_{i-1}^x\geq c$.
\end{lemma}

The ranges for $k$ in (i),(ii),(iii), define the "initial"
subinterval, the "middle" subinterval and the "final" subinterval of a "thick" interval
discussed in the outline of the proof of Proposition
\ref{Euclidean_near_CAT(0)}. The last assertion, in the language of the outline, states that a "very thick" interval has either its "initial" or "final" subinterval non--"thin".
\medskip

\proof First we give the proof of (i)--(iii) under an additional
assumption that for all $i\leq k\leq j$ we have $p_k\neq p'_k$
(recall that this does not depend on the choice of $p'_k$). The outline of the proof with this assumption was already given at the beginning of the section.
\medskip

To start, observe that from Lemma \ref{preparation_to_triplanes} and Lemma
\ref{properties_of_special_characteristic discs}(i,ii) we get immediately the
following.
\medskip

\textbf{Corollary.} There exist
$i\leq l'\leq m'\leq j$ such that
\\ (1) for $i\leq k< l'$ the simplex $\chi_k$ is
$\partial$--right,
\\ (2) for $l'\leq k\leq m'$ the simplices $\chi_k$ are alternatingly edges and vertices and their barycenters
lie on a~straight vertical line $\lambda$ in $\Delta$; moreover for $l'<k<m'$ the simplices $\chi_k$ do not meet $v_k,w_k$,
\\ (3) for $m'< k\leq j$ the simplex $\chi_k$ is
$\partial$--left.

\medskip
Recall that the restriction to $\Delta$ of the $CAT(0)$ diagonal
$\gamma'$ (c.f. Definition \ref{euslidean diagonal}) in the characteristic disc containing $\Delta$ crosses transversally each
vertical line in $\Delta$, by Lemma \ref{transversality} (since $(j+1)-(i-1)>2$).
Let $l'\leq l\leq m'$ be maximal satisfying
$(\gamma'\cap v_kw_k)^x\leq \lambda^x-c-\frac{1}{2}$ for
$l'\leq k<l$. Similarly, let $l'\leq m\leq m'$ be minimal satisfying $(\gamma'\cap
v_kw_k)^x\geq \lambda^x+c+\frac{1}{2}$ for $m <k\leq m'$.

We prove that assertion (i) is satisfied with $l$ as above. First
consider $i\leq k< l'$. Then assertion (i) follows from assertion
(1) of the corollary. Now suppose that $l'\leq k< l$. Then, by
definitions of $l$ and $\rho_k$, if $\rho_k$ is a vertex, then
$\rho_k^x\leq \lambda^x-c-\frac{1}{2}$, and if $\rho_k$ is an edge
then the horizontal coordinates of its vertices are $\leq
\lambda^x-c$. Moreover, in case the latter inequality is an equality, we have that
$\chi_k$ is a vertex. In all cases $\chi_k$ lies to the
right of $\rho_k$ and the distance between them is $\geq c$, as
desired. Analogically, assertion (iii) holds with $m$ as above.

Now we prove assertion (ii). Consider $l\leq k\leq m$. If $l=m=l'$ or $l=m=m'$, then (ii)
follows immediately. Otherwise, by definition of $m,l$ we
have $(\gamma'\cap v_lw_l)^x>\lambda^x-c-\frac{1}{2}$
and $(\gamma'\cap v_mw_m)^x<\lambda^x+c+\frac{1}{2}$, hence
$\lambda^x-c-\frac{1}{2}<(\gamma'\cap v_kw_k)^x<\lambda^x+c+\frac{1}{2}$. By
definition of $\rho_k$, via similar considerations as in the
previous paragraph, we have that $\diam (\rho_k\cup\chi_k)\leq c+1$ and $|\rho_k,\chi_k|\leq c$.
By the former inequality we have that $p'_k$ are at distance $\leq c+1$ from $r_k$. (Record the latter one, i.e.
$|\rho_k,\chi_k|\leq c$, which we will need later in the proof.)

We would like to compute the differences between
$(\mathcal{S}^{-1}_{pr}(p'_k))^x$, when we vary $l\leq k\leq m$. These differences are equal to the differences between $(\mathcal{S}^{-1}_{ps}(p'_k))^x$ in $\Delta_{ps}$, where $\mathcal{S}_{ps}$ (resp. $\Delta_{ps}$) is the
partial characteristic mapping (resp. partial characteristic disc) for $(p_k)_{k=l}^m,(s_k)_{k=l}^m$. To see this, it is enough to apply Lemma
\ref{comparing characteristic discs} with $(p'_k),(p_k)$ in place of $(p_k),(\tilde{p}_k)$, where we use our additional assumption $p_k\neq p'_k$.

We claim that $(\mathcal{S}^{-1}_{ps}(p'_k))^x$ vary at most by $\frac{1}{2}$ for $l\leq k\leq m$. Indeed, by our additional assumption and assertion (2) of the corollary we have, for $l<k<m$, that $p_k\neq p'_k,\ s_k\neq s'_k,\ t_k\neq t'_k$. Thus we can apply Lemma \ref{preparation_to_triplanes} with $(s_k), (p_k), (t_k)$ in place of $(\sigma_k),(\tau_k),(p_k)$ to obtain, for $l\leq k\leq m$, that the barycenters of $\mathcal{S}_{ps}^{-1}(p'_ks'_k)$ lie on a common vertical line in $\Delta_{ps}$. This justifies the claim.

Thus $(\mathcal{S}^{-1}_{pr}(p'_k))^x$ vary at most by $\frac{1}{2}$, for $l\leq k\leq m$. Let $\mu$ be the greater among (at most two) values attained by $(\mathcal{S}^{-1}_{pr}(p'_k))^x$.
By the previous estimates we have that $(\mathcal{S}^{-1}_{pr}(r_k))^x\leq \mu + c+1$. On the other hand, we have $\mu\leq (\mathcal{S}^{-1}_{pr}(r_k))^x$.
Hence we obtain that the
differences between $(\mathcal{S}^{-1}_{pr}(r_k))^x$ are $\leq c+1$,
as desired.

\medskip
Now we must remove the additional assumption that for all $i\leq k\leq j$
we have $p_k\neq p_k'$. We have now only the last assertion of Lemma
\ref{preparation_to_triplanes} at our disposal.

Let $i\leq i_1\leq j_1< i_2\leq j_2< \ldots < i_q\leq j_q\leq j$, where $j_h<i_{h+1}-1$ for $1\leq h< q$,
be such that only for $i_h\leq k \leq j_h$ our
additional assumption is satisfied.
For all other $i\leq k\leq j$, in particular, for $k=i_h-1,j_h+1$ (where $1\leq h\leq q$), except possibly for $i_1-1$ if it equals $i-1$, and $j_q+1$ if it equals $j+1$, we have $|\chi_k,\rho_k|\geq |p_k,\delta_k|-1\geq c+3$. Thus for $k=i_h,j_h$, except possibly for $i_1$ if it equals $i$ and for $j_q$ if it equals $j$, we have, by Lemma \ref{properties_of_euclidean_diagonal} and by the last assertion of Lemma \ref{preparation_to_triplanes}, that  $|\chi_k,\rho_k|\geq c+1$. So for all $k$ not contained in the (open) intervals $(i_h,j_h)$ we have  $|\chi_k,\rho_k|\geq c+1$.

Put for a moment $j_0=i,\ i_{q+1}=j$. By the previous paragraph, by Lemma \ref{properties_of_euclidean_diagonal} and by the last assertion of Lemma \ref{preparation_to_triplanes}, for any $0\leq h\leq q$ and all $j_h\leq k\leq i_{h+1}$, either $\rho_k$ lies always between $\chi_k$ and $v_k$, or $\rho_k$ lies always between $\chi_k$ and $w_k$.

Now let us analyze what happens for a~fixed $1\leq h\leq q$ for $i_h\leq k\leq j_h$.
Apply our argument under the additional assumption $p_k=p'_k$ to $i=i_h,\ j=j_h$.
Observe that if $|\chi_{i_h}, \rho_{i_h}|\geq c+1$ (which holds unless possibly $h=1$ and $i_1=i$) and $\chi_{i_h}$ lies between $\rho_{i_h}$ and $v_{i_h}$, then we have that $l=m=i_h$ (otherwise we have recorded that $|\rho_k,\chi_k|\leq c$ for $l\leq k\leq m$). Similarly, if
$|\chi_{j_h}, \rho_{j_h}|\geq c+1$ (which holds unless possibly $h=q$ and $j_q=j$) and $\chi_{j_h}$ lies between $\rho_{j_h}$ and $w_{j_h}$,
 then $l=m=j_h$. In particular, those two situations cannot happen simultaneously, and if any of them happens, then either assertion (i) or assertion (iii) is valid for all $i_h\leq k\leq j_h$.

 Summarizing, there can be at most one $h$ such that $l\neq j_h$ and $m\neq i_h$. If there is no such $h$, then either assertion (i) or assertion (iii) holds for all $i\leq k\leq j$ and we are done. If not,
define $l,m$ as in the previous argument for $i=i_h,j=j_h$. They satisfy assertions (i,ii,iii), as required.

\medskip
Finally, we prove the last assertion. Pick $\lambda, l,m$ as above. Let $\Delta$ be the characteristic disc for $(i-1,j+1)$ and let $\gamma'$ be its $CAT(0)$ diagonal.
Since the maximal
thickness for $(\sigma_k),(\tau_k)$ of the layers
from $i$ to $j$ is $\geq 2c+4$, then by Lemma \ref{properties_of_special_characteristic discs}(i,ii), we have that
$v_{j+1}-w_{i-1}\geq 2c+1$.
Thus we can assume w.l.o.g. that $\lambda^x -w_{i-1}^x \geq c+\frac{1}{2}$. Thus $\lambda^x-(\gamma'\cap v_{i}w_{i})^x\geq c+\frac{1}{2}$ and $l>i$.
Observe that $\lambda$ goes through the barycenter of $\chi_l$, hence $w_{l-1}^x\geq \lambda^x-\frac{1}{2}$
so $w_{l-1}^x-w_{i-1}^x\geq c$, as desired.
\qed
\medskip

The next lemma in particular guarantees that in a "thick" interval, the vertices $\mathcal{S}_{pr}^{-1}(r_k)$ for $k$ in the "final" subinterval outside the "post--final" subinterval form a coarse vertical line. We consider it, together with the previous lemma, the heart of the proof of Proposition \ref{Euclidean_near_CAT(0)}. Below we put $\Delta$ to be the characteristic disc for the thick interval containing $i,\ldots,j$ for $(\sigma_k),(\tau_k)$. Let $v_k,w_k$ be its boundary vertices, etc.

\begin{lemma}
\label{joining 0}
Suppose that for some $i\leq j$ and for all $i\leq k\leq j$
the layer $k$ is thick for $(\sigma_k),(\tau_k)$, $|p_k,\delta_k|\geq c+2\geq 7$ and $\chi_k$ is either $\partial$--left or $c$--left from $\rho_k$. If $(\gamma'\cap v_{j+1}w_{j+1})^x=v_{j+1}^x+\frac{1}{2}$, then $v_{j+1}^x -v_i^x<c$.
\end{lemma}

\proof
By contradiction. Roughly, the idea is the following. If $v_{j+1}^x$ is relatively large w.r.t. $v_i^x$, this means that the directed geodesic $(\sigma_k)$ performs in the layers $i,\ldots,j$ an unexpected turn towards $(\tau_k)$. On the other hand, there is plenty of room in the partial characteristic disc $\Delta_{pt}$ for $(p_k),(\tau_k)$, since $p_k$ are far away from $\delta_k$, hence away from $\sigma_k$. By assumption on $\chi_k$ the corresponding characteristic image $\mathcal{S}_{pt}(\Delta_{pt})$ almost passes through $\sigma_k$. We can then see through $\Delta_{pt}$ that $(\sigma_k)$ actually goes vertically for all consecutive $i\leq k\leq j$. This yields a contradiction.
\medskip

Formally, suppose $v_{j+1}^x -v_i^x\geq c$. By increasing $i$, if necessary, we may assume that
$i$ is maximal $\leq j$ satisfying $v_{j+1}^x -v_i^x\geq c$. Hence $v_{j+1}^x -v_i^x= c$.

We claim that for
all $i\leq k\leq j$ we have that $\chi_k$ is $\partial$--left.
Indeed, by maximality of $i$ we have $(\gamma'\cap v_{j+1}w_{j+1})^x -v_k^x\leq c+\frac{1}{2}$.
By Lemma \ref{transversality} we have that $(\gamma'\cap v_kw_k)^x-(\gamma'\cap v_{j+1}w_{j+1})^x <0$.
Putting these inequalities together implies that $|v_k,\rho_k|\leq c$. Hence if $\chi_k$ is $c$--left from $\rho_k$, then it equals $v_k$, thus it is also $\partial$--left, as required. Thus we have proved the claim. Moreover,
$|v_k,\rho_k|\leq c$ together with $|p_k,\delta_k|\geq c+2$ gives also that $|p_k,\sigma_k|\geq 2$ and $p_k\neq t'_k$ for $i\leq k\leq j$.

Denote $h_k=\mathcal{S}^{-1}(t'_k)\in \chi_k$. By the claim we have $|v_kh_k|\leq 1$.
Let $\Delta_{pt}$ be the characteristic disc for the thick interval $(i_{pt},j_{pt})$ for $(p_k),(t_k)$ containing $i\leq k\leq j$ and let $\mathcal{S}_{pt}$ be the corresponding characteristic mapping (we have $|p_kt_k|=|p_kt'_k|+|t'_kt_k|\geq 2$, since $\chi_k$ is $\partial$--left). Denote $\tilde{v}_k=\mathcal{S}^{-1}_{pt}(p_k),\tilde{w}_k=\mathcal{S}^{-1}_{pt}(t_k)$.
Let $\tilde{h}_k=\mathcal{S}_{pt}^{-1}(t'_k)$. Since for $i\leq k\leq j$ we have $|t_kt'_k|\geq 1$, by Lemma \ref{comparing characteristic discs} the differences between $h^x_k$ (coordinates in $\Delta$) and $\tilde{h}^x_k$ (coordinates in $\Delta_{pt}$) agree.

Now observe that $t'_k$ spans a simplex with $\sigma_{k}$ by the claim, Lemma \ref{position of directed_geodesic_simplices_in_characteristic_image} and Lemma \ref{properties of characteristic surfaces}(iii,iv), for all $i\leq k\leq j$.
Denote $\phi=\Span\{t'_i,\sigma_i\}$. Denote by $\phi_i=\phi,\phi_{i+1},\ldots $ the simplices of the directed geodesic from $\phi$ to $\tau$. Denote by $\beta_k$ the simplices of the directed geodesic
from $t'_i$ to $\tau$. By Lemma \ref{inclusions} we have $\beta_k\subset\phi_k\supset\sigma_k$ for $k-i$ even, and
$\beta_k\supset\phi_k\subset\sigma_k$ for $k-i$ odd.
Denote by $\alpha_{k}$ the simplices of the directed geodesic in $\Delta_{pt}$ from $\tilde{h}_i$ to $\tilde{v}_{j_{pt}}\tilde{w}_{j_{pt}}$.

First we prove that for all $i\leq k\leq  j$ we have $\tilde{v}_k\notin \alpha_{k}$.
For $k=i$ this follows from $p_i\neq t'_i$. For $k>i$ we argue by contradiction.
Let $i<k_0\leq j$ be minimal such that $\tilde{v}_{k_0}\in \alpha_{k_0}$. Observe that $\Delta_{pt}$ is actually a partial characteristic disc for $(p_k),(\tau_k)$ and $(\tau_k)$ is the directed geodesic from $\tau$ to $\sigma$. Hence, similarly as in Lemma \ref{directed geodesics in characteristic discs}, for $i\leq k\leq k_0$ the simplices $\alpha_{k}$ are alternatingly
vertices and edges, with barycenters on a common vertical line. Moreover, by minimality of $k_0$, we have that $\alpha_{k_0}$ is an edge. By Lemma \ref{projection_in_geodesic} and Lemma \ref{inclusions} (applied alternatingly for consecutive layers exactly as in the proof of Proposition \ref{main_proposition_to_euclidean_geodesics}), we have that $\beta_{k}\subset \mathcal{S}_{pt}(\alpha_k)$ for $k-i$ even and $\mathcal{S}_{pt}(\alpha_k)
\subset\beta_{k}$ for $k-i$ odd, for all $i\leq k\leq k_0$. In particular, since $\alpha_i$ is a vertex and $\alpha_{k_0}$ is an edge, we have that
$p_{k_0}\in \mathcal{S}_{pt}(\alpha_{k_0})\subset \beta_{k_0}\supset \phi_{k_0} \subset \sigma_{k_0}$. But this contradicts $|p_{k_0},\sigma_{k_0}|\geq 2$. Hence we proved that for all $i\leq k\leq j$ we have $\tilde{v}_k\notin \alpha_{k}$.

From the above proof we also get that for all $i\leq k\leq j$ we have $\beta_{k}\subset \mathcal{S}_{pt}(\alpha_k)$ for $k-i$ even and $\mathcal{S}_{pt}(\alpha_k)
\subset\beta_{k}$ for $k-i$ odd,
and the simplices $\alpha_{k}$ are alternatingly vertices and edges, with barycenters on a common vertical line.
Since $t'_k$ and $\sigma_k$ span a simplex, this implies that $t'_k\in B_2(\mathcal{S}_{pt}(\alpha_k))$, hence $\tilde{h}_k\in B_2(\alpha_k)$, for $i\leq k\leq j$.
Since the barycenters of $\alpha_k$ lie on a common vertical line through $\tilde{h}_i$, we conclude that $|\tilde{h}^x_i-\tilde{h}^x_k|\leq 2\frac{1}{2}$ for $i\leq k\leq j$, in particular for $k=j$.
But $\tilde{h}^x_{j}-\tilde{h}^x_i=h^x_{j}-h^x_i\geq c-1\frac{1}{2}$.  This contradicts $c\geq 5$.
\qed
\medskip

We immediately get the following corollary, which excludes the possibility of adjoining a non--"thin" "final" subinterval of a "thick" interval to the beginning of a thin layer for $(\sigma_k),(\tau_k)$.

\begin{cor}
\label{joining}
Suppose that for some $i\leq j$ the layer $j+1$ is thin for $(\sigma_k),(\tau_k)$, and for all $i\leq k\leq j$  the layer $k$ is thick for $(\sigma_k),(\tau_k)$, $|p_k,\delta_k|\geq c+2\geq 7$ and $\chi_k$ is either $\partial$--left or $c$--left from $\rho_k$. Then $v_{j+1}^x -v_i^x<c$.
\end{cor}

The next preparatory lemma takes care of the "thin" intervals for $(\sigma_k),(\tau_k)$. Let $d$ be a positive integer.

\begin{lemma}
\label{thinvertical}
Suppose that for some $i\leq j$ the layers $i,j$
for $(\sigma_k),(\tau_k)$ have thickness $\leq d$ and for all $i\leq k\leq j$ the layer $k$ for $(\sigma_k),(\tau_k)$ has thickness $\leq 2c+3$ and $|p_k,\delta_k|\geq 2c+4$. Then the differences between $(\mathcal{S}_{pr}^{-1}(r_k))^x$ are $\leq c+2d +3\frac{1}{2}$.
\end{lemma}

We can also obtain an estimate independent of $c$ on the differences between $(\mathcal{S}_{pr}^{-1}(r_k))^x$. However, we will not need it.
\medskip

\proof We can define $p'_k$ as usual (even for thin layers). Observe that we have
$p_k\neq p'_k, \ |p_k, \sigma_k|\geq 2$, and $|p_k, \tau_k|\geq 2$, for $i\leq k\leq j$.
Let $\tilde{s}_k\in\sigma_k,\tilde{t}_k\in\tau_k$ realize maximal
distances from $p_k$ to $\sigma_k,\tau_k$, respectively.
Let $\Delta_{ps},\Delta_{pt},\mathcal{S}_{ps},\mathcal{S}_{pt}$
denote the characteristic discs and mappings for $(p_k),(\sigma_k)$
and $(p_k),(\tau_k)$, respectively, for the thick intervals
containing all $i\leq k\leq j$. Since $p_k\neq p'_k$, we have
by Lemma \ref{comparing
characteristic discs} that the differences between
$(\mathcal{S}^{-1}_{ps}(p_k))^x$, between
$(\mathcal{S}^{-1}_{pt}(p_k))^x$, and between
$(\mathcal{S}^{-1}_{pr}(p_k))^x$ agree, if we vary $k$ among $i\leq k\leq j$.

For $i\leq k\leq j$ denote $\dot{s}_k=\mathcal{S}^{-1}_{ps}(\tilde{s}_k), \dot{t}_k=\mathcal{S}^{-1}_{pt}(\tilde{t}_k)$.
Let $i\leq k_1<k_2\leq j$.
By Lemma \ref{properties_of_special_characteristic discs}(i,ii) we
have that
$\dot{s}_{k_1}^x-\dot{s}_{k_2}^x\geq
-\frac{1}{2}$ and
$\dot{t}_{k_1}^x-\dot{t}_{k_2}^x\leq
\frac{1}{2}$. In particular, $\dot{s}_{k_2}^x-\dot{s}_{j}^x\geq
-\frac{1}{2}$ and $\dot{s}_{i}^x-\dot{s}_{k_1}^x\geq
-\frac{1}{2}$, for $i\leq k_1<k_2\leq j$.
Hence
$$\dot{s}_{k_2}^x-\dot{s}_{k_1}^x\geq
\dot{s}_{j}^x-\dot{s}_{i}^x-1\geq
\dot{t}_{j}^x-\dot{t}_{i}^x-1-2d\geq -2d-1\frac{1}{2}.$$
Analogically, $$\dot{t}_{k_2}^x-\dot{t}_{k_1}^x\leq 2d+1\frac{1}{2}.$$

It will be convenient for us to assume that the
coordinates in $\Delta_{ps},\Delta_{pt}$ agree on
$\mathcal{S}^{-1}_{ps}(p_k)$ and $\mathcal{S}^{-1}_{pt}(p_k)$, so that we
can compare coordinates of points in $\Delta_{ps}$ and
$\Delta_{pt}$. With this convention, for any $i\leq k_1,k_2\leq j$
we have that
$\dot{s}_{k_1}^x-\dot{t}_{k_2}^x\geq
\dot{s}_{j}^x-\dot{t}_{j}^x-1\geq -d-1$.
Analogically $\dot{s}_{k_1}^x-\dot{t}_{k_2}^x\leq d+1$.
So altogether the differences between all the
$\dot{s}_{k}^x,\dot{t}_{k}^x$,
where $i\leq k\leq j$, are $\leq 2d+1\frac{1}{2}$. In particular, if we denote by $a$ the minimum over $k$ 
of $\dot{s}_{k}^x,\dot{t}_{k}^x$ and by $b$ the maximum over $k$ of $\dot{s}_{k}^x,\dot{t}_{k}^x$, we get $b-a\leq 2d+1\frac{1}{2}$.

For a fixed $k$, since the thickness of the layer $k$ is $\leq 2c+3$, we have that $|\tilde{s}_kr_k|\leq c+1$ or $|\tilde{t}_kr_k|\leq c+1$, hence $$\min\{|p_k\tilde{s}_k|, |p_k\tilde{t}_k|\}\leq |p_kr_k|+ c+1,$$
thus $r_k^x\geq a-(c+1)$. On the other hand,
we have $$|p_kr_k|\leq \max\{|p_k\tilde{s}_k|,|p_k\tilde{t}_k|\}+1,$$
hence $r_k^x\leq b+1$.
This altogether implies that the differences between $(\mathcal{S}^{-1}_{pr}(r_k))^x$ are
\begin{align*}
\leq (c+1)+\Big(2d+1\frac{1}{2}\Big)+1=c+2d +3\frac{1}{2}.
\end{align*}
\qed
\medskip

Finally, we prove the following easy lemma, which will be needed also later in
Section \ref{Contracting property}.
\begin{lemma}
\label{Simple_splitting} Let $\Delta$ be a generalized characteristic disc for
$(i,j)$. Let $\gamma$ be a $CAT(0)$ geodesic in $\Delta$ connecting some points
in $v_iw_i, v_jw_j$. For $i \leq k\leq j$ let $h_k\in v_kw_k$ be some points at
distance $\leq \frac{1}{2}$ from $\gamma\cap v_kw_k$. Let
$\Delta_{split}\subset \Delta$ be the generalized characteristic disc for
$(i,j)$ with $w_k$ substituted with $h_k$. Then the $CAT(0)$ geodesic $h_ih_j$ in
$\Delta_{split}$ is $1$--close to the piecewise linear boundary path $h_ih_{i+1}\ldots h_j$.
\end{lemma}
\proof
For $i\leq k\leq j$ let $h'_k$ be the points on $v_kw_k$ with $(h'_k)^x=\max \{(\gamma \cap v_kw_k)^x,h_k^x\}$.
Let $\Delta_{cut}\subset \Delta$ be the generalized characteristic disc for $(i,j)$ with $w_k$ substituted with $h'_k$. Then $\gamma$ is also a $CAT(0)$ geodesic in $\Delta_{cut}$. By Lemma \ref{geodesics stable under perturbing boubdary} applied to $\Delta_{split}\subset \Delta_{cut}$ we have that the $CAT(0)$ geodesic $h_ih_j$ in $\Delta_{split}$ is $\frac{1}{2}$--close to $\gamma$, hence $1$--close to the path $h_ih_{i+1}\ldots h_j$.
\qed
\medskip

Now we are ready to put together all pieces of the puzzle.

\medskip\par\noindent\textbf{Proof of Proposition \ref{Euclidean_near_CAT(0)}.}\ignorespaces
\ Put $c=5$. For the layers $k$ such that $|p_kr_k|\leq 7c+13$ there is
nothing to prove. Now suppose that for some $i'< j'$,
where $j'-i'\geq 2$, we have $|p_{i'}r_{i'}|=|p_{j'}r_{j'}|=7c+13$
and for $i'<k<j'$ we have $|p_kr_k|\geq 7c+14$, hence $|p_k,\delta_k|\geq
7c+13$. In particular, $p_k$ are as far from $\delta_k$ as required in Lemma \ref{joining 0} and Corollary \ref{joining}.

Let $\Delta_{pr}$ be the partial characteristic disc for $(i',j')$ for $(p_k),(r_k)$, and let $\mathcal{S}_{pr}$ be the corresponding partial characteristic mapping. Denote $u_k=\mathcal{S}^{-1}_{pr}(r_k)$.
\medskip

\textbf{Step 1.} There exist $i'\leq l\leq m \leq j'$ such that
\\ (1) for $i'\leq k< l$ the layer $k$ is thick for $(\sigma_t),(\tau_t)$, some 1--skeleton geodesic $p_ks_k$ intersects $\delta_k$, and $(\gamma'\cap v_kw_k)^x<w_k^x-\frac{1}{2}$ (in the appropriate characteristic disc for $(\sigma_t),(\tau_t)$, with the usual notation $v_k,w_k$, etc.),
\\ (2) among $l\leq k\leq m$ the differences between $u_k^x$ are $\leq 7c+11\frac{1}{2}$,
\\ (3) for $j'\geq k> m$ the layer $k$ is thick for $(\sigma_t),(\tau_t)$, some 1--skeleton geodesic $p_kt_k$ intersects $\delta_k$, and $(\gamma'\cap v_kw_k)^x>v_k^x+\frac{1}{2}$.
\medskip

This is the division into the "pre--initial" interval, the union of the central intervals, and the "post--final" interval in the language of the outline of the proof.
\medskip

Let us justify Step 1.
First consider the simple case that there are no thin layers for $(\sigma_k),(\tau_k)$ among the layers $i'\leq k\leq j'$.
Then Lemma \ref{triplanes} applied to $i=i',\ j=j'$ gives us a pair of numbers $l',m'$, which satisfies assertions (1) and (3) of Step 1 (with $l',m'$ in place of $l,m$), except for the statements on the position of $\gamma'$ (we will refer to these as \emph{incomplete} assertions (1),(3)).

Let $l<l'$ be minimal $\geq i'$ such that $(\gamma'\cap v_{l}w_{l})^x=w_{l}^x-\frac{1}{2}$ (if there is no such $l$, in particular, if $l'=i'$, then we put $l=l'$). Similarly, let
$m>m'$ be maximal $\leq j'$ such that $(\gamma'\cap v_{m}w_{m})^x=v_{m}^x+\frac{1}{2}$ (if there is no such $m$, in particular, if $m'=j'$, then we put $m=m'$). Obviously, $l,m$ satisfy complete assertions (1) and (3) of Step 1. To prove that they satisfy assertion (2), we need the following.
\medskip

\textbf{Claim.}
Among $l\leq k\leq l'-1$ the differences between $u_k^x$ are $\leq c+1$.
Analogically, among $m'+1\leq k\leq m$ the differences between $u_k^x$ are $\leq c+1$.
\medskip

To justify the claim, we need to introduce some notation. Up to the end of the proof of the claim we consider $l\leq k\leq l'-1$.
Observe that the layers $k$ for $(p_k),(s_k)$ are thick, since by incomplete assertion (1) we have that $|p_ks_k|>|p_k,\delta_k|$.
Denote by $\Delta, \mathcal{S}$ (resp. $\Delta_{ps},\mathcal{S}_{ps}$) the characteristic disc and mapping for the thick interval containing $k$ for $(\sigma_t),(\tau_t)$ (resp. for $(p_t),(s_t)$).
For each $k$ choose a~vertex $\overline{h}_k$ in $\delta_k\cap p_ks_k$ closest to $p_k$. By Proposition \ref{stable_geodesics_in characteristic image} we have that $\overline{h}_k\in \mathcal{S}_{ps}(\Delta_{ps})$.
Denote $h_k=\mathcal{S}^{-1}(\overline{h}_k),\ \tilde{h}_k=\mathcal{S}^{-1}_{ps}(\overline{h}_k)$.
Since by incomplete assertion (1) we have $s_k\neq s'_k$, Lemma \ref{comparing characteristic discs} gives that
the differences between $-h_k^x$ and between $\tilde{h}_k^x$ agree (the sign changes since $(s_k)$ plays the role of the left boundary component in $\mathcal{S}(\Delta)$ and the right one in $\mathcal{S}_{ps}(\Delta_{ps})$).
By Lemma \ref{comparing characteristic discs} applied to $\Delta_{ps}$ and $\Delta_{pr}$, and since $|p_kr_k|=|p_k\overline{h}_k|$ or
$|p_kr_k|=|p_k\overline{h}_k|+1$, we have that the differences between $u_k^x$
differ at most by 1 from the differences between $\tilde{h}_k^x$.
Hence the differences between $u_k^x$
differ at most by 1 from the differences between $-h_k^x$.

Now we can proceed with justifying the claim.
By Lemma \ref{joining 0} we have that $w^x_{l'-1}-w_{l}^x<c$, hence $(\gamma'\cap v_{l}w_{l})^x\geq w^x_{l'-1}-c$. Thus, by Lemma \ref{transversality}, we have
$(\gamma'\cap v_{k}w_{k})^x\geq w^x_{l'-1}-c$ for all $k$. This implies, by the definition of $\rho_k$, that $h^x_k\geq w^x_{l'-1}-c-\frac{1}{2}$. On the other hand, by Lemma \ref{properties_of_special_characteristic discs} we have that $w^x_k\leq w^x_{l'-1}+\frac{1}{2}$, hence we have $h^x_k\leq w^x_{l'-1}-\frac{1}{2}$. Thus the differences between $h_k^x$ are $\leq c$, hence the differences between $u_k^x$ are $\leq c+1$. This justifies the first assertion of the claim. The second one is proved analogically.
\medskip

Now we can finish the proof of Step 1 in the simple case that there are no thin layers for $(\sigma_k), (\tau_k)$, among the layers $i'\leq k\leq j'$. To prove assertion (2), we need to compare $u_{k_1}^x$ and $u_{k_2}^x$, for $l\leq k_1<k_2\leq m$. Assume, which is the worst possible case, that $l\leq k_1\leq l'-1$ and $m'+1\leq k_2\leq m$. By Lemma \ref{triplanes}(ii) and by the claim we have
\begin{align*}
|u_{k_1}^x-u_{k_2}^x|&\leq |u_{k_1}^x-u_{l'-1}^x| +\frac{1}{2}+|u_{l'}^x-u_{m'}^x| +\frac{1}{2} +|u_{m'+1}^x-u_{k_2}^x|\leq
\\
&\leq (c+1)+\frac{1}{2}+ (c+1) +\frac{1}{2}+(c+1),
\end{align*}
which is even better then the required estimate. This ends the proof of Step 1 in the simple case.
\medskip

Now consider the complex case that that there is a thin layer among the layers $i'\leq k\leq j'$. Let $(l_0,m_0)$ be a maximal (w.r.t. inclusion) interval,
with $i'\leq l_0\leq m_0\leq j'$, such that the layers $l_0,m_0$ are thin for $(\sigma_k),(\tau_k)$ and for $l_0< k< m_0$ the layer $k$ has thickness $\leq 2c+3$ (possibly $l_0=m_0$). This is the "proper thin" interval of the outline of the proof.

First we argue that for $i'\leq k<l_0$ and $m_0<k\leq j'$ the layer $k$ is thick.
Otherwise, suppose w.l.o.g. that $k_0$ is maximal $< l_0$ such that the layer $k_0$ is thin. Then, by maximality of $(l_0,m_0)$, the thick interval
$(k_0,l_0)$ contains some $k$ such that the layer $k$ has thickness $\geq 2c+4$. Thus by the last assertion of Lemma \ref{triplanes} applied to $i=k_0+1,j=l_0-1$ we get $k_0< l\leq m< l_0$ so that
either  $m<l_0-1$ and $v_{l_0}^x-v_{m+1}^x\geq c$, or $l>k_0+1$ and $w_{l-1}^x-w_{k_0}^x\geq c$. In both cases this contradicts
Corollary \ref{joining} applied respectively to $i=m+1, j=l_0-1$, or to $i=l-1, j=k_0+1$ with the roles of $v,w$ interchanged and the order on naturals inversed. Thus we have proved
that for $i'\leq k<l_0$ and $m_0<k\leq j'$ the layer $k$ is thick for $(\sigma_k),(\tau_k)$.

Now we can apply Lemma \ref{triplanes} to $i=i',\ j=l_0-1$. Denote by $l',m'$ the pair
of numbers given by its assertion. By Corollary \ref{joining} we have that $v_{l_0}^x-v_{k}^x< c$ for $m'+1\leq k\leq l_0$.
Similarly, we apply Lemma \ref{triplanes}
to $i=m_0+1,\ j=j'$ and denote by $l'',m''$ the pair of numbers given by
its assertion. By Corollary \ref{joining} we have $w_{k}^x-w_{m_0}^x< c$ for $m_0\leq k\leq l''-1$. Hence, by Lemma \ref{properties_of_special_characteristic discs}(i,ii), the thickness of the layer $k$, for $m'+1\leq k \leq l_0$ and for  $m_0\leq k\leq l'-1$, is $\leq c+1$.

Define, similarly as before, $l<l'$ to be minimal $\geq i'$ such that $(\gamma'\cap v_{l}w_{l})^x=w_{l}^x-\frac{1}{2}$ (if there is no such $l$, in particular, if $l'=i'$, then we put $l=l'$), in appropriate characteristic disc. Similarly, let
$m>m''$ be maximal $\leq j'$ such that $(\gamma'\cap v_{m}w_{m})^x=v_{m}^x+\frac{1}{2}$ (if there is no such $m$, in particular, if $m''=j'$, then we put $m=m''$).

For $l,m$ as above we have that assertion (1) follows from Lemma \ref{triplanes}(i) and assertion (3) follows
from Lemma \ref{triplanes}(iii).
As for assertion (2), assume, which is the worst possible case, that $l\leq k_1\leq l'-1$ and $m''+1\leq k_2\leq m'$. Combining Lemma \ref{thinvertical}
applied to $i=m'+1,\ j=l''-1,\ d=c+1$ with Lemma \ref{triplanes}(ii) and with the claim above (which is also valid in this complex case) we get
\begin{align*}
|u_{k_1}^x-u_{k_2}^x|&\leq |u_{k_1}^x-u_{l'-1}^x| +\frac{1}{2} +|u_{l'}^x-u_{m'}^x|+\frac{1}{2}+|u_{m'+1}^x-u_{l''-1}^x| +
\\
&+\frac{1}{2} +|u_{l''}^x-u_{m''}^x|+\frac{1}{2}+|u_{m''+1}^x-u_{k_2}^x|\leq
\\
&\leq (c+1)+\frac{1}{2}+(c+1)+\frac{1}{2}+\Big(c+2d+ 3\frac{1}{2}\Big)+
\\
&+\frac{1}{2}+(c+1)+\frac{1}{2}+(c+1)=7c+11\frac{1}{2},
\end{align*}
as required. Thus we have completed the proof of Step 1.
\medskip

\textbf{Step 2.}
$\gamma_{pr}$ is $99$--close to $(u_k)$.
\medskip

For the layers $i'\leq k<l$ define $\Delta, \mathcal{S}, \Delta_{ps},\mathcal{S}_{ps}$ and $h_k\in \rho_k\subset\Delta,\ \tilde{h}_k\in \Delta_{ps}, \ \overline{h}_k=\mathcal{S}(h_k)=\mathcal{S}_{ps}(\tilde{h}_k)$ like in Step 1 (which is possible by assertion (1) of Step 1).
Recall that the differences between $u_k^x$
differ at most by 1 from the differences between $-h_k^x$.
In particular, since for $i'\leq k_1<k_2<l$ we have $h_{k_1}^x-h_{k_2}^x\leq \frac{1}{2}$ (by Lemma \ref{transversality} and the definition of $\rho_k$), it follows that $u_{k_2}^x-u_{k_1}^x\leq 1\frac{1}{2}$. Analogically we choose vertices $h_k\in \rho_k$ (in appropriate characteristic disc) for $m<k\leq j'$, so that $|p_kr_k|=|p_k\overline{h}_k|$ or
$|p_kr_k|=|p_k\overline{h}_k|+1$. Hence for $m<k_2<k_1\leq j'$ we have $u_{k_2}^x-u_{k_1}^x\leq 1\frac{1}{2}$.
\medskip

Let $l\leq k_0\leq m$  be such that $u_{k_0}^x$
is minimal. Let $\alpha$ be a~vertical line segment in $\Delta_{pr}$ from
the layer $\max\{l-1,i'+1\}$ to the layer $\min \{m+1,j'-1\}$ at distance $2$ to the left
from $u_{k_0}$. By assertion (2) of Step 1 and by the fact that $|p_kr_k|\geq 7c+14$ this line segment
is really  contained in $\Delta_{pr}$. Let $\beta_1, \beta_2$ be $CAT(0)$
geodesics in $\Delta_{pr}$ connecting $u_{i'},u_{j'}$
to the endpoints of $\alpha$. Since $u_{k_2}^x-u_{k_1}^x\leq 1\frac{1}{2}$ for $i'\leq k_1<k_2<l$ and $m<k_2<k_1\leq j'$, we have for all $i'\leq k\leq j'$ that $u_k^x> \alpha^x$.
Hence the region in $\Delta_{pr}$ to the right of the concatenation $\beta_1\alpha\beta_2^{-1}$ is~convex, and thus contains the $CAT(0)$ geodesic in $\Delta_{pr}$ joining $u_{i'}$ with $u_{j'}$.
\medskip

We claim that $\beta_1$ is $(7c+16)$--close
to $(u_k)$. Indeed, if $l=i'$ or $l-1=i'$, then this is easy. Otherwise,
let $i'\leq k\leq l-1$. Let $\Delta''\subset \Delta|_{i'}^{l-1}$ be the generalized characteristic disc for $(i',l-1)$ obtained from $\Delta'|_{i'}^{l-1}$ (the modified characteristic disc, in which $\gamma'$ is a $CAT(0)$ geodesic) by substituting $w'_k$ with $w''_k$, such that $(w''_k)^x=h_k^x+1$. Denote $\gamma'$ restricted to the layers from $i'$ to $l-1$ by $\gamma'|_{i'}^{l-1}$. We have $\gamma'|_{i'}^{l-1}\subset \Delta''$ and by assertion (1) of Step 1 we have that $\gamma'|_{i'}^{l-1}$ is in $\Delta''$ a $CAT(0)$ geodesic. Let $\Delta'_{ps}\subset \Delta_{ps}|_{i'}^{l-1}$ be the generalized characteristic disc
for $(i',l-1)$ obtained from $\Delta_{ps}|_{i'}^{l-1}$ by deleting $\frac{1}{2}$--horizontal neighborhood of the boundary component corresponding to $(s_k)$. Observe that there is an (orientation reversing) embedding $e''\colon \Delta'' \rightarrow \Delta'_{ps}$, and that $e''(\gamma')$ is still a $CAT(0)$ geodesic in $\Delta'_{ps}$. Moreover, $e''(h_k)=\tilde{h}_k$, so that $|e''(\gamma'\cap v_kw_k)\tilde{h}_k|\leq \frac{1}{2}$.

Let $\Delta_{ph}\subset \Delta_{ps}|_{i'}^{l-1}$ be the generalized characteristic disc for $(i',l-1)$ obtained from $\Delta'_{ps}$ by splitting along $\tilde{h}_k$ (in fact $\Delta_{ph}$ is the partial characteristic disc for $(p_k),(\overline{h}_k)$).
By Lemma \ref{Simple_splitting}
the $CAT(0)$ geodesic $\tilde{h}_{i'}\tilde{h}_{l-1}$ in $\Delta_{ph}$ is $1$--close
to the boundary path $(\tilde{h}_k)$.
Now recall that there is an embedding $e \colon \Delta_{ph} \rightarrow \Delta_{pr}$, such that
 $|e(\tilde{h}_k)u_k|\leq 1$.
Let us compute the distances between the endpoints of the image under $e$ of the $CAT(0)$ geodesic
$\tilde{h}_{i'}\tilde{h}_{l-1}$ and the endpoints of $\beta_1$ in $\Delta_{pr}$.
The distance between $e(\tilde{h}_{i'})$ and $u_{i'}$ is $\leq 1$, and the distance between the second pair of endpoints is $\leq 2+(7c+11\frac{1}{2})+\frac{1}{2}$ by assertion (2) of Step 1. Hence,
by Lemma \ref{geodesics stable under perturbing boubdary}, we have that $e(\tilde{h}_{i'}\tilde{h}_{l-1})$ is $(7c+14)$--close to $\beta_1$.
Recall that $e(\tilde{h}_{i'}\tilde{h}_{l-1})$ is 1--close to $e\big((\tilde{h}_k)\big)$, which is 1--close to $(u_k)$.
Altogether,
$\beta_1$ is $\big((7c+14)+1+1\big)$--close to $(u_k)$, as desired. Thus we have justified the claim. Analogically, $\beta_2$ is $(7c+16)$--close to $(u_k)$.
\medskip

From the claim and since, by assertion (2) of Step 1, $\alpha$ is $(7c +14)$--close to $(u_k)$,
it follows that the two boundary components of the convex region in $\Delta_{pr}$ to the right of $\beta_1\alpha\beta_2^{-1}$
are $(7c+16)$--close.
Hence the $CAT(0)$ geodesic $u_{i'}u_{j'}$ in $\Delta_{pr}$ is $(7c+16)$--close
to $(u_k)$. Now consider the~$CAT(0)$ geodesic $\gamma_{pr}$ in
$\Delta_{pr}$ (which appears in the statement of the proposition)
restricted to the layers from $i'$ to $j'$. Since its endpoints are at
distance $\leq 7c+13$ from the endpoints of $u_{i'}u_{j'}$ (this is
because $|p_{i'}r_{i'}|=7c+13=|p_{j'}r_{j'}|$), we get (by Lemma \ref{geodesics stable under perturbing boubdary}) that $\gamma_{pr}$ is $(14c+29)$--close to $u_k$, as desired (recall that $c=5$). \qed

\section{Contracting}
\label{Contracting property}

In this section we prove the following consequence of Proposition
\ref{Euclidean_near_CAT(0)}, which summarizes the contracting properties of
Euclidean geodesics.

\begin{theorem}[Theorem \ref{3}]
\label{contracting} Let $s,s',t$ be vertices in a systolic complex $X$ such
that $|st|=n,|s't|=n'$. Let $(r_k)_{k=0}^n,(r'_k)_{k=0}^{n'}$ be 1--skeleton geodesics such that
$r_k\in \delta_k,r'_k\in \delta'_k$, where $(\delta_k),(\delta'_k)$ are
Euclidean geodesics for $t,s$ and for $t,s'$ respectively. Then for all $0\leq
c\leq 1$ we have $|r_{\lfloor cn\rfloor}r'_{\lfloor cn'\rfloor}|\leq
c|ss'|+C$, where $C$ is a universal constant.
\end{theorem}

In the proof we need three easy preparatory lemmas. 


\begin{lemma}
\label{splitting}
Let $D$ be a~2--dimensional systolic complex (in particular $CAT(0)$ with the standard piecewise Euclidean metric).
Let $x,y$  be
vertices in $D$. Then there exists a~1--skeleton geodesic $\omega$ in $D$ joining $x,y$
such that
if $D_0$ is the union with $\omega$ of a~connected component of $D\setminus\omega$, then the $CAT(0)$ geodesic $xy$
in $D_0$ is $1$--close to $\omega$.
\end{lemma}
\proof Let $L_i$ be the layers in $D$ between $x,y$ and let $L$ be the span in $D$
of the union of $L_i$. Observe that $L$ is convex in $CAT(0)$ sense in $D$. Hence
the $CAT(0)$ geodesic $xy$ in $D$ is contained in $L$. Now similarly as in Definition \ref{euslidean diagonal} define
vertices $\omega_i\in L_i$ to be the vertices nearest to $xy\cap L_i$ (possibly
non-unique). Analogically as in Lemma \ref{properties_of_euclidean_diagonal}
one proves that $\omega_i,\omega_{i+1}$ are neighbors, hence $(\omega_i)$ form
a path $\omega$, which is  a~1--skeleton geodesic. By the construction we have
$|\omega_i,xy\cap L_i|\leq \frac{1}{2}$ (here $|\cdot,\cdot|$ denotes the distance along the straight line). For a~fixed $D_0$ the $CAT(0)$
geodesic $xy$ in $D_0$ is contained in $L\cap D_0$, hence it is
$1$--close to $\omega$ by Lemma \ref{Simple_splitting} applied to
$L$. \qed

\begin{lemma}
\label{moving geodesic} Let $\Delta$ be a generalized characteristic disc for
$(i,j)$. Let $\Delta_{split}\subset \Delta$ be a generalized characteristic
disc for $(i,j)$ with $w_k$ substituted with $\dot{w}_k$ for some $\dot{w}_k\in
v_kw_k$. Let $\gamma, \dot{\gamma}$ be $CAT(0)$ geodesics with common endpoints
in the layers $i,j$ in $\Delta,\Delta_{split}$, respectively. Then $\dot{\gamma}\cap v_kw_k$ is not
farther from $v_k$ then $\gamma\cap v_kw_k$.
\end{lemma}
\proof Let $\Delta_0\subset \Delta$ be the characteristic disc for $(i,j)$ with
$w_k$ substituted with $\gamma\cap v_kw_k$. Then $\Delta_0\cap \Delta_{split}$
is convex in $\Delta_{split}$ and we are done. \qed

\begin{lemma}
\label{tales} Let $T$ be a~$CAT(0)$ (i.e. simply connected) subspace of $\mathbb{E}^2$, whose boundary
is an embedded loop which consists of three geodesic (in $T$) segments $\alpha,\beta,\gamma$, where $\alpha$ is contained in a~straight line in $\mathbb{E}^2$. Denote
$x=\beta\cap\gamma$. Let $\eta$ be a~geodesic in $T$ contained in a~straight
line parallel to $\alpha$ with endpoints on $\beta,\gamma$. Let $c$ denote the
ratio of the distances in $\mathbb{E}^2$ between $x$ and the line containing
$\eta$ and between $x$ and the line containing $\alpha$. Then
$\frac{|\eta|}{|\alpha|}\leq c$.
\end{lemma}

\proof
Let $y_1,y_2\in \mathbb{E}^2$ be points on the line containing $\eta$ colinear with $x$ and the endpoints of $\alpha$. By the
Tales Theorem we have $\frac{|y_1y_2|}{|\alpha|}= c$. On the other hand, since $\beta,\gamma$ are geodesics in $T$,
we get that $\eta\subset y_1y_2$.
\qed
\medskip

We are now ready for the endgame.

\medskip\par\noindent\textbf{Proof of Theorem \ref{contracting} (Theorem \ref{3}).}\ignorespaces
\
Let $m$ be maximal satisfying $r_m=r'_m$.
First assume that $\lfloor
cn\rfloor\leq m$ or $\lfloor cn'\rfloor\leq m$, say $\lfloor
cn'\rfloor\leq m$. Then $|r_{\lfloor cn'\rfloor}r'_{\lfloor cn'\rfloor}|\leq 198$. Indeed, let $\Delta$ be the characteristic disc for $(r_i),(r'_i)$ between $t$ and $r_m=r'_m$, for the thick interval containing $\lfloor cn'\rfloor$ (if layer $\lfloor cn'\rfloor$ is thin then there is nothing to prove). Then by Proposition \ref{Euclidean_near_CAT(0)} applied to $(r_i)_{i=0}^n$ and $r'_0,\ldots r'_m, r_{m+1},\ldots r_n$ we get that the $CAT(0)$ geodesic in $\Delta$ joining the barycenters of the two outermost edges is $99$--close to the boundary component corresponding to $(r_i)$. Similarly we get that this $CAT(0)$ geodesic is $99$--close to the second boundary component. Altogether we get that $|r_{\lfloor cn'\rfloor}r'_{\lfloor cn'\rfloor}|\leq 198$, as desired.
This yields
\begin{align*}
|r_{\lfloor cn\rfloor}r'_{\lfloor cn'\rfloor}|&
\leq |r_{\lfloor cn\rfloor}r_{\lfloor cn'\rfloor}|+|r_{\lfloor cn'\rfloor}r'_{\lfloor cn'\rfloor}|\leq
|\lfloor cn\rfloor-\lfloor cn'\rfloor|+ 198<\\
& < c|n-n'|+199\leq c|ss'|+199,
\end{align*}
as required. So from now on we assume that $\lfloor
cn\rfloor> m$ and $\lfloor cn'\rfloor> m$.
\medskip

Let $k$ be minimal such that $r_k$ lies on some 1--skeleton geodesic $ss'$. Now
let $k'$ be minimal such that $r'_{k'}$ lies on some 1--skeleton geodesic
$r_ks'$. Consider various 1--skeleton geodesics $\psi$ connecting $r_k$ with
$r_{k'}$. The loops $r_mr_{m+1}\ldots r_k \psi r'_{k'} r'_{k'-1}\ldots r'_m$
are embedded by the choice of $m,k,k'$. Consider a surface $S\colon
D\rightarrow X$ of minimal area spanned on such a~loop (we allow $\psi$ to
vary). By minimality of the area $D$ is systolic, hence $CAT(0)$ w.r.t. the
standard piecewise Euclidean metric. Denote the preimages of $r_i,r'_i,\psi$
in $D$ by $x_i,x'_i,\alpha$ respectively. We attach to $D$ at
$x_k,x'_{k'},x_m=x'_m$ three simplicial paths $\beta,\beta',\zeta$ of lengths
$n-k,n'-k',m$ respectively and denote obtained in this way simplical (and
$CAT(0)$) complex by $D'$. Denote the vertices in $D'\setminus D$ by
$x_n,\ldots,x_{k+1}$, by $x'_{n'},\ldots,x'_{k'+1}$, and by
$x_0=x'_0,\ldots,x_{m-1}=x'_{m-1}$ in $\beta,\beta',\zeta$ respectively.

By minimality of the area of $D$, the path $\beta \alpha \beta'^{-1}$ is a
$CAT(0)$ geodesic in $D'$. Let $D_1,D_2$ be simplicial spans in $D'$ of the
unions of all 1--skeleton geodesics from $x_0$ to $x_n$ and from $x'_0$ to
$x'_{n'}$ respectively. Observe that $D_1, D_2$ are convex (in $CAT(0)$ sense) in
$D'$, hence the $CAT(0)$ geodesics in $D'$ from $x_0$ to $x_n$ and from $x'_0$ to
$x'_{n'}$ agree with $CAT(0)$ geodesics joining those pairs of points in $D_1,
D_2$, respectively. By Proposition \ref{Euclidean_near_CAT(0)}, $(x_i)$ is
99--close (in $D_1$) to the $CAT(0)$ geodesic $x_0x_n$ and $(x'_i)$ is 99--close
(in $D_2$) to the $CAT(0)$ geodesic $x'_0x'_{n'}$.

Our goal, which
immediately implies Theorem \ref{contracting} (Theorem \ref{3}), is to get an estimate
$|x_{\lfloor cn\rfloor}x'_{\lfloor cn'\rfloor}|\leq c|x_nx'_{n'}|+C$ with
some universal constant $C$.
\medskip

We claim that for any three consecutive vertices $v,w,u$ on $\alpha$ we have
that $|x_0w|=|x_0v|+1$ implies $|x_0u|=|x_0w|+1$. We prove this claim by
contradiction. If $|x_0u|=|x_0w|-1$ then, by Lemma \ref{projection lemma},
$u,v$ are neighbors contradicting the fact that $vwu$ is a~1--skeleton geodesic. If
$|x_0u|=|x_0w|$, then by Lemma \ref{projection lemma} there exists a~vertex
$z\in D$ in the projection of the edge $wu$ onto $B_{|x_0v|}(x_0)$.
Again by Lemma \ref{projection lemma}, we have that $|zv|\leq 1$. Thus the
defect at $w$ is $\geq 1$, contradicting the minimality of the area of $D$. This
justifies the claim.

The claim implies that $\alpha$ is a~concatenation $\alpha_1\alpha_0\alpha_2$,
where vertices in $\alpha_0$ are at constant distance from $x_0$ and $\alpha_1,\alpha_2$ are
contained in 1--skeleton geodesic rays in $D'$ issuing from $x_0$. We apply Lemma
\ref{splitting} to obtain a special 1--skeleton geodesic $\omega$ in $D'$
connecting $x_0$ to $\alpha_1\cap \alpha_0$. Let $\widetilde{D}_1$ be the union
of $\omega$ and all of the components of $D'\setminus \omega$ containing some $x_i$ (i.e. on one "side" of $\omega$).
Denote by $\widetilde{D}_1^c$ the union of $\omega$ with the other components
of $D'\setminus \omega$. Denote by $\omega'$ a~1--skeleton geodesic connecting
$x_0$ to $\alpha_0\cap \alpha_2$ given by Lemma \ref{splitting} applied do
$\widetilde{D}_1^c$. Let $\widetilde{D}_2$ be the union of $\omega'$ with the
components of $\widetilde{D}_1^c\setminus \omega'$ containing some $x'_i$. Denote the
union of $\omega'$ with the other components of $\widetilde{D}_1^c\setminus
\omega'$ by $\widetilde{D}_0$.

Note that, since $\widetilde{D}_1\subset D_1, \widetilde{D}_2\subset D_2$, by
Lemma \ref{moving geodesic} we have that $(x_i)$ is 99--close to the $CAT(0)$
geodesic $x_0x_n$ in $\widetilde{D}_1$ and $(x'_i)$ is 99--close to the $CAT(0)$
geodesic $x'_0x'_{n'}$ in in $\widetilde{D}_2$. Moreover, by Lemma
\ref{splitting} and Lemma \ref{moving geodesic}, the $CAT(0)$ geodesics in
$\widetilde{D}_0,\widetilde{D}_1,\widetilde{D}_2$ joining the endpoints of
$\omega,\omega'$ are $1$--close (in particular $99$--close) to $\omega,\omega'$, respectively.
Moreover, vertices in $\alpha_0$ are at constant distance from $x_0$ in $\widetilde{D}_0$,
and $\omega\alpha_1^{-1},\omega'\alpha_2$ are 1--skeleton geodesics in
$\widetilde{D}_1, \widetilde{D_2}$,  respectively. Thus substituting
$D'=\widetilde{D}_0,\widetilde{D}_1,\widetilde{D}_2$ we have reduced the proof
of Theorem \ref{contracting} (up to replacing $C$ with $3C$) to the following two special cases:
\medskip

\noindent \textbf{(i)} vertices in $\alpha$ are at a~constant distance from $x_0$ (hence from $x_m$) or

\noindent \textbf{(ii)} $n'=k'$ and $\alpha x'_{k'}\ldots x'_0$ is a~1--skeleton geodesic.
\medskip

Observe that it is now possible that $x_i=x'_i$ for $i>m$.
Let $m'$ be maximal such that $x_{m'}=x'_{m'}$. If $\lfloor cn \rfloor\leq m'$ or $\lfloor cn' \rfloor\leq m'$, say the latter, then, since the $CAT(0)$ geodesics $x_0x_n, x'_0x'_{n'}$ in $D'$ coincide on $x_0x_{m'}$, we get that $|x_{\lfloor cn' \rfloor}x'_{\lfloor cn' \rfloor}|\leq 99+99=198$, hence
\begin{align*}
|x_{\lfloor cn \rfloor}x'_{\lfloor cn' \rfloor}|&\leq |x_{\lfloor cn \rfloor}x_{\lfloor cn' \rfloor}|+
|x_{\lfloor cn' \rfloor}x'_{\lfloor cn' \rfloor}|
\leq \\
&\leq |\lfloor cn \rfloor-\lfloor cn' \rfloor|+198
< c|n-n'|+199\leq c|x_nx'_{n'}|+199,
\end{align*}
as desired. So from now on we can assume that $\lfloor cn\rfloor> m',\lfloor
cn'\rfloor> m'$, and we can replace the component of  $D'\setminus x_{m'}$
containing $x_0$ with a simplicial path of length $m'$. Let $D$ be as before
the maximal subcomplex of $D'$ which is a topological disc.
\medskip

First suppose that we are in case (i). Observe that (up to increasing $C$ by 2) we can assume that $n=k$ and $n'=k'$.
This is because once we proved our estimate for $n=k, n'=k'$ we can concatenate an estimate
realizing path $x_{\lfloor ck\rfloor}x'_{\lfloor ck'\rfloor}$ with the paths
$x_{\lfloor ck\rfloor}\ldots x_{\lfloor cn\rfloor}$ and $x'_{\lfloor ck'\rfloor}\ldots x'_{\lfloor cn'\rfloor}$, obtaining a path from
$x_{\lfloor cn\rfloor}$ to $x'_{\lfloor cn'\rfloor}$
of length
\begin{align*}
({\lfloor cn\rfloor}&-{\lfloor ck\rfloor})+ |x_{\lfloor ck\rfloor}x'_{\lfloor ck'\rfloor}| +({\lfloor cn'\rfloor}-{\lfloor ck'\rfloor})< \\
&< (c(n-k)+1)+(c|x_kx_{k'}|+C)+(c(n'-k')+1)=\\
&=(c(|x_nx_k|)+1)+(c|x_kx_{k'}|+C)+(c(|x_{n'}x_{k'}|)+1)= c|x_nx_{n'}|+(C+2),
\end{align*}
as required.

We claim that $D$ is flat and the interior vertices of $\alpha$
have defect 0. Indeed, observe that the defects at the interior vertices of
$\alpha$ and at the interior vertices of $D$ are $\leq 0$, whereas the defect
at $x_{m'}=x'_{m'}$ is $\leq 2$. Hence, by Gauss--Bonnet Lemma
\ref{Gauss-Bonnet}, it is enough to prove that the sums of the defects at the
vertices of each of the paths $x_{m'+1}\ldots x_k$ and $x'_{m'+1}\ldots x_{k'}$ are $\leq 2$.
Suppose otherwise, w.l.o.g., that the sum of the defects at the vertices of $x_{m'+1}\ldots
x_k$ is $\geq 3$. Denote the vertex following $x_k$ on $\alpha$ by $y$. Then
$|x_{m'}y|\leq |x_{m'+1}x_k|$, hence $|x_0y|<|x_0x_k|$, which contradicts the
assumptions of case (ii). Thus we have proved the claim. In particular, $\alpha$ is contained in a straight line in $D\subset \mathbb{E}^2_{\Delta}$ and $k=k'$.

Define $\eta$ to be the~path in $D$ starting at $x_{\lfloor
ck\rfloor}$ reaching  $x'_{\lfloor ck\rfloor}$ contained (in $D\subset
\mathbb{E}^2_{\Delta}$) in a straight line parallel to $\alpha$. Let $\xi_1,\xi_2$ be $CAT(0)$ geodesics in $D$ joining $x_k$ with $x_{m'}$ and $x'_{k}$ with $x'_{m'}=x_{m'}$, respectively. Let $z_i=\eta\cap \xi_i$, for $i=1,2$. We have $|x_{\lfloor
ck\rfloor}z_1|\leq 99$ and $|z_2x'_{\lfloor
ck\rfloor}|\leq 99$ (again exceptionally $|\cdot,\cdot|$ denotes the distance along the straight line).
Let $m''$ be maximal such that $\xi_1\cap x_{m''}x'_{m''}=\xi_2\cap x_{m''}x'_{m''}$. Then for all $i\leq m''$ we have
$\xi_1\cap x_{i}x'_{i}=\xi_2\cap x_{i}x'_{i}$. In particular, if $\lfloor ck \rfloor\leq m''$, then $z_1=z_2$ and $|\eta|\leq 198$, as desired. If $\lfloor ck \rfloor> m''$, then we apply
Lemma \ref{tales} with $T\subset D$ the geodesic triangle with vertices $x_k,x'_k,\xi_1\cap x_{m''}x'_{m''}=\xi_2\cap x_{m''}x'_{m''}$. We get that $|\eta|\leq c|x_kx'_k|+198$, as desired.
\medskip

Now suppose that we are in case (ii). Like in case (i) (up to increasing $C$ by 1) we can assume that $n=k$.
Since the boundary of $D$ is a union of two geodesics, by Gauss--Bonnet Lemma
\ref{Gauss-Bonnet}, $D$ is flat. Consider an embedding $D\subset \mathbb{E}^2_{\Delta}$ such that the layers (denoted by $L_k$) between $x_{m'}=x'_{m'}$ and $x_k$ in $\mathbb{E}^2_{\Delta}$ are horizontal and $x_i$ are to the left from $x'_i$, for $i\leq k'$. By minimality of area, $\alpha$ is contained in a straight line in $D\subset \mathbb{E}^2_{\Delta}$. Like in case (i), let $\xi_1,\xi_2$ be $CAT(0)$ geodesics in $D$ joining $x_k$ with $x_{m'}$ and $x'_{k'}$ with $x'_{m'}=x_{m'}$, respectively. Similarly like in the previous case, let $m''$ be maximal such that $\xi_1\cap L_{m''}=\xi_2\cap L_{m''}$. Denote $u=\xi_1\cap L_{m''}=\xi_2\cap L_{m''}$.
By the same argument as after the choice of $m'$, we can assume that $\lfloor ck' \rfloor > m''$.
Let $z_1=\xi_1\cap L_{\lfloor ck\rfloor}, \ z_2=\xi_2\cap L_{\lfloor ck'\rfloor}$.
Let $y_1\in L_{\lfloor ck\rfloor}\cap D$ be the vertex with minimal possible $y_1^x$ but $\geq z_1^x$. Similarly, let
$y_2\in L_{\lfloor ck'\rfloor}\cap D$ be the vertex with maximal possible $y_2^x$ but $\leq z_2^x$. We claim that $|y_1y_2|=\lfloor ck\rfloor-\lfloor ck'\rfloor$.

Before we justify the claim, observe that it already implies the theorem.
Indeed, the claim gives
\begin{align*}
|x_{\lfloor ck\rfloor}x'_{\lfloor ck'\rfloor}|&\leq |x_{\lfloor ck\rfloor}y_1|+|y_1y_2|+|y_2x'_{\lfloor ck'\rfloor}|\leq \\
&\leq 99+(\lfloor ck\rfloor-\lfloor ck'\rfloor)+99<\\
&< 99+(c(k'-k)+1)+99=c|x_kx'_{k'}|+199,
\end{align*}
as desired.

Finally, let us justify the claim. We need to show that $y_2^x-y_1^x\leq \frac{\lfloor ck\rfloor-\lfloor ck'\rfloor}{2}$. By the choice of $m''$ we have that $z_1, z_2$ lie in the Euclidean triangle in $\mathbb{E}^2_{\Delta}$ with vertices $x_k, x'_{k'}, u$. Denote by $u_1$ (resp. $u_2$) the vertex on the edge $ux_k$ (resp. $ux'_{k'}$) of this triangle in $L_{\lfloor ck\rfloor}$ (resp. $L_{\lfloor ck'\rfloor}$).
Assume w.l.o.g. that $\frac{\lfloor ck' \rfloor}{k'}\geq \frac{\lfloor ck \rfloor}{k}$.
Denote then by $u_*$ the vertex on the edge $ux_k$ dividing this edge in same proportion as the proportion in which $u_2$ divides $ux'_{k'}$.
By the Tales Theorem, and since $u_1u_*\subset ux_k$ forms with the vertical direction angle $\leq 30^\circ$, we have that
\begin{align*}
u_2^x-u_1^x &\leq (u_2^x-u_*^x)+(u_*^x-u_1^x)< c\big(x_k^x-(x'_{k'})^x\big) +\frac{1}{2}=\\
&= \frac{ck-ck'}{2} + \frac{1}{2}<\frac{\lfloor ck\rfloor-\lfloor ck'\rfloor}{2} +1,
\end{align*}
hence
\begin{align*}
y_2^x-y_1^x &\leq z_2^x-z_1^x\leq u_2^x-u_1^x<\frac{\lfloor ck\rfloor-\lfloor ck'\rfloor}{2} +1.
\end{align*}
Thus, since $y_2^x-y_1^x$ and $\frac{\lfloor ck\rfloor-\lfloor ck'\rfloor}{2}$ differ by an integer (because $y_1,y_2$ are vertices in $\mathbb{E}^2_{\Delta}$), we have
$y_2^x-y_1^x\leq \frac{\lfloor ck\rfloor-\lfloor ck'\rfloor}{2}$, as desired.
This ends the proof of the claim and of the whole theorem. \qed
\medskip

If we followed the constants carefully, we would get that Theorem
\ref{contracting} (Theorem \ref{3}) is satisfied with any $C\geq 208$.

\section{Final remarks}
\label{final remarks}
In this section we state some additional results on the compactification $\overline{X}$, for which we do not provide proofs.
\medskip

$EZ$--structures explored by Farrell--Lafont \cite{FL} in relation to the Novikov conjecture concern only the torsion--free group case.
To get similar results (Novikov conjecture) for a group $G$ with torsion one needs to construct an
appropriate compactification (which we will also call an $EZ$--structure) of a \emph{classifying
space for proper $G$--actions}, denoted $\underline EG$.
$\underline EG$ is a contractible space with a proper $G$ action such that,
for every finite subgroup $F$ of $G$, the set $\underline EG^F\subset \underline EG$ of points fixed by $F$ (the \emph{fixed point set} of $F$) is contractible (in particular non--empty).

For a systolic group $G$ acting geometrically on a systolic complex
$X$ it is possible, that $X$ is an $\underline EG$.
This is not known yet however, due to the fact that it is not known whether the fixed point theorem holds for finite groups acting on systolic complexes.
The best result in this direction is the theorem of Przytycki
\cite{PP2} saying that, for every finite group $F$ acting on a systolic complex $X$, there is a non--empty
$F$--invariant subcomplex of $X$ of diameter at most $5$.
Using this, Przytycki proved \cite{PP3} that the Rips complex $R_5(X)$ of $X$ is an $\underline EG$.
Thus we can get the desired $EZ$--structure
by compactifying $R_5(X)$, using the following result analogous to \cite[Lemma 1.3]{Be} (and whose proof follows the lines of the proof of the latter).

\begin{lem}
\label{BeRips}
Let $(\cx,\partial X)$ be an $EZ$--structure on $G$ and let $G$ act geometrically on a contractible space $Y$. Then there is a natural $EZ$--structure $(Y\cup \partial X)$ on $G$.
\end{lem}

Thus in our case we get a compactification $\overline{R_5(X)}=
R_5(X)\cup \partial X$ of
$R_5(X)$ by adjoining our boundary $\partial X$ of $X$ to the Rips complex $R_5(X)$.
We claim that the following holds.

\begin{claim}
\label{Ripscpt}
Let a group $G$ act geometrically by simplicial automorphisms
on a systolic complex $X$. Let $\overline{R_5(X)}=R_5(X)\cup \partial X$ be the compactification of $R_5(X)$ obtained by applying
Lemma \ref{BeRips} to $(\cx,\partial X)$.
Then the following hold:
\begin{enumerate}
 \item $\overline{R_5(X)}$ is a Euclidean retract (ER),
 \item $\partial X$ is a $Z$--set in $\overline{R_5(X)}$,
 \item for every compact set $K\subset R_5(X)$,
$(gK)_{g\in G}$ is a null sequence,
 \item the action of $G$ on $R_5(X)$ extends to an action,
by homeomorphisms, of $G$ on $\overline{R_5(X)}$,
 \item for every finite subgroup $F$ of $G$, the fixed point set $\overline{R_5(X)}^F$ is contractible,
 \item for every finite subgroup $F$ of $G$, the fixed point set $R_5(X)^F$ is dense in $\overline{R_5(X)}^F$.
\end{enumerate}
\end{claim}

Assertions 1--4 follow from Lemma \ref{BeRips}. Assertion 6 is also easy to prove. The only difficulties in proving Claim \ref{Ripscpt} concern assertion 5. To obtain it one has to introduce good geodesics in the Rips complex
and reprove Lemma \ref{BMsys} with $\overline{R_5(X)}^F$ in place
of $\cx$.

Combining Claim \ref{Ripscpt} above and Theorem 4.1 of Rosenthal \cite{R},
we immediately get the following.

\begin{claim}
\label{NovTor}
The Novikov conjecture holds for systolic groups.
\end{claim}

Note that if the fixed point theorem holds for finite groups acting on $X$, then, by \cite{PP3}, we have that $X$ is $\underline{E}G$. Then
in Claim \ref{Ripscpt} we can substitute  $\overline{R_5(X)}$  with $\overline{X}$ and it is easier to prove assertion 5 in this case.
Then we can apply Theorem 4.1 from \cite{R} directly to $\cx$ to obtain
Claim \ref{NovTor}.
\medskip

Now we turn to the question of determining our boundary in some specific cases. We have already mentioned the case of hyperbolic systolic groups in Remark \ref{Gromovbd}. Now we consider the $CAT(0)$ case.
After making it through the second part of the article, the reader should not be surprised by the following.
\begin{claim}
\label{boudaries_equal}
If $X$ is a two--dimensional simplicial complex, which is $CAT(0)$ (which is equivalent with systolic in dimension two), then its compactification by the $CAT(0)$ visual boundary is homeomorphic in a natural way with our $\overline{X}$.
\end{claim}

For example, this implies that our boundary of a systolic Euclidean plane is a circle. The argument for Claim \ref{boudaries_equal} is that our compactification is constructed using Euclidean geodesics in systolic complexes, which in this case are coarsely $CAT(0)$ geodesics. 

The next claim concerns the following construction, which has not yet appeared in the literature. Namely Elsner and Przytycki had developed a way to turn equivariantly any $\mathcal{VH}$--complex which is $CAT(0)$ into a systolic complex (that is how they observed that the abelian product of two free groups is systolic). Although the resulting complex is usually not 2--dimensional, the only higher dimensional simplices that appear are used to deal with branching at the vertical edges. This is why we believe that the $CAT(0)$ visual boundary of the original $\mathcal{VH}$--complex is homeomorphic in a natural way with our boundary of the resulting systolic complex.

In particular, this would imply that there is a systolic group acting geometrically on two systolic complexes whose (our) boundaries are not homeomorphic. Namely, in the family of torus complexes defined by Croke--Kleiner \cite{CrK} the complexes with $\alpha=\frac{\pi}{2}$ and $\alpha=\frac{\pi}{3}$ have universal covers with non--homeomorphic $CAT(0)$ visual boundaries.
At the same time, there is a torus complex with $\alpha=\frac{\pi}{3}$, whose universal cover is 2--dimensional systolic while there also is a torus complex with $\alpha=\frac{\pi}{2}$, whose universal cover is a $\mathcal{VH}$--complex, which is $CAT(0)$.


\medskip


\begin{bibdiv}
\begin{biblist}

\bib{JS+}{article}{
title={Infinite groups with fixed point properties},
author={Arzhantseva, G.},
author={Bridson, M. R.},
author={Januszkiewicz, T.},
author={Leary, I. J.},
author={Minasyan, A.},
author={\'Swi\k{a}tkowski, J.},
status={submitted},
eprint={arXiv:0711.4238v1 [math.GR]}
}

\bib{Be}{article}{
author ={Bestvina, M.},
title  ={Local homology properties of boundaries of groups},
journal={Michigan Math. J.},
volume ={43},
date   ={1996},
number ={1},
pages  ={123--139}
}

\bib{BeMe}{article}{
author ={Bestvina, M.},
author ={Mess, G.},
title  ={The boundary of negatively curved groups},
journal={J. Amer. Math. Soc.},
volume ={4},
date   ={1991},
number ={3},
pages  ={469--481}
}

\bib{B}{article}{
author ={Bowditch, B.H.},
title  ={Cut points and canonical splittings of hyperbolic groups},
journal={Acta Math.},
volume ={180},
date   ={1998},
number ={2},
pages  ={145--186}
}


\bib{CP}{article}{
author ={Carlsson, G.},
author ={Pedersen, E.K.},
title  ={Controlled algebra and the Novikov conjectures for
$K$-- and $L$--theory},
journal={Topology},
volume ={34},
date   ={1995},
number ={3},
pages  ={731--758}
}

\bib{C}{article}{
author={Chepoi, V.},
title={Graphs of some ${\rm CAT}(0)$ complexes},
journal={Adv. in Appl. Math.},
volume={24},
date={2000},
number={2},
pages={125--179}
}

\bib{CrK}{article}{
author ={Croke, C.B.},
author ={Kleiner, B.},
title  ={Spaces with nonpositive curvature and their ideal boundaries},
journal={Topology},
volume ={39},
date   ={2000},
number ={3},
pages  ={549--556}
}

\bib{Da}{article}{
author ={Dahmani, F.},
title  ={Classifying spaces and boundaries for relatively hyperbolic groups},
journal={Topology},
volume ={(3)86},
date   ={2003},
number ={3},
pages  ={666--684}
}

\bib{D}{article}{
author ={Dranishnikov, A.N.},
title  ={On Bestvina--Mess formula},
pages  ={77--85},
book   ={   title  ={Topological and asymptotic aspects of group theory},
            series ={Contemp. Math.},
            volume ={394},
            date   ={2006},
            publisher={Amer. Math. Soc., Providence, RI}
            }
}

\bib{Du}{book}{
    title     ={Topology},
    author    ={Dugundji, J.},
    publisher ={Allyn and Bacon, Inc., Boston, Mass.},
    date      ={1966}
}

\bib{E}{article}{
    title     ={Flats and flat torus theorem in systolic spaces},
    author    ={Elsner, T.},
    status    ={submitted}
}
\bib{E3}{article}{
    title     ={Systolic spaces with isolated flats},
    author    ={Elsner, T.},
    status    ={in preparation}
}

\bib{FL}{article}{
author ={Farrell, F.T.},
author ={Lafont, J.--F.}
title  ={EZ--structures and topological applications},
journal={Comment. Math. Helv.},
volume ={80},
date   ={2005},
number ={1},
pages  ={103--121}
}


\bib{H}{article}{
    title     ={Complexes simpliciaux hyperboliques de grande dimension},
    author    ={Haglund, F.},
    journal   ={Prepublication Orsay},
    volume    ={71},
    date      ={2003},
    status    ={preprint}
}

\bib{SH}{article}{
    title     ={Separating quasi--convex subgroups in 7--systolic groups},
    author    ={Haglund, F.},
    author    ={\'Swi\k{a}tkowski, J.},
    journal   ={Groups, Geometry and Dynamics},
    volume    ={2},
    date      ={2008},
    number    ={2},
    pages     ={223--244}
}

\bib{JS}{article}{
    title     ={Simplicial Nonpositive Curvature},
    author    ={Januszkiewicz, T.},
    author    ={\'Swi\k{a}tkowski, J.},
    journal   ={Publ. Math. IHES},
    volume    ={104},
    number    ={1},
    date      ={2006},
    pages     ={1\ndash85}
}
\bib{JS2}{article}{
    title     ={Filling invariants of systolic complexes and groups},
    author    ={Januszkiewicz, T.},
    author    ={\'Swi\k{a}tkowski, J.},
    journal   ={Geometry \& Topology},
    volume    ={11},
    date      ={2007},
    pages     ={727\ndash758}
}

%

\bib{PS}{article}{
    title     ={Boundaries and JSJ decompositions of $CAT(0)$--groups},
    author    ={Papasoglu, P.},
    author    ={Swenson, E.},
    status    ={submitted},
    eprint    ={arXiv:math/0701618v1 [math.GR]}
}


\bib{PP2}{article}{
    title     ={The fixed point theorem for simplicial nonpositive curvature},
    author    ={Przytycki, P.},
    journal   ={Mathematical Proceedings of Cambridge Philosophical Society},
    volume    ={144},
    number    ={03},
    date      ={2008},
    pages     ={683--695}
}
\bib{PP3}{article}{
    title     ={\underline{E}G for systolic groups},
    author    ={Przytycki, P.},
    status    ={to appear},
    journal={Comment. Math. Helv.}
}
\bib{R}{article}{
    title     ={Split injectivity of the Baum--Connes assembly map}
    author    ={Rosenthal, D.},
    eprint    ={arXiv: math/0312047}
}
\end{biblist}
\end{bibdiv}
\end{document}